\documentclass[10pt, oneside]{article}

\usepackage{amsfonts}
\usepackage{amssymb} 
\usepackage{amsmath,amsfonts,amssymb,amsthm}
\usepackage{caption}
\usepackage{subfigure}
\usepackage{cite}
\usepackage{xcolor}
\usepackage{graphicx}
\usepackage{subfigure}
\usepackage{float}
\usepackage{paralist}
\usepackage{indentfirst}
\usepackage{authblk}
\usepackage[colorlinks,
linkcolor=red,
citecolor=blue
]{hyperref}

\usepackage[utf8]{inputenc}  % 源文件是 UTF-8
\usepackage[T1]{fontenc}
\usepackage{lmodern}         % 避免 Type 3 字体
\input{glyphtounicode}

\usepackage[T1]{fontenc}

\newtheorem{theorem}{Theorem}[section]
\newtheorem{lemma}{Lemma}[section]
\newtheorem{proposition}{Proposition}[section]

\newtheorem{remark}{Remark}[section]

\newtheorem{defn}{Definition}[section]

\def\var{\varepsilon}

\def\bma#1\ema{{\allowdisplaybreaks\begin{aligned}#1\end{aligned}}}

\numberwithin{equation}{section}

\usepackage[title,titletoc]{appendix}

\usepackage{titlesec}

\begin{document}

\allowdisplaybreaks

\title{{\textbf{
Global well-posedness and decay estimates of strong solutions to a non-conservative compressible two-fluid model  in  $L^p$ critical framework}}}

\author{ Zhiwen He$^1$, Shiyu Li$^2$, Jiayan Wu$^{3,*}$,  and Ting Zhang$^4$}

\date{}

\maketitle

\begin{abstract} 
We investigate the global existence and long-time behavior of   solutions   to a general multidimensional non-conservative compressible two-fluid model satisfying the algebraic capillary pressure relation \(P^{+}(\rho^+)-P^{-}(\rho^-)=f(\alpha^{-}\rho^{-})\) in the $L^p$ critical framework.
Our main contributions are twofold. First, under the sharp stability condition \(-\frac{s_{-}^{2}(1,1)}{\alpha^{-}(1,1)}<f^{\prime}(1)<0,\)
we establish the global existence and uniqueness of solutions for small perturbations around a linearly stable constant equilibrium.  In particular, our result removes the additional negativity assumption required in Evje--Wang--Wen [Arch. Ration. Mech. Anal. 221:1285--1316, 2016].  
Second, assuming in addition that the low-frequency of the initial perturbation is small in a lower-regularity Besov space, we derive optimal decay rates for the solution.

\end{abstract}

\renewcommand{\thefootnote}{}
%\footnote{*  Corresponding author\\ \hspace{1em}  Email addresses: zhiwen\_he@zju.edu.cn (Zhiwen He),  lishiyu0110@zstu.edu.cn (Shiyu Li), wujiayan@scut.edu.cn (Jiayan Wu), zhangting79@zju.edu.cn (Ting Zhang)}
\footnotemark\footnotetext{\textbf{2020 MSC.} 76T10; 76N10.}
\footnotetext{\textbf{Keywords.} Non-conservative compressible two-fluid model;  Small initial data; Critical regularity}

\maketitle
\section{Introduction }

\subsection{Background and Motivation }
Multifluid flows are ubiquitous in nature and arise in a wide range of industrial settings, including nuclear engineering, chemical processing, and oil and gas production. While they share some features with single-fluid models, multifluid systems involve strong interphase coupling and more intricate structures, giving rise to distinctive phenomena and presenting substantial analytical challenges in understanding their dynamics.

We consider the following generic two-phase compressible system in $\mathbb{R}^d$ with $d\geq2$:
\begin{equation}\label{model-main}
\left\{\begin{array}{l}
\partial_{t}\left(\alpha^{\pm} \rho^{\pm}\right)+\operatorname{div}\left(\alpha^{\pm} \rho^{\pm} u^{\pm}\right)=0, \\
\partial_{t}\left(\alpha^{\pm} \rho^{\pm} u^{\pm}\right)+\operatorname{div}\left(\alpha^{\pm} \rho^{\pm} u^{\pm} \otimes u^{\pm}\right)
+\alpha^{\pm} \nabla P^{\pm}\left(\rho^{\pm}\right)=\operatorname{div}\left(\alpha^{\pm} \tau^{\pm}\right)%+\kappa^\pm\alpha^\pm\rho^\pm\nabla \Delta (\alpha^\pm\rho^\pm)
,  \\
P^{+}\left(\rho^{+}\right)-P^{-}\left(\rho^{-}\right)=f(\alpha^-\rho^{-}),
\end{array}\right.
\end{equation}
where   $0\leq \alpha^{\pm}\leq 1$  are the volume
fractions of the fluid $+$ and fluid $-$, satisfying $\alpha^++\alpha^-=1$; $\rho^{\pm}(x, t) \geq 0, u^{\pm}(x, t)$ and
$P^{\pm}\left(\rho^{\pm}\right)$
denote the densities, the velocities of each phase, and the two
pressure functions, respectively.  $f(\cdot)$ is a smooth function with $f(\cdot)\in C^3([0,\infty])$, which is strictly decreasing near the equilibrium. Moreover,
$\tau^{\pm}$ are the viscous stress tensors
\begin{equation}\label{1.2}
\tau^{\pm}\triangleq\mu^{\pm}\left(\nabla u^{\pm}+\nabla^{t}
u^{\pm}\right)+\lambda^{\pm} \operatorname{div} u^{\pm} \mathrm{Id},
\end{equation}
where the constants $\mu^{\pm}$ and $\lambda^{\pm}$ are shear and
bulk viscosity coefficients satisfying the physical conditions:
\begin{align}\label{system2-2}
\mu^{\pm}>0 \mbox{ and } 2 \mu^{\pm}+d \lambda^{\pm} \geq 0,
\end{align}
which
implies that $\mu^{\pm}+\lambda^{\pm}>0.$ %$\kappa^\pm $ are the capillary coefficients. 
The pressure functions $P^{\pm}$ are assumed to satisfy the general conditions
\begin{align}
P^{\pm}(s)\in C^{\infty}(\mathbb{R}_{+}),\quad \partial_{s} P^{\pm}(s)>0\quad\text{for}\quad s>0.\label{system3}
\end{align}
In this paper, we are concerned with the Cauchy problem of the system \eqref{model-main} subject to the initial data
\begin{align}
	&(\alpha^{\pm},\rho^{\pm}, u^{\pm})(x,0)=(\alpha^{\pm}_0,\rho^{\pm}_0, u^{\pm}_0)(x)\rightarrow (\bar{\alpha}^{\pm},\bar{\rho}^{\pm},0),\quad |x|\rightarrow \infty.\label{system4}
\end{align}
Here, $\bar{\alpha}^{\pm}$ and $\bar{\rho}^{\pm}$ are positive constants.

This system \eqref{model-main} is commonly referred to as a two-fluid flow model with an algebraic closure; see, for instance, \cite{Brennenbook,BDGG2010,BHL2012,Prospererttibook,Ishiibook2,WYZ2018} and the references therein for further background. In the seminal work \cite{BDGG2010},  
Bresch–Desjardins–Ghidaglia–Grenier  investigated a model closely related to \eqref{model-main}. More
precisely, they considered the setting where (i) the two phases share a common pressure, $P^+=P^-=P$; (ii) the viscosities depend linearly on the densities, $\mu(\rho^\pm)=\mu^\pm\rho^\pm$ and $\lambda^\pm=0;$ and (iii)  capillarity contributions (Korteweg-type terms) involving third-order derivatives of $\alpha^\pm\rho^\pm$ are included. An informative overview of this model, emphasizing its key features and underlying mechanisms, can be found in the introduction of \cite{BDGG2010}.  From the analytical viewpoint, the two-fluid system presents several difficulties, including:
\begin{itemize}
\item \textbf{Partial dissipation.} The mass conservation equations lack dissipation, whereas the momentum equations contain viscous dissipation.

\item \textbf{Degeneracy to single-phase regions.} Transitions to regions where $\alpha^{+} \rho^{+}$ or $\alpha^{-} \rho^{-}$ vanishes may occur when either the volume fractions $\alpha^{\pm}$ or the densities
$\rho^{\pm}$ become zero.

\item \textbf{Non-conservative structure.} The non-conservative terms $\alpha^{\pm} \nabla P^{\pm}$ in the momentum equations prevent a straightforward adaptation of methods developed for single-phase models.
\end{itemize}

In the common-pressure case, Bresch-Desjardins-Ghidaglia-Grenier \cite{BDGG2010} proved the existence of global weak solutions
on the periodic domain for  $1<\gamma^\pm <6$. Later, 
Bresch-Huang-Li  \cite{BHL2012} showed the existence
of global weak solutions in one dimension, when the capillary effects are absent   and
 $\gamma^\pm >1$. Under the hypotheses that the  initial fraction densities$ \alpha^\pm\rho^\pm-1 $ and $u^\pm$ are small in $H^{l}\cap L^1$ for some integer $l\geq 3,$   Cui–Wang–Yao–Zhu \cite{CWYZ2016} derived  the global stability  of classical solutions to the three-dimensional Cauchy problem with capillarity effects  under the special density–dependent viscosities and equal coefficients
\begin{align}\label{1.6}
   \mu^\pm(\rho^{\pm})=\nu\rho^\pm\quad\text{with}\quad \nu>0, \quad \lambda^{\pm}(\rho^{\pm})=0,\quad \kappa^\pm=\kappa>0.
\end{align}
Li–Wang–Wu–Zhang \cite{LWWZ2023} recently proved the global stability of model \eqref{model-main} for general constant viscosities as in \eqref{1.2} and unequal capillary coefficients, by exploiting the dissipative structure of the system and making full use of several key observations. Meanwhile, it should be emphasized that both  \cite{CWYZ2016} and \cite{LWWZ2023} hinge essentially on a low-frequency condition, i.e., the smallness of the 
$L^1$-norm of the initial perturbation.  
Subsequently, Wu-Yao-Zhang \cite{WYZMathAnn} obtained the global stability for the 3D Cauchy problem with sufficiently small initial data in
$L^1\cap H^3$, and  
very recently,  \cite{TWYZ2025} established   the vanishing capillary limit based on that result. Shou-Wu-Yao-Zhang \cite{SWYZ}  further removed the $L^1$ condition in  \cite{WYZMathAnn} and established the global existence and large-time behavior of solutions for data close to equilibrium in the low-regularity space $\dot{B}^{-1/2}_{2,1}\cap\dot{B}^{5/2}_{2,1}$. In a related recent work, Wang--Wen \cite{WangWen2026} established global
solvability for initial perturbations that are sufficiently small in
$H^3(\mathbb{R}^3)$, with
$(P_0-\bar P,u_0^+,u_0^-)\in L^p(\mathbb{R}^3)$ for
$1\leq p<\frac65$, without requiring the corresponding $L^1$-norm to be small. Xu–Chi \cite{Xu-Nonlinearity} studied the Cauchy problem for a compressible generic two-fluid model with capillarity effects in any dimension $d\geq 2$, and proved the global existence and uniqueness of solutions in critical Besov spaces. Moreover, under a mild additional decay assumption involving only the low-frequency part of the initial data, they established the optimal time-decay rates for the resulting global solutions.  

When considering \eqref{model-main} with $(f\neq 0)$, capillary pressure forces are incorporated; such effects are standard in the modeling of two-phase flow in porous media. For further background on this model and its variants, we refer to \cite{Bearbook,Ishiibook,Prospererttibook}  and the references therein.
In the excellent work \cite{EvjeWangWang2016}, Evje--Wang--Wen investigated the
two--fluid model \eqref{model-main}. As a matter of
fact, they made the following assumptions on pressures:
\begin{equation}\label{uneq-pressure}
P^+(\rho^+)-P^-(\rho^-)=(\rho^+)^{\bar{\gamma}^+}-(\rho^-)^{\bar{\gamma}^-}=f(\alpha_-\rho_-),
\end{equation}
where $f$ is the so-called capillary pressure function which belongs to $C^3([0,
\infty))$, and is a strictly decreasing function near the
equilibrium satisfying
\begin{equation}-\frac{s_{-}^{2}(1,1)}{\alpha^{-}(1,1)}<f^{\prime}(1)<\frac{\bar{\eta}-s_{-}^{2}(1,1)}{\alpha^{-}(1,1)}<0,\label{condition-f}
\end{equation}
where $\bar{\eta}$ is a positive, small fixed constant, and
$s_{\pm}^{2}\triangleq\frac{\mathrm{d} P^{\pm}}{\mathrm{d}
\rho^{\pm}}\left(\rho^{\pm}\right)=\bar{\gamma}^{\pm}
\frac{P^{\pm}\left(\rho^{\pm}\right)}{\rho^{\pm}}$ represent the
sound speed of each phase respectively. Under the assumptions
\eqref{uneq-pressure} and \eqref{condition-f} on pressures, they obtained global
existence and decay rates of the solutions when the initial
perturbation is sufficiently small in $H^3$. 
It should be mentioned that the technical smallness of $\bar{\eta}$ is essential for closing the higher–order energy estimates;  
see \cite{EvjeWangWang2016}, pp. 1302–1308 for details. Building on this, Wang--Wen--Yao \cite{WWY} further developed the global well-posedness theory and demonstrated exponential stability for the corresponding initial--boundary value problem. 
These works revealed that the assumptions regarding the capillary pressure $f$ exert a decisive stabilizing effect on the dynamics. In addition, based on  \cite{EvjeWangWang2016}, Lai-Wen-Yao \cite{LWY2017} established the vanishing–capillarity limit for the non-conservative compressible two–fluid model \eqref{model-main} with capillarity effects.
In the opposite direction, Wu-Yao-Zhang \cite{WuYaoZhangin} proved instability when $f^{\prime}(1)>0$. Wu-Wu-Zhang \cite{WuWuZhang-S} improved condition \eqref{condition-f} to  \begin{align}
    -\frac{s_{-}^{2}(1,1)}{\alpha^{-}(1,1)}<f^{\prime}(1)<0,\label{cond-f2}
\end{align} and on this basis established the global existence and optimal large-time behavior for a generic compressible two-fluid model \eqref{model-main}.  Meanwhile, Shou-Wu-Yao-Zhang \cite{SWYZ2} studied the global existence and long-time behavior of large solutions in the high-capillarity regime for a general multidimensional non-conservative compressible two-fluid model \eqref{model-main} with capillarity effects  under condition \eqref{1.6} and  \eqref{cond-f2}. Here we emphasize that conditions in Refs. \cite{WuWuZhang-S} and \cite{SWYZ2} were improved from condition \eqref{condition-f} to condition \eqref{cond-f2} by means of different approaches.

For the compressible Navier–Stokes equations, the global-in-time existence of strong solutions (for small perturbations around equilibrium) was established in \cite{MR1779621,CPAM2010,MR2679372}. Later, Danchin–Xu \cite{Danchin-Xu-ARMA2017} derived optimal large-time decay estimates  under a mild additional decay assumption on low-frequency of the initial data. Subsequently, Xu \cite{CMP2019Xu} further relaxed the low-frequency assumption by enlarging the admissible range of indices, thereby extending the applicability of the decay theory. Recently,  Brandolese-Shou-Xu-Zhang   \cite{MR4792291} established a sharp decay characterization for the
compressible Navier-Stokes equations in the critical $L^p$ framework.
In this paper,  under the minimal regularity assumption of the initial data, we establish the global well-posedness of solutions for  \eqref{model-main}. Observe that the system \eqref{model-main} is (formally) invariant, for any $l>0$, under the scaling
\begin{align*}
\alpha^\pm (t,x)\rightsquigarrow \alpha^\pm(l^{2}t,lx),\quad  \rho^\pm(t,x)\rightsquigarrow \rho^\pm(l^{2}t,lx) \ \ \hbox{and} \ \ u^\pm(t,x)\rightsquigarrow l u^\pm(l^{2}t,lx)
\end{align*}
up to a change of the pressure term $P^\pm $ into $ l^{2}P^\pm$. This motivates the use of so-called \emph{critical} functional spaces, whose norms are invariant under the above scaling. 

The main purpose of the present paper is to develop a global well-posedness and decay theory for the non-conservative compressible two-fluid system \eqref{model-main} in the $L^{p}$-critical framework. Our contributions are twofold. First, compared with the existing $L^{2}$-based theories for complessible Navier--Stokes equation or two-fluid models, our approach allows the velocity component to have negative regularity when $p>d$, and hence covers highly oscillatory initial velocities. Compared with the classical energy framework of Evje--Wang--Wen \cite{EvjeWangWang2016}, we work under the sharp stability condition 
$$-\frac{s_-^2(1,1)}{\alpha^-(1,1)}<f'(1)<0,$$
without imposing the additional technical smallness assumption in their stability condition. The main  difficulty is that the two density equations are transported by two different velocities, while the momentum equations are coupled through non-conservative pressure terms. Therefore, neither a single Lagrangian transformation nor two independent Lagrangian transformations seem sufficient to close the critical estimates. To overcome this difficulty, we introduce two effective velocities (see \eqref{TEvelo}) adapted to the coupled linearized system and exploit the algebraic stability structure encoded in
$$\gamma_1\gamma_4-\gamma_2\gamma_3>0.$$
This allows us to recover the missing dissipation of the density variables and to close the global estimates in critical Besov spaces.

Second, under an additional low-frequency assumption on the initial perturbation in a lower-regularity Besov space $\dot{B}^{\sigma_1}_{2,\infty}$  $(\sigma_0\leq\sigma_1<\frac{d}{2}-1 \ \text{with}\ \sigma_0=\frac{d}{2}-\frac{2d}{p})$, we derive optimal time-decay estimates for the global solutions. This gives a complete description of the large-time behavior of small solutions in the $L^p$-critical framework and extends the decay theory for compressible Navier--Stokes equations to the present non-conservative two-fluid system with algebraic capillary pressure closure.

%It should be stressed that the existing results for the system \eqref{model-main}—including the global existence and large-time decay of strong solutions—are largely established in an $L^{2}$-based setting, most notably in critical Besov spaces built on  $L^{2}$. In contrast, to the best of our knowledge, there are no results that provide both global well-posedness and time-decay estimates for \eqref{model-main} in the $L^{p}$ framework. This observation motivates the present work and   the goal of the present paper is twofold:
%\begin{itemize}
%\item \textbf{Global well-posedness in critical spaces.} We prove the global existence and uniqueness of solutions to \eqref{model-main} for small initial data around the equilibrium state in   the $L^{p}$ framework, under the stability condition 
%$$
%-\frac{s_{-}^{2}(1,1)}{\alpha^{-}(1,1)}<f^{\prime}(1)<0,
%$$
%which is weaker than \eqref{condition-f} in \cite{EvjeWangWang2016}.
%\item \textbf{Optimal large-time behavior.} Under the additional smallness assumption that the low-frequency of  initial perturbation lies in a lower–regularity Besov space $\dot{B}^{\sigma_1}_{2,\infty}$  $(\sigma_0\leq\sigma_1<\frac{d}{2}-1 \ \text{with}\ \sigma_0=\frac{d}{2}-\frac{2d}{p})$,  we derive optimal decay rates of the global solution to  \eqref{model-main}. 
%\end{itemize}

\subsection{Reformulation}
In this subsection, we devote ourselves to reformulating the system
\eqref{model-main}. It follows from the pressure relation  \eqref{model-main}$_3$ that
\begin{equation}\label{1.9}
\mathrm{d} P^+-\mathrm{d} P^-= \mathrm{d} f(\alpha^-\rho^-)
\end{equation}
and 
\begin{equation} 
\mathrm{d} P^+ = s_{+}^{2}\mathrm{d}\rho^+,\  \mathrm{d} P^- = s_{-}^{2}\mathrm{d}\rho^-,
\end{equation}
 where $s_{\pm}^{2}\triangleq\frac{\mathrm{d} P^{\pm}}{\mathrm{d}
\rho^{\pm}}\left(\rho^{\pm}\right)$ represent the
sound speed of each phase respectively. As in \cite{BDGG2010,EvjeWangWang2016}, we
introduce the fractional densities
\begin{equation}\label{1.10}
R^{\pm}=\alpha^{\pm} \rho^{\pm},
\end{equation}
which, together with the relation: $\alpha^++\alpha^-=1$, yields
\begin{equation}\label{1.11}
\mathrm{d} \rho^{+}=\frac{1}{\alpha_{+}}\left(\mathrm{d}
R^{+}-\rho^{+} \mathrm{d} \alpha^{+}\right), \quad \mathrm{d}
\rho^{-}=\frac{1}{\alpha_{-}}\left(\mathrm{d} R^{-}+\rho^{-}
\mathrm{d} \alpha^{+}\right).
\end{equation}
Combining \eqref{1.9} and \eqref{1.10}, we  obtain
\begin{equation}\label{1.12}
\mathrm{d} \alpha^{+}=\frac{\alpha^{-} s_{+}^{2}}{\alpha^{-}
\rho^{+} s_{+}^{2}+\alpha^{+} \rho^{-} s_{-}^{2}} \mathrm{d}
R^{+}-\frac{\alpha^{+}\alpha^- }{\alpha^{-} \rho^{+}
s_{+}^{2}+\alpha^{+} \rho^{-}
s_{-}^{2}}\left(\frac{s_{-}^{2}}{\alpha^-}+f^{\prime}\right)
\mathrm{d} R^{-}. \end{equation}
 Substituting \eqref{1.12} into \eqref{1.11}, we deduce
the following identities:
$$
\mathrm{d} \rho^{+}=\frac{\rho^{+} \rho^{-}
s_{-}^{2}}{R^{-}\left(\rho^{+}\right)^{2}
s_{+}^{2}+R^{+}\left(\rho^{-}\right)^{2} s_{-}^{2}}\left(\rho^{-}
\mathrm{d} R^{+}+\left(\rho^++\rho^+\frac{\alpha^-f^{\prime}}{s_-^2}\right)\mathrm{d} R^{-}\right),
$$
and
$$
\mathrm{d} \rho^{-}=\frac{\rho^{+} \rho^{-}
s_{+}^{2}}{R^{-}\left(\rho^{+}\right)^{2}
s_{+}^{2}+R^{+}\left(\rho^{-}\right)^{2} s_{-}^{2}}\left(\rho^{-}
\mathrm{d} R^{+}+\left(\rho^+-\rho^-\frac{\alpha^+f^{\prime}}{s_+^2}\right)\mathrm{d} R^{-}\right),
$$
which, together with \eqref{1.9}, yields the common pressure
differentials 
$$
\mathrm{d} P^+=\mathcal{C}^2\left(\rho^{-} \mathrm{d}
R^{+}+\left(\rho^++\rho^+\frac{\alpha^-f^{\prime}}{s_-^2}\right) \mathrm{d} R^{-}\right),\quad \mathrm{d} P^-=\mathcal{C}^2\left(\rho^{-} \mathrm{d}
R^{+}+\left(\rho^+-\rho^-\frac{\alpha^+f^{\prime}}{s_+^2}\right) \mathrm{d} R^{-}\right), 
$$
where
$$
\mathcal{C}^2\triangleq\frac{s_{-}^{2} s_{+}^{2}}{\alpha^{-} \rho^{+}
s_{+}^{2}+\alpha^{+} \rho^{-} s_{-}^{2}}.
$$
Next, since  $\alpha^++\alpha^-=1$, we have  
\begin{align} \label{1.15}
\frac{R^+}{\rho^+} +\frac{R^-}{\rho^-} =1\quad\text{and hence}\quad \rho^- \;=\;\frac{R^-\,\rho^+}{\rho^+ - R^+}\,. 
\end{align}
From the pressure relation \eqref{model-main}$_3$, define
$$
\varphi(\rho^+) \;=\; P^+(\rho^+)\;-\;P^-\!\Bigl(\frac{R^-\,\rho^+}{\rho^+ - R^+}\Bigr)\;-\;f(R^-),
$$
and determine $\rho^+$ from $\varphi(\rho^+)=0.$

Fixing $R^\pm$ and  differentiating with respect to $\rho^+$ gives
$$
\varphi'(\rho^+)
=\;(P^+)'(\rho^+)
\;+\;(P^-)' \!\Bigl(\tfrac{R^-\,\rho^+}{\rho^+ - R^+}\Bigr)\;\frac{R^-\,R^+}{(\rho^+ - R^+)^2}
\;=\;
s_+^2
\;+\;
s_-^2\,\frac{R^-\,R^+}{(\rho^+ - R^+)^2},
$$
where $s_\pm^2=(P^\pm)'(\cdot)$ denotes the sound speed of each phase separately.  By construction $R^+<\rho^+<\infty$, it is clear that $\varphi'(\rho^+)>0$ on $(R^+,\infty)$ while $\varphi$ maps $(R^+,\infty)$ onto $(-\infty,\infty)$.  By the implicit function theorem, there is a unique solution
$$
\rho^+=\rho^+(R^+,R^-)\in (R^+,\infty).
$$
Finally, combining this with $\alpha^++\alpha^-=1$, \eqref{1.11}  and \eqref{1.15}, we set
$$
\rho^\pm=\rho^\pm(R^+,R^-), 
\quad
\alpha^\pm=\alpha^\pm(R^+,R^-),
\quad
\mathcal{C}=\mathcal{C}(R^+,R^-),
$$
%where the mappings $\phi_\pm$, $\psi_\pm$ and $\varsigma$ are given explicitly   on
see p. 614 of \cite{BDGG2010} for more details.

Therefore, we can rewrite  \eqref{model-main}  in the following equivalent form
\begin{equation}\label{model-main1}
\left\{\begin{array}{l}
\partial_{t} R^{+}+  u^{+}\cdot\nabla R^+ =-R^+\operatorname{div} u^{+}, \\
R^{+}(\partial_{t} u^{+} +  u^{+} \cdot\nabla u^{+} )+\alpha^{+} \mathcal{C}^2\left[ \rho^-\nabla R^++\left(\rho^++\rho^+\frac{\alpha^-f^{\prime}}{s_-^2}\right)\nabla R^-\right] \\
\hspace{2.5cm}=\operatorname{div}\left\{\alpha^{+}\left[\mu^{+}\left(\nabla u^{+}+\nabla^{t} u^{+}\right)
+\lambda^{+} \operatorname{div} u^{+} \operatorname{Id}\right]\right\}, \\
\partial_{t} R^{-}+u^{-}\cdot\nabla R^- =-R^-\operatorname{div} u^{-},\\
R^{-}(\partial_{t} u^{-} +  u^{-} \cdot\nabla u^{-} )+\alpha^{-} \mathcal{C}^2\left[ \rho^-\nabla R^++\left(\rho^+-\rho^-\frac{\alpha^+f^{\prime}}{s_+^2}\right)\nabla R^-\right] \\
\hspace{2.5cm}=\operatorname{div}\left\{\alpha^{-}\left[\mu^{-}\left(\nabla
u^{-}+\nabla^{t} u^{-}\right)+\lambda^{-} \operatorname{div} u^{-}
\operatorname{Id}\right]\right\}
\end{array}\right.
\end{equation}
with the initial data
\begin{align}\label{reformulat2}
	&(R^{+}, R^{-}, u^{+},u^{-})(x,0)=(R^{+}_{0}, R^{-}_{0}, u^{+}_{0},u^{-}_{0})(x)\rightarrow (\bar{R}^+,\bar{R}^-,0,0),\quad |x|\rightarrow \infty,
\end{align}
where $\bar{R}^\pm:=\bar{\alpha}^\pm \bar{\rho}^\pm$. Without loss of generality, in what follows, we shall assume
$
\bar{\alpha}^\pm \bar{\rho}^\pm=  1.
$
To facilitate the proof of Theorem  \ref{Thm1}, we need some reformulations of \eqref{model-main1}. More precisely, taking a change of variables by 
$$n^\pm=R^\pm -1,$$
the system \eqref{model-main1} can be reformulated as
\begin{equation}  \label{equ:fornu-pm-S}
\left\{\begin{array}{l}
\partial_{t} n^{+}+\operatorname{div} u^{+}=S_1^+, \\
\partial_{t} u^{+}+\beta_1 \nabla n^{+}+\beta_2 \nabla n^{-}-\nu_{1}^{+} \Delta u^{+}-
\nu_{2}^{+} \nabla \operatorname{div} u^{+}=S_2^+, \\
\partial_{t} n^{-}+\operatorname{div} u^{-}=S_1^-, \\
\partial_{t} u^{-}+\beta_3 \nabla n^{+}+\beta_4 \nabla n^{-}-\nu_{1}^{-} \Delta u^{-}-
\nu_{2}^{-} \nabla \operatorname{div} u^{-}=S_2^-,
\end{array}\right.
\end{equation}
where
%$\beta_1=\frac{\mathcal{C}^{2}(1,1)\rho^-(1,1)}{\rho^+(1,1)}$, $\beta_2=\mathcal{C}^{2}(1,1)+\frac{\mathcal{C}^{2}(1,1)\alpha^-(1,1)f'(1)}{s_-^2(1,1)},$ $\beta_3=\mathcal{C}^{2}(1,1),$ $\beta_4=\frac{\mathcal{C}^{2}(1,1)\rho^+(1,1)}{\rho^-(1,1)}-\frac{\mathcal{C}^{2}(1,1)\alpha^+(1,1)f'(1)}{s_+^2(1,1)}$, $\nu_1^{\pm}=\frac{\mu^{\pm}}{\rho^{\pm}(1,1)},$ $\nu_2^{\pm}=\frac{\mu^{\pm}+\lambda^{\pm}}{\rho^{\pm}(1,1)}>0$ 
\begin{equation}\label{1.43}
	\left\{
	\begin{aligned}
    &\beta_1=\frac{\mathcal{C}^{2}(1,1)\rho^-(1,1)}{\rho^+(1,1)},\quad \quad \beta_2=\mathcal{C}^{2}(1,1)+\frac{\mathcal{C}^{2}(1,1)\alpha^-(1,1)f'(1)}{s_-^2(1,1)},\\
    & \beta_3=\mathcal{C}^{2}(1,1),\quad\quad \quad\quad\quad\quad\beta_4=\frac{\mathcal{C}^{2}(1,1)\rho^+(1,1)}{\rho^-(1,1)}-\frac{\mathcal{C}^{2}(1,1)\alpha^+(1,1)f'(1)}{s_+^2(1,1)}
	\end{aligned}
	\right.
\end{equation}
and $\nu_1^{\pm}=\frac{\mu^{\pm}}{\rho^{\pm}(1,1)},$ $\nu_2^{\pm}=\frac{\mu^{\pm}+\lambda^{\pm}}{\rho^{\pm}(1,1)}>0$
and the source terms are 
\begin{equation}
    \left\{
\begin{aligned}
    S_1^{+}=&-n^+\text{div}u^+-u^+\cdot\nabla n^+,\\
    S_2^+=&-g^+(n^+,n^-)\nabla n^+ -\tilde{g}^+(n^+,n^-)\nabla n^- -(u^+\cdot \nabla )u^+\\
    &+\mu^+ h^+(n^+,n^-)\nabla n^+\nabla u^+ +\mu^+k^+(n^+,n^-)\nabla n^-\nabla u^+\\
    &+\mu^+ h^+(n^+,n^-)\nabla n^+\nabla^t u^+ +\mu^+k^+(n^+,n^-)\nabla n^-\nabla^t u^+\\
    &+\lambda^+h^+(n^+,n^-)\nabla n^+\text{div}u^+ Id +\lambda^+k^+(n^+,n^-)\nabla n^-\text{div}u^+Id\\
    &+\mu^+l^+(n^+,n^-)\Delta u^++(\mu^++\lambda^+)l^+(n^+,n^-)\nabla\text{div}u^+,\\
    S_1^-=&-n^-\text{div}u^--u^-\cdot\nabla n^-,\\
    S_2^-=&-g^-(n^+,n^-)\nabla n^- -\tilde{g}^-(n^+,n^-)\nabla n^+ -(u^-\cdot \nabla )u^-\\
    &+\mu^- h^-(n^+,n^-)\nabla n^+\nabla u^- +\mu^-k^-(n^+,n^-)\nabla n^-\nabla u^-\\
    &+\mu^- h^-(n^+,n^-)\nabla n^+\nabla^t u^- +\mu^-k^-(n^+,n^-)\nabla n^-\nabla^t u^-\\
    &+\lambda^-h^-(n^+,n^-)\nabla n^+\text{div}u^- Id +\lambda^-k^-(n^+,n^-)\nabla n^-\text{div}u^-Id\\
    &+\mu^-l^-(n^+,n^-)\Delta u^-+(\mu^-+\lambda^-)l^-(n^+,n^-)\nabla\text{div}u^-,\\
\end{aligned}
    \right.
\end{equation}
where the nonlinear functions of $(n^+,n^-)$ are  defined  by
\begin{equation}
    \left\{
    \begin{aligned}
        g^+(n^+,n^-)&=\frac{(\mathcal{C}^{2}\rho^-)(1+n^+,1+n^-) }{\rho^+(1+n^+,1+n^-)}-\frac{(\mathcal{C}^{2}\rho^-)(1,1) }{\rho^+(1,1)},\\
        g^-(n^+,n^-)&= \frac{(\mathcal{C}^{2}\rho^+)(1+n^+,1+n^-) }{\rho^-(1+n^+,1+n^-)}  -\frac{(\mathcal{C}^{2}\rho^+)(1,1)}{\rho^-(1,1)}+\frac{(\mathcal{C}^{2}\alpha^+)(1,1)f'(1)}{s_+^2(1,1)}\\
        &\quad-\frac{(\mathcal{C}^{2}\alpha^+)(1+n^+,1+n^-)f'(1+n^-)}{s_+^2(1+n^+,1+n^-)} ,
    \end{aligned}
    \right.
\end{equation}
\begin{equation}
    \left\{
    \begin{aligned}
        \Tilde{g}^+(n^+,n^-)=& \mathcal{C}^{2}(1+n^+,1+n^-)-\mathcal{C}^{2}(1,1)+\frac{(\mathcal{C}^{2}\alpha^-)(1,1)f'(1)}{s_-^2(1,1)}\\
        &-\frac{(\mathcal{C}^{2}\alpha^-)(1+n^+,1+n^-)f'(1+n^-)}{s_-^2(1+n^+,1+n^-)},\\
         \Tilde{g}^-(n^+,n^-)=& \mathcal{C}^{2}(1+n^+,1+n^-)-\mathcal{C}^{2}(1,1),
    \end{aligned}
    \right.
\end{equation}

\begin{equation}
    \left\{
    \begin{aligned}
        h^+(n^+,n^-)&= \frac{\mathcal{C}^{2}\alpha^-(1+n^+,1+n^-)}{(1+n^+)s_-^2(1+n^+,1+n^-)},\\
         h^-(n^+,n^-)&= \frac{\mathcal{C}^{2} (1+n^+,1+n^-)}{(\rho^-s_-^2)(1+n^+,1+n^-)},
    \end{aligned}
    \right.
\end{equation}
\begin{equation}
    \left\{
\begin{aligned}
    k^+(n^+,n^-)&=-\left[ \frac{\mathcal{C}^{2} (1+n^+,1+n^-)}{(s_+^2\rho^+)(1+n^+,1+n^-)}+\frac{f'(1+n^-)(\alpha^-\mathcal{C}^{2})(1+n^+,1+n^-)}{(\rho^+s_+^2s_-^2)(1+n^+,1+n^-)}\right],\\
     k^-(n^+,n^-)&=\frac{\mathcal{C}^{2}\alpha^+(1+n^+,1+n^-)}{(1+n^-)s_+^2(1+n^+,1+n^-)} + \frac{f'(1+n^-)(\alpha^+\mathcal{C}^{2})(1+n^+,1+n^-)}{(\rho^-s_+^2s_-^2)(1+n^+,1+n^-)},
\end{aligned}
    \right.
\end{equation}
and 
\begin{equation}
    \left\{
    \begin{aligned}
        l^+(n^+,n^-)&= \frac{1}{\rho^+(1+n^+,1+n^-)}-\frac{1}{\rho^+(1,1)},\\
         l^-(n^+,n^-)&=  \frac{1}{\rho^-(1+n^+,1+n^-)}-\frac{1}{\rho^-(1,1)}.
    \end{aligned}
    \right.
\end{equation}
We introduce the rescaled variables 
\begin{align*}
    \tilde{n}^+= \beta_1n^+, \  \tilde{n}^-= \beta_4n^-, \ \tilde{u}^+=\sqrt{\beta_1}u^+,\ \tilde{u}^-= \sqrt{\beta_4}u^-,
\end{align*}
and subsequently omit the tildes for simplicity,
and setting 
\begin{align*}
    \gamma_1=\sqrt{\beta_1},\ \gamma_2=\frac{\beta_2\sqrt{\beta_1}}{\beta_4}, \gamma_3=\frac{\beta_3\sqrt{\beta_4}}{\beta_1},\ \gamma_4=\sqrt{\beta_4},
\end{align*}
and 
\begin{align*}
    \gamma_+=\sqrt{\frac{\gamma_1}{\gamma_2}}, \  \gamma_-=\sqrt{\frac{\gamma_4}{\gamma_3}},
\end{align*}
the system \eqref{equ:fornu-pm-S} is reformulated as 
\begin{equation}  \label{equ:fornu-pm}
\left\{\begin{array}{l}
\partial_{t} n^{+}+\gamma_1\operatorname{div} u^{+}=\mathcal{F}_1^+, \\
\partial_{t} u^{+}+\gamma_1 \nabla n^{+}+\gamma_2 \nabla n^{-}-\nu_{1}^{+} \Delta u^{+}-
\nu_{2}^{+} \nabla \operatorname{div} u^{+}=\mathcal{F}_2^+, \\
\partial_{t} n^{-}+\gamma_4\operatorname{div} u^{-}=\mathcal{F}_1^-, \\
\partial_{t} u^{-}+\gamma_3 \nabla n^{+}+\gamma_4 \nabla n^{-}-\nu_{1}^{-} \Delta u^{-}-
\nu_{2}^{-} \nabla \operatorname{div} u^{-}=\mathcal{F}_2^-,
\end{array}\right.
\end{equation}
with the initial data
\begin{align}\label{initial-nu}
    (n^+,u^+,n^-,u^-)|_{t=0}=(n_0^+,u_0^+,n_0^-,u_0^-)(x)  \rightarrow  (0,\vec{0},0,\vec{0}),\  \mbox{as }\  |x|\to\infty,
\end{align}
where the initial data are rescaled accordingly.
Here the source terms are
\begin{equation}
    \left\{
\begin{aligned}
    \mathcal{F}_1^+&=\beta_1S_1^+\left(\frac{n^+}{\beta_1},\frac{u^+}{\sqrt{\beta_1}}\right),\\
    \mathcal{F}_2^+&=\sqrt{\beta_1}S_2^+\left(\frac{n^+}{\beta_1},\frac{u^+}{\sqrt{\beta_1}},\frac{n^-}{\beta_4},\frac{u^-}{\sqrt{\beta_4}}\right),\\
    \mathcal{F}_1^-&=\beta_4S_1^-\left(\frac{n^-}{\beta_4},\frac{u^-}{\sqrt{\beta_4}}\right),\\
    \mathcal{F}_2^-&=\sqrt{\beta_4}S_2^-\left(\frac{n^+}{\beta_1},\frac{u^+}{\sqrt{\beta_1}},\frac{n^-}{\beta_4},\frac{u^-}{\sqrt{\beta_4}}\right).
\end{aligned}
    \right.
\end{equation}
A direct computation gives
\begin{align}
    \gamma_1\gamma_4-\gamma_2\gamma_3=\frac{\beta_1\beta_4-\beta_2\beta_3}{\sqrt{\beta_1\beta_4}}=-\frac{\mathcal{C}^{2}(1,1)f^{\prime}(1)}{\sqrt{\beta_1\beta_4}\rho^+(1,1)}>0.
\end{align}

\subsection{Main results}
Our first result establishes the global well-posedness of the Cauchy problem for the two-fluid system \eqref{equ:fornu-pm} in the critical $L^
p$ Besov space  with small initial data under the stability condition
 \begin{equation}\label{stability-condi}
    -\frac{s_{-}^{2}(1,1)}{\alpha^{-}(1,1)}<f^{\prime}(1)<0.
\end{equation}
We denote 
\begin{equation}
    \begin{aligned}
        \mathcal{E}_{2,p}^{\frac{d}{p},\pm}=\{(n^\pm,u^\pm)\in \mathcal{C}(\mathbb{R}_{+};\dot{\mathcal{B}}_{2,\,p}^{\frac{d}{2}-1,\,\frac{d}{p}})\cap L^1(\mathbb{R}_{+};\dot{\mathcal{B}}_{2,\,p}^{\frac{d}{2}+1,\,\frac{d}{p}} )\times\mathcal{C}(\mathbb{R}_{+};\dot{\mathcal{B}}_{2,\,p}^{\frac{d}{2}-1,\,\frac{d}{p}-1} )\cap L^1(\mathbb{R}_{+};\dot{\mathcal{B}}_{2,\,p}^{\frac{d}{2}+1,\,\frac{d}{p}+1 )} \}
    \end{aligned}
\end{equation}
with norm
\begin{align}
\|(n^\pm,u^\pm)\|_{\mathcal{E}_{2,p}^{\frac{d}{p}}}=\|n^\pm \| _{\widetilde{L}^{\infty}_t \dot{\mathcal{B}}_{2,\,p}^{\frac{d}{2}-1,\,\frac{d}{p}} \cap L^{1}_t\dot{\mathcal{B}}_{2,\,p}^{\frac{d}{2}+1,\,\frac{d}{p}} } 
        &+\|u^\pm\| _{\widetilde{L}^{\infty}_t\dot{\mathcal{B}}_{2,\,p}^{\frac{d}{2}-1,\,\frac{d}{p}-1}\cap  L^{1}_t\dot{\mathcal{B}}_{2,\,p}^{\frac{d}{2}+1,\,\frac{d}{p}+1}},
\end{align}
where $\dot{\mathcal{B}}_{2,\,p}^{s,\,\sigma}$ is defined in Definition \ref{defforB}.
We are ready to formulate our first main theorem.
\begin{theorem}\label{Thm1}
Let $d\geq2$ and   $ p$ satisfy  
\begin{align}\label{range-p}
    2\leq p\leq \min\{4,\frac{2d}{d-2}\} \mbox{ and, additionally, }p\not=4 \mbox{ if } d=2.
 \end{align}  
Under the stability condition \eqref{stability-condi}, there exists a small constant $\epsilon_0=\epsilon_0(p,\,d,\, \lambda,\,  \mu,\, P)>0$ such that if
\begin{equation}
    \begin{aligned}
     \|(n_0^+,n_0^-)\|_{\dot{\mathcal{B}}_{2,\,p}^{\frac{d}{2}-1,\,\frac{d}{p}}}  + \| ( u_0^+,u_0^-) \|_{ \dot{\mathcal{B}}_{2,\,p}^{\frac{d}{2}-1,\,\frac{d}{p}-1}}\leq \epsilon_0,
    \end{aligned}
\end{equation}
then the Cauchy problem \eqref{equ:fornu-pm}-\eqref{initial-nu} has a global   solution $(n^+,n^-,u^+,u^-)$ such that 
\begin{equation} 
\left\{
    \begin{aligned}
       & n^{\pm}\in \mathcal{C}(\mathbb{R}_{+};\dot{\mathcal{B}}_{2,\,p}^{\frac{d}{2}-1,\,\frac{d}{p}})\cap L^1(\mathbb{R}_{+};\dot{\mathcal{B}}_{2,\,p}^{\frac{d}{2}+1,\,\frac{d}{p}} ), \\
       & u^{\pm}\in \mathcal{C}(\mathbb{R}_{+};\dot{\mathcal{B}}_{2,\,p}^{\frac{d}{2}-1,\,\frac{d}{p}-1} )\cap L^1(\mathbb{R}_{+};\dot{\mathcal{B}}_{2,\,p}^{\frac{d}{2}+1,\,\frac{d}{p}+1} ).
    \end{aligned}
    \right.
\end{equation}
Moreover, we have, for any $t>0$,
\begin{equation}
    \begin{aligned}
        \| (n^+,u^+)  \| _{ \mathcal{E}_{2,p}^{\frac{d}{p},+} } 
        &+\|(n^-,u^-)\| _{\mathcal{E}_{2,p}^{\frac{d}{p},-} }
        \leq\, C ( \|(n_0^+,n_0^-)\|_{\dot{\mathcal{B}}_{2,\,p}^{\frac{d}{2}-1,\,\frac{d}{p}}}  + \| ( u_0^+,u_0^-) \|_{ \dot{\mathcal{B}}_{2,\,p}^{\frac{d}{2}-1,\,\frac{d}{p}-1}}).
    \end{aligned}
\end{equation}
for some constant $C=C(p,\,d,\, \lambda,\,  \mu,\, P)> 0$.

In particular, if $2\leq p\leq d$, then the uniqueness holds in $\mathcal{E}_{2,p}^{\frac{d}{p},\pm}$.

\end{theorem}
\begin{remark}
Compared with Theorem 1.3 in \cite{CPAM2010}, Theorem \ref{Thm1} further generalizes the single phase flow to the compressible non-conservative two-fluid  system \eqref{model-main}.
In physical dimensions $d= 2, 3$, the value of $p$ enables us to
consider the case $p > d$ for which the velocity regularity exponent $\frac{d}{p}-1$ becomes negative. Consequently, Theorem \ref{Thm1} can be applied to large highly oscillating initial velocities.
\end{remark}
Let us point out new ingredients in the proof of Theorem \ref{Thm1}. For usual compressible N-S equations, the major difficulty stems from the
convection term in the density equation, as it may cause a loss of one derivative of the density. To overcome it, previous proofs heavily relied on a paralinearized
version combined with a Lagrangian change of variables. For the compressible two--fluid model with constant viscosity coefficients and unequal pressure functions $P^{+}\left(\rho^{+}\right)\not=P^{-}\left(\rho^{-}\right)$, we seem unable to use a single Lagrangian transformation to absorb the convective term of the density. Moreover, due to the non-conservative terms, we cannot apply two separate Lagrangian transformations associated with the two velocities to obtain the desired estimates. 
To address this, we introduce the following “two effective velocities”, 
\begin{equation}\label{TEvelo}
\left\{
    \begin{aligned}
        w^{+}=\nabla (-\Delta)^{-1}\left( \frac{\gamma_1}{\nu_1^++\nu_2^+}n^+ +\frac{\gamma_2}{\nu_1^++\nu_2^+}n^--\text{div}u^+\right),\\
          w^{-}=\nabla (-\Delta)^{-1}\left( \frac{\gamma_3}{\nu_1^-+\nu_2^-}n^+ +\frac{\gamma_4}{\nu_1^-+\nu_2^-}n^--\text{div}u^-\right).
    \end{aligned}
    \right.
\end{equation}
In fact, the definition of $w^\pm$ is a little different as that in \cite{MR2823668}. On the other hand, $(w^+,w^-)$ satisfies the following heat equations
\begin{equation} 
 \left\{
   \begin{array}{ll}
     \partial_t w^+-(\nu_1^++\nu_2^+)\Delta w^+  =g_1,\\
     \partial_t w^- -(\nu_1^++\nu_2^+)\Delta w^-  =g_2,
   \end{array}
 \right.
\end{equation}
where $g_1$ and $g_2$ are defined in \eqref{def:g1g2}. Meanwhile, $(n^+,n^-)$ satisfies  the following damping system
\begin{equation}
\left\{
    \begin{aligned}
      \partial_t n^++u^+\cdot\nabla n^+ +   \frac{\gamma_1}{\nu_1^++\nu_2^+}n^+ +\frac{\gamma_2}{\nu_1^++\nu_2^+}n^-=-n^+\text{div}u^+-\text{div}w^+,\\
        \partial_t n^-+u^-\cdot\nabla n^-+  \frac{\gamma_3}{\nu_1^-+\nu_2^-}n^+ +\frac{\gamma_4}{\nu_1^-+\nu_2^-}n^-=-n^-\text{div}u^--\text{div}w^-.
    \end{aligned}
    \right.
\end{equation}
To control the $\widetilde{L}^{\infty}_t(\dot{B}_{p,1}^{\frac{d}{p} })\cap  L^{1}_t(\dot{B}_{p,1}^{\frac{d}{p}})$ norm of $(n^+,n^-)$, we make full use of $\gamma_+^2\gamma_-^2>1$, which can be derived from $\gamma_1\gamma_4-\gamma_2\gamma_3>0$ due to $f^{\prime}(1)<0.$ Indeed, we make a little analysis on condition $\gamma_+^2\gamma_-^2>1$, i.e.  there exists constant $\tilde{\eta}>0$, such that $\gamma_1>\tilde{\eta}\gamma_2$ and $\tilde{\eta}\gamma_4>\gamma_3.$ This is enough to bound the $\widetilde{L}^{\infty}_t(\dot{B}_{p,1}^{\frac{d}{p} })\cap  L^{1}_t(\dot{B}_{p,1}^{\frac{d}{p}})$ norm of $(n^+,n^-)$.

A natural question after Theorem \ref{Thm1} concerns the large-time behavior of the global solutions constructed above. One may wonder what their asymptotic behavior is as $t\to \infty$.
When the smallness of the low-frequency part of the initial data in $\dot{B}^{\sigma_0}_{2,\infty}$ is proposed, we can also derive the optimal time-decay rates of the global solution to the Cauchy problem.
\begin{theorem}\label{Thm2}
    Let $(n^+,n^-,u^+,u^-)$  be the corresponding global solution to the Cauchy problem \eqref{equ:fornu-pm}-\eqref{initial-nu} that was constructed in Theorem \ref{Thm1} and $p$ satisfies \eqref{range-p}. Let the real number $\sigma_1$ satisfy
    \begin{align}
        \sigma_0\leq\sigma_1<\frac{d}{2}-1 \ \text{with}\ \sigma_0=\frac{d}{2}-\frac{2d}{p}.
    \end{align}
    There exists a positive constant $c=c(p,d,\lambda,\mu,P^\pm)$ such that if in addition the initial
data $(n_0^+,n_0^-,u_0^+,u_0^-)$ fulfill 
\begin{align}
    \mathcal{Y}_p(0)\triangleq\|(n_0^+,n_0^-,u_0^+,u_0^-)\|^{\ell}_{ \dot{B}_{2,\infty}^{\sigma_1}}\leq c,
\end{align}
then it holds that 
\begin{align}\label{ineqtime}
    \mathcal{Y}_p(t)\lesssim\mathcal{Y}_p(0)+\|(\nabla n_0^+,\nabla n_0^-,u_0^+,u_0^-)\|^h_{\dot{B}_{p,1}^{\frac{d}{p}-1}}   \mbox{  for all } t\geq 0,
\end{align}
where the functional $ \mathcal{Y}_p(t)$ is defined by
\begin{equation}
    \begin{aligned}
    \mathcal{Y}_p(t)\triangleq&\sup_{\sigma\in[\sigma_1+\varepsilon,\frac{d}{2}+1]}\|\langle \tau \rangle^{\frac{\sigma-\sigma_1}{2}}(n^+,n^-,u^+,u^-)\|_{L_t^{\infty}(\dot{B}_{2,1}^{\sigma})}^{\ell}+ \|\langle \tau \rangle^{\alpha}(n^+,n^-) \|_{\widetilde{L}_t^{\infty}(\dot{B}_{p,1}^{\frac{d}{p}})}^{h} \\
    &+\|\langle \tau\rangle^{\alpha}(u^+,u^-)\|_{\widetilde{L}_t^{\infty}(\dot{B}_{p,1}^{\frac{d}{p}-1})}^{h}+\| \tau^{\alpha}\nabla(u^+,u^-)\|_{\widetilde{L}_t^{\infty}(\dot{B}_{p,1}^{\frac{d}{p}})}^{h}
\end{aligned}
\end{equation}
with $\alpha=\frac{d}{2}+\frac{1}{2}-\sigma_1-\epsilon $ for $\epsilon>0$ sufficiently small.
\end{theorem}
\begin{remark}
    Theorem \ref{Thm2}  concerns the large-time smoothness of solutions for the non-conservative compressible two-fluid system, revealing the regularizing effect of the viscosity  with  small initial data. It is shown that the solution becomes {\emph{globally smooth}} and  its derivatives decay at the rates that coincide with the optimal rates for the heat equation. 
In particular,  it follows from the embedding $ \dot{B}_{2,1}^{\sigma+d(\frac12-\frac1p)}  \hookrightarrow \dot{B}_{p,1}^{\sigma}  $ that
\begin{align*}
\big\|\Lambda^\sigma (n^+,n^-)\big\|_{L^p}
\;\lesssim \;(\mathcal{Y}_p(0)+\|(\nabla n_0^+,\nabla n_0^-,u_0^+,u_0^-)\|^h_{\dot{B}_{p,1}^{\frac{d}{p}-1}}  )\langle t\rangle^{-\tfrac{\sigma-\widetilde{\sigma}_1}{2}} \  \mbox{ if }\  \widetilde{\sigma}_1<\sigma\leq\frac{d}{p}
\end{align*}
and 
\begin{align*}
\big\|\Lambda^\sigma (u^+,u^-)\big\|_{L^p}
\;\lesssim \;(\mathcal{Y}_p(0)+\|(\nabla n_0^+,\nabla n_0^-,u_0^+,u_0^-)\|^h_{\dot{B}_{p,1}^{\frac{d}{p}-1}}  )\langle t\rangle^{-\tfrac{\sigma-\widetilde{\sigma}_1}{2}} \ \mbox{ if }\  \widetilde{\sigma}_1<\sigma\leq\frac{d}{p}+1,
\end{align*}
where $\widetilde{\sigma}_1=\sigma_1- d(\frac12-\frac1p)$ with $p\geq 2$.
\end{remark}

\subsection{Notations and outline}\label{notations}
 We list some notations that are used frequently throughout the paper.
For simplicity, $C$ denotes a generic positive constant. $A\lesssim B$ ($A\gtrsim B$) means that $A\leq C B$ ($A\geq C B$),  while $A\sim B$ means that both $A\lesssim  B$
and $A\gtrsim B$. 
For $A=(A_{i,j})_{1\leq i,j\leq d}$ and $B=(B_{i,j})_{1\leq i,j\leq d}$ two $d\times d$ matrices, we denote $A:B={\rm{Tr}}\,AB=\sum_{i,j}A_{i,j}B_{j,i}$. 
 For any Banach space $X$ and the functions $g,h\in X$, let $\|(g,h)\|_{X}:=\|g\|_{X}+\|h\|_{X}$. For $p\in[1, \infty]$ and
$T>0$, the notation $L^p(0, T; X)$ or $L^p_T(X)$ designates the set of measurable functions $f: [0, T]\to X$ with $t\mapsto\|f(t)\|_X$ in $L^p(0, T)$, endowed with the norm $\|\cdot\|_{L^p_{T}(X)} :=\|\|\cdot\|_X\|_{L^p(0, T)}$, and $C([0,T];X)$ denotes the set of continuous functions $f: [0, T]\to X$. Let $\mathcal{F}(f)=\widehat{f}$ and $\mathcal{F}^{-1}(f)=\breve{f}$ be the Fourier transform of $f$ and its inverse. We denote $\mathcal{P}:={\rm Id}+(-\Delta)^{-1}\nabla\operatorname{div}$ as the incompressible projector and $\mathcal{Q}=I-\mathcal{P}$ as the compressible projector, respectively.

\vspace{2mm}

For the convenience of readers,
 we shall collect some basic facts on  Littlewood-Paley theory in this section to make the paper as self-contained as possible. The reader may refer to Chapters 2 and 3 in \cite{bahouri1} for more details.

Let $\chi(\xi)$ be 
 a smooth, radial, non-increasing function that is  compactly supported in $B(0,\frac{4}{3})$ and satisfies $\chi(\xi)=1$ in $B(0,\frac{3}{4})$. Then $\varphi(\xi):=  \chi(\frac{\xi}{2})-\chi(\xi)$ verifies
$$
\sum_{j\in \mathbb{Z}}\varphi(2^{-j}\xi)=1\quad\mbox{for}\ \ \xi\neq 0, \text{and}\quad  \operatorname{Supp} \varphi\subset \Big\{\xi\in\mathbb{R}^{d}~|~\frac{3}{4}\leq |\xi|\leq \frac{8}{3}\Big\}.
$$
For any $j\in \mathbb{Z}$, we define the homogeneous dyadic blocks $\dot{\Delta}_{j}$ by
$$
\dot{\Delta}_{j}u:=  \mathcal{F}^{-1}\big{(} \varphi(2^{-j}\cdot )\mathcal{F}(u) \big{)}=2^{jd}h(2^{j}\cdot)\star u\quad\text{with}\quad h:=  \mathcal{F}^{-1}\varphi,
$$
which satisfy $\dot{\Delta}_{j}\dot{\Delta}_{\ell}u=0$ if $|j-\ell|\geq2$.  We also define the low-frequency cut-off
$$
\dot{S}_j u:= \mathcal{F}^{-1}\big{(} \varphi(2^{-j}\cdot )\mathcal{F}(u) \big{)}=\sum_{j'\leq j-1}\dot{\Delta}_j u.
$$
We have the Littlewood-Paley decomposition $u=\sum_{j\in \mathbb{Z}}\dot{\Delta}_{j}u$
that holds in $\mathcal{S}_h'$, i.e.,  tempered distributions $\mathcal{S}'$ such that $\lim\limits_{j\rightarrow-\infty}\|\dot{S}_ju\|_{L^{\infty}}=0$ holds for any $u\in \mathcal{S}_h'$. 

We denote by $\mathcal{Z}'(\mathbb{R}^n)$  the dual space of 
$$
\mathcal{Z}(\mathbb{R}^n)\triangleq\{f\in \mathcal{S}(\mathbb{R}^n): \partial^{\alpha}\hat{f}(0), \forall \alpha\in (\mathbb{N} \cup 0)^n\}.
$$
Firstly, we give the definition of the homogeneous Besov space.
\begin{defn}
For $s\in \mathbb{R}$ and $1\leq p,r\leq \infty$, the  homogeneous Besov norm $\|\cdot\|_{\dot{B}^{s}_{p,r}}$ is defined by
$$
\|u\|_{\dot{B}^{s}_{p,r}}:=  \|\{2^{js}\|\dot{\Delta}_{j}u\|_{L^{p}}\}_{j\in\mathbb{Z}}\|_{\ell^{r}}.
$$
\end{defn}
Secondly, we introduce the hybrid Besov space that will be used in this paper.
\begin{defn}\label{defforB}
For $s,\sigma\in \mathbb{R}$ and $1\leq p \leq \infty$. The   hybrid   Besov space  $\mathcal{\dot{B}}_{2,\,p}^{s,\,\sigma}$ is defined by  
$$\mathcal{\dot{B}}_{2,\,p}^{s,\,\sigma}\triangleq\{f \in \mathcal{Z}'(\mathbb{R}^n): \|f\|_{\mathcal{\dot{B}}_{2,\,p}^{s,\,\sigma}}<\infty\}$$
with 
$$
\|f\|_{\mathcal{\dot{B}}_{2,\,p}^{s,\,\sigma}}:=  \sum_{2^j\leq R_0} 2^{js}\|\dot{\Delta}_{j}f\|_{L^{2}}+\sum_{2^j> R_0} 2^{j\sigma}\|\dot{\Delta}_{j}f\|_{L^{p}},
$$
where $R_0$ is a fixed and sufficiently large constant.
\end{defn}
Since we are concerned with time-dependent functions valued in Besov spaces, the following space-time mixed norm is usually mentioned: 
$$\|u\|_{L^q_T\mathcal{\dot{B}}_{2,\,p}^{s,\,\sigma}}\triangleq\Big\|\|u(t,\cdot)\|_{\mathcal{\dot{B}}_{2,\,p}^{s,\,\sigma}}\Big\|_{L^q_T}.$$
Here, we introduce another space-time mixed Besov norm, which is referred to
Chemin-Lerner’s spaces. The definition is given by in such way
$$\|u\|_{\widetilde{L}^q_T\mathcal{\dot{B}}_{2,\,p}^{s,\,\sigma}}\triangleq \sum_{2^j\leq R_0}2^{js} \|\dot{\Delta}_j u \|_{L^q(0,T;L^2)} +\sum_{2^j> R_0}2^{j\sigma} \|\dot{\Delta}_j u \|_{L^q(0,T;L^p)}.$$
The index $T$ will be omitted if $T =+\infty$ and we shall denote by $\widetilde{\mathcal{C}}_b(\mathcal{\dot{B}}_{2,\,p}^{s,\,\sigma})$ the subset of functions $\widetilde{L}^{\infty}(\mathcal{\dot{B}}_{2,\,p}^{s,\,\sigma})$ which are continuous from $\mathbb{R}_+$ to $\mathcal{\dot{B}}_{2,\,p}^{s,\,\sigma}$. It is easy to check
that $\widetilde{L}_T^1\mathcal{\dot{B}}_{2,\,p}^{s,\,\sigma}=L^1_T\mathcal{\dot{B}}_{2,\,p}^{s,\,\sigma} $ and $\widetilde{L}_T^q\mathcal{\dot{B}}_{2,\,p}^{s,\,\sigma} \subseteq L^q_T\mathcal{\dot{B}}_{2,\,p}^{s,\,\sigma} $ for $q>1$.

Also, for a tempered distribution $f$, we denote 
$$f^\ell\triangleq\sum_{2^j\leq R_0}\dot{\Delta}_j f, \qquad f^h\triangleq f-f^\ell, $$
and 
$$ \|f\|^\ell_{\dot{B}_{p,1}^s}=\sum_{2^j\leq R_0} 2^{js}\|\dot{\Delta}_j f\|_{L^p},\qquad\|f\|^h_{\dot{B}_{p,1}^s}=\sum_{2^j> R_0} 2^{j\sigma}\|\dot{\Delta}_j f\|_{L^p},$$
$$ \|f\|^\ell_{\widetilde{L}_T^q\dot{B}_{p,1}^s}=\sum_{2^j\leq R_0} 2^{js}\|\dot{\Delta}_j f\|_{L^q(0,T;L^p)},\qquad\|f\|^h_{\widetilde{L}^q_T\dot{B}_{p,1}^s}=\sum_{2^j> R_0} 2^{j\sigma}\|\dot{\Delta}_j f\|_{L^q(0,T;L^p)}.$$
Here and below, all functional spaces
will be considered in
$\mathbb{R}^{d}$, we shall omit the space dependence for simplicity. Due to Bernstein's inequalities, there are  well-known  properties in Besov spaces.
\begin{lemma}
Let $s,\sigma\in\mathbb{R}$ and $1\leq p\leq +\infty$. Then we have
\begin{itemize}
    \item $\mathcal{\dot{B}}_{2,\,p}^{s_2,\,\sigma}\subset \mathcal{\dot{B}}_{2,\,p}^{s_1,\,\sigma}$ for $s_1\geq s_2$ and
    $\mathcal{\dot{B}}_{2,\,p}^{s,\,\sigma_2}\subset \mathcal{\dot{B}}_{2,\,p}^{s,\,\sigma_1}$ for $\sigma_1\leq \sigma_2$.
    
    \item Interpolation: For $s_1,s_2,\sigma_1,\sigma_2\in\mathbb{R}$ and $\theta\in[0,1]$, we have
    \[
    \|f\|_{\mathcal{\dot{B}}_{2,\,p}^{\theta s_1+(1-\theta)s_2,\,\theta\sigma_1+(1-\theta)\sigma_2}}
    \leq
    \|f\|_{\mathcal{\dot{B}}_{2,\,p}^{s_1,\,\sigma_1}}^{\theta}
    \|f\|_{\mathcal{\dot{B}}_{2,\,p}^{s_2,\,\sigma_2}}^{1-\theta}.
    \]
    
    \item Embedding:
   $ \mathcal{\dot{B}}_{2,\,p}^{n/2,\,n/p} \hookrightarrow L^\infty ,$
    $\dot{B}_{2,\,1}^{s}\hookrightarrow \mathcal{\dot{B}}_{2,\,p}^{s,\,s-n/2+n/p}
    \hookrightarrow \dot{B}_{p,\,1}^{\,s-n/2+n/p},
    \qquad \text{for } p\geq 2. $
\end{itemize}
\end{lemma}

\begin{lemma}\label{lemma62}
The following properties hold{\rm:}
\begin{itemize}
\item{} For any $s\in\mathbb{R}$, $1\leq p_{1}\leq p_{2}\leq \infty$ and $1\leq r_{1}\leq r_{2}\leq \infty$, it holds
\begin{equation}\nonumber
\begin{aligned}
\dot{B}^{s}_{p_{1},r_{1}}\hookrightarrow \dot{B}^{s-d\left(\frac{1}{p_{1}}-\frac{1}{p_{2}}\right)}_{p_{2},r_{2}}.
\end{aligned}
\end{equation}
\item{} For any $1\leq p\leq q\leq\infty$, we have the following chain of continuous embedding{\rm:}
\begin{equation}\nonumber
\begin{aligned}
\dot{B}^{0}_{p,1}\hookrightarrow L^{p}\hookrightarrow \dot{B}^{0}_{p,\infty}\hookrightarrow \dot{B}^{\sigma}_{q,\infty} \quad \text{for} \quad\sigma=-d\left(\frac{1}{p}-\frac{1}{q}\right)<0.
\end{aligned}
\end{equation}

\item{}
For any $\sigma\in \mathbb{R}$, the Fourier multiplier $\Lambda^{\sigma}$ with symbol $|\xi|^{\sigma}$ is an isomorphism from $\dot{B}^{s}_{p,r}$ to $\dot{B}^{s-\sigma}_{p,r}$.
\item{} Let $1\leq p_{1},p_{2},r_{1},r_{2}\leq \infty$, $s_{1}\in\mathbb{R}$ and $s_{2}\in\mathbb{R}$ satisfy
    \begin{align}
    s_{2}<\frac{d}{p_{2}}\quad\text{\text{or}}\quad s_{2}=\frac{d}{p_{2}}~\text{and}~r_{2}=1.\nonumber
    \end{align}
    The space $\dot{B}^{s_{1}}_{p_{1},r_{1}}\cap \dot{B}^{s_{2}}_{p_{2},r_{2}}$ endowed with the norm $\|\cdot \|_{\dot{B}^{s_{1}}_{p_{1},r_{1}}}+\|\cdot\|_{\dot{B}^{s_{2}}_{p_{2},r_{2}}}$ is a Banach space. 
\end{itemize}
The following Morse-type product estimates in Besov spaces play a fundamental role in the analysis of nonlinear terms:
\begin{lemma}\label{Lemma1.2}
	The following statements hold:
	\begin{itemize}
		\item{} Let $s>0$ and $1\leq p,r\leq \infty$. Then $\dot{B}^{s}_{p,r}\cap L^{\infty}$ is a algebra and
		\begin{equation}\label{uv1}
			\|uv\|_{\dot{B}^{s}_{p,r}}\lesssim \|u\|_{L^{\infty}}\|v\|_{\dot{B}^{s}_{p,r}}+ \|v\|_{L^{\infty}}\|u\|_{\dot{B}^{s}_{p,r}}.
		\end{equation}
		\item{}
		Let $s_{1}$, $s_{2}$ and $p$ satisfy $2\leq p\leq \infty$, $-\frac{d}{p}<s_{1},s_{2}\leq \frac{d}{p}$ and $s_{1}+s_{2}>0$. Then we have
		\begin{equation}\label{uv2}
			\|uv\|_{\dot{B}^{s_{1}+s_{2}-\frac{d}{p}}_{p,1}}\lesssim \|u\|_{\dot{B}^{s_{1}}_{p,1}}\|v\|_{\dot{B}^{s_{2}}_{p,1}}.
		\end{equation}
		\item{} Assume that $s_{1}$, $s_{2}$ and $p$ satisfy $2\leq p\leq \infty$, $-\frac{d}{p}<s_{1}\leq \frac{d}{p}$, $-\frac{d}{p}\leq s_{2}<\frac{d}{p}$ and $s_{1}+s_{2}\geq0$. Then 
it holds
		\begin{equation}\label{uv3}
			\|uv\|_{\dot{B}^{s_{1}+s_{2}-\frac{d}{p}}_{p,\infty}}\lesssim \|u\|_{\dot{B}^{s_{1}}_{p,1}}\|v\|_{\dot{B}^{s_{2}}_{p,\infty}}.
		\end{equation}
        \item {} Let $s_{1}$, $s_{2}$, $p_1$ and $p_2$ be such that
$$
s_1+s_2>0,\quad s_1\leq \frac{d}{p_1},\quad s_2\leq \frac{d}{p_2},\quad s_1\geq s_2,\quad \frac{1}{p_1}+\frac{1}{p_2}\leq 1.
$$
Then we have
\begin{align}\label{uv4}
\|uv\|_{\dot{B}^{s_{2}}_{q,1}}\lesssim \|u\|_{\dot{B}^{s_{1}}_{p_1,1}}\|v\|_{\dot{B}^{s_{2}}_{p_2,1}}\quad \text{with}\quad \frac{1}{q}=\frac{1}{p_1}+\frac{1}{p_2}-\frac{s_1}{d}.
\end{align}
	\end{itemize}

\end{lemma}

\end{lemma}
At last, we present
an elementary fact which will be used later.
\begin{lemma}\emph{\cite{Danchin-Xu-ARMA2017}}\label{Lem1.2}
    Let $X:[0,T]\to \mathbb{R}_+$ be a continuous function such that $X^p$ is differentiable for some $p\geq 1$ and satisfies
    \begin{equation}
        \begin{aligned}
            \frac{1}{p}\frac{d}{dt}X^p+BX^p\leq AX^{p-1}
        \end{aligned}
    \end{equation}
for some constant $B\geq 0$ and measurable function $A:[0,T]\to \mathbb{R}_+$.

Define $X_{\epsilon}=(X^p+\epsilon^p)^{1/p}$ for $\epsilon>0.$ Then it holds that 
\begin{equation}
    \begin{aligned}
        \frac{d}{dt}X_{\epsilon}+BX_{\epsilon}\leq A +B\epsilon.
    \end{aligned}
\end{equation}
\end{lemma}

\section{Global existence for the  non-conservative compressible two-fluid}
In order to prove the global existence of the
solution to the Cauchy problem \eqref{equ:fornu-pm}, we need to establish the uniform-in-time a
priori estimates below.
We denote the energy norm as follows:
\begin{equation}\label{def-E}
    \begin{aligned}
        \mathcal{X}_p(t) &\triangleq \| (n^+,n^-,u^+,u^-) \|^{\ell}_{\widetilde{L}^{\infty}_t(\dot{B}_{2,1}^{\frac{d}{2}-1})\cap L^{1}_t(\dot{B}_{2,1}^{\frac{d}{2}+1}) }+ \|(n^+,n^-)\|^h_{\widetilde{L}^{\infty}_t(\dot{B}_{p,1}^{\frac{d}{p} })\cap  L^{1}_t(\dot{B}_{p,1}^{\frac{d}{p}})} \\
        &\quad+\|(u^+,u^-)\|^h_{\widetilde{L}^{\infty}_t(\dot{B}_{p,1}^{\frac{d}{p}-1})\cap  L^{1}_t(\dot{B}_{p,1}^{\frac{d}{p}+1})} 
    \end{aligned}
\end{equation}
and 
\begin{equation}\label{def-X0}
    \begin{aligned}
        \mathcal{X}(0)\triangleq \|(n_0^+,n_0^-,u_0^+,u_0^-) \|^{\ell}_{\dot{B}_{2,1}^{\frac{d}{2}-1}}+\|(n_0^+,n_0^-)\|^h_{\dot{B}_{p,1}^{\frac{d}{p} }}+\|(u^+_0,u^-_0)\|^h_{\dot{B}_{p,1}^{\frac{d}{p}-1}}.
    \end{aligned}
\end{equation}
\begin{proposition}\label{Prop2.1}
 Let $p$ satisfies \eqref{range-p}. For given time $T > 0$, let $(n^+,n^-,u^+,u^-) $ be a strong solution to   the Cauchy problem \eqref{equ:fornu-pm} on $t \in [0, T]$ such that
  \begin{align}
      \|(n^+,n^-)\|_{L^\infty([0,T]\times \mathbb{R}^d)}\leq \frac12.
  \end{align}
Then  
   \begin{align}
      \mathcal{X}_p(t) \leq C_0 \mathcal{X}(0), \quad t \in (0, T), 
  \end{align}
where $C_0 > 1$ is a constant independent of   $T$, and  $\mathcal{X}_p(t)$ and   $\mathcal{X}(0)$  are defined by \eqref{def-E} and   \eqref{def-X0}, respectively.
\end{proposition}
Therefore, the proof of Proposition \ref{Prop2.1} consists of Lemmas \ref{Lem2.1}--\ref{Lem2.2} below.
\subsection{Low-frequency estimate}
\begin{lemma}\label{Lem2.1}
Let  $T>0$ be any given time, and $(n^+,n^-,u^+,u^-) $ be any strong solution to the Cauchy problem \eqref{equ:fornu-pm} for $t\in(0,T)$. Then it holds that
 \begin{equation}
    \begin{aligned}
        \| (n^+,n^-,u^+,u^-) \|^{\ell}_{\widetilde{L}^{\infty}_t(\dot{B}_{2,1}^{\frac{d}{2}-1})\cap L^{1}_t(\dot{B}_{2,1}^{\frac{d}{2}+1}) }\lesssim\| (n_0^+,n_0^-,u_0^+,u_0^-) \|^{\ell}_{ \dot{B}_{2,1}^{\frac{d}{2}-1}}+(1+\mathcal{X}_p(t) +\mathcal{X}_p^2(t))\mathcal{X}_p^2(t).
    \end{aligned}
\end{equation}
\end{lemma}
\begin{proof}
Denote $n_j^\pm=\dot{\Delta}_jn^\pm$, $u_j^\pm=\dot{\Delta}_j u^\pm$ for simplicity. By applying $\dot{\Delta}_j$ to \eqref{equ:fornu-pm}, we have 
\begin{equation}  \label{equ:fornu-pm-local}
\left\{\begin{array}{l}
\partial_{t} n_j^{+}+\gamma_1\operatorname{div} u_j^{+}=\dot{\Delta}_j\mathcal{F}_1^+, \\
\partial_{t} u_j^{+}+\gamma_1 \nabla n_j^{+}+\gamma_2 \nabla n_j^{-}-\nu_{1}^{+} \Delta u_j^{+}-
\nu_{2}^{+} \nabla \operatorname{div} u_j^{+}=\dot{\Delta}_j\mathcal{F}_2^+, \\
\partial_{t} n_j^{-}+\gamma_4\operatorname{div} u_j^{-}=\dot{\Delta}_j\mathcal{F}_1^-, \\
\partial_{t} u_j^{-}+\gamma_3 \nabla n_j^{+}+\gamma_4 \nabla n_j^{-}-\nu_{1}^{-} \Delta u_j^{-}-
\nu_{2}^{-} \nabla \operatorname{div} u_j^{-}=\dot{\Delta}_j\mathcal{F}_2^-.
\end{array}\right.
\end{equation}
Multiplying \eqref{equ:fornu-pm-local}$_1$, \eqref{equ:fornu-pm-local}$_2$, \eqref{equ:fornu-pm-local}$_3$ and \eqref{equ:fornu-pm-local}$_4$ by $\frac{\gamma_1}{\gamma_2}n_j^+$, $\frac{\gamma_1}{\gamma_2}u_j^+$, $\frac{\gamma_4}{\gamma_3}n_j^-$ and $\frac{\gamma_4}{\gamma_3}u_j^-$ respectively, we derive
\begin{equation}
    \begin{aligned}
        \frac{1}{2}\frac{d}{dt}&\left\{  \|\gamma_+ n_j^++(\gamma_+)^{-1}n_j^-\|^2_{L^2}+\widetilde{\gamma}\|n_j^-\|^2_{L^2}  +\frac{\gamma_1}{\gamma_2}\|u_j^+\|^2_{L^2}+\frac{\gamma_4}{\gamma_3}\|u_j^-\|^2_{L^2}     \right\}\\
        &+\frac{\gamma_1}{\gamma_2}\left(\nu_1^+\|\nabla u_j^+\|^2_{L^2}+\nu_2^+\|\text{div} u_j^+\|^2_{L^2}\right)+\frac{\gamma_4}{\gamma_3}\left(\nu_1^-\|\nabla u_j^-\|^2_{L^2}+\nu_2^-\|\text{div} u_j^-\|^2_{L^2}\right)\\
        &\lesssim \|\dot{\Delta}_j(\mathcal{F}_1^+,\mathcal{F}_2^+,\mathcal{F}_1^-,\mathcal{F}_2^-)\|_{L^2}\|(n_j^\pm,u_j^\pm)\|_{L^2}.
    \end{aligned}
\end{equation}

To capture the dissipation with respect to $(n^+,n^-)$, we    multiply $\nabla$ \eqref{equ:fornu-pm-local}$_1$, \eqref{equ:fornu-pm-local}$_2$,    $\nabla$ \eqref{equ:fornu-pm-local}$_3$ and \eqref{equ:fornu-pm-local}$_4$ by $\frac{1}{\gamma_2}u_j^+$, $\frac{1}{\gamma_2}\nabla n^+_j$,  $\frac{1}{\gamma_3}u_j^-$  and $\frac{1}{\gamma_3}\nabla n^-$ respectively, we obtain
\begin{equation}
    \begin{aligned}
       \frac{ 1}{\gamma_2} \frac{d}{dt}&(u_j^+|\nabla n_j^+) +\frac{\gamma_1}{\gamma_2}\|\nabla n_j^+\|^2_{L^2}+(\nabla n_j^-|\nabla n_j^+) -\frac{\gamma_1}{\gamma_2}\|\text{div}u_j^+\|^2_{L^2}\\
        &\quad\quad-\frac{\nu_1^+}{\gamma_2}(\Delta u_j^+,\nabla n_j^+)-\frac{\nu_2^+}{\gamma_2}(\nabla\text{div} u_j^+,\nabla n_j^+)\\
        &\lesssim (\|\nabla\dot{\Delta}_j \mathcal{F}_1^+\|_{L^2} +\|\text{div}\dot{\Delta}_j \mathcal{F}_2^+\|_{L^2})\| (n_j^+,u_j^+)\|_{L^2}.
    \end{aligned}
\end{equation}
Similarly, we have
\begin{equation}
    \begin{aligned}
        \frac{ 1}{\gamma_3}\frac{d}{dt}&(u_j^-|\nabla n_j^-) +(\nabla n_j^+|\nabla n_j^-) +\frac{\gamma_4}{\gamma_3}\|\nabla n_j^+\|^2_{L^2}-\frac{\gamma_4}{\gamma_3}\|\text{div}u_j^-\|^2_{L^2}\\
        &\quad\quad-\frac{\nu_1^-}{\gamma_3}(\Delta u_j^-,\nabla n_j^-)-\frac{\nu_2^-}{\gamma_3}(\nabla\text{div} u_j^-,\nabla n_j^-)\\
        &\lesssim (\|\nabla\dot{\Delta}_j \mathcal{F}_1^-\|_{L^2} +\|\text{div}\dot{\Delta}_j \mathcal{F}_2^+\|_{L^2})\| (n_j^-,u_j^-)\|_{L^2}.
    \end{aligned}
\end{equation}
Now, we multiply two small constants $\eta_1$ (to be determined) to (2.3) and  (2.4), respectively, and then add the resulting equations with (2.2) together. Consequently, we are led to the following inequality
    \begin{align}\label{eq:0}
        &\frac{1}{2}\frac{d}{dt} \left\{\|\gamma_+ n_j^++(\gamma_+)^{-1}n_j^-\|^2_{L^2}+\widetilde{\gamma}\|n_j^-\|^2_{L^2}  +\frac{\gamma_1}{\gamma_2}\|u_j^+\|^2_{L^2}+\frac{\gamma_4}{\gamma_3}\|u_j^-\|^2_{L^2}\right.\nonumber\\
       &+ \left.\frac{2\eta_1}{\gamma_2}(u_j^+|\nabla n_j^+) + \frac{2\eta_1}{\gamma_3}(u_j^-|\nabla n_j^-) \right\}\nonumber\\
        &+\frac{\gamma_1}{\gamma_2}\left(\nu_1^+\|\nabla u_j^+\|^2_{L^2}+\nu_2^+\|\text{div} u_j^+\|^2_{L^2}\right)+\frac{\gamma_4}{\gamma_3}\left(\nu_1^-\|\nabla u_j^-\|^2_{L^2}+\nu_2^-\|\text{div} u_j^-\|^2_{L^2}\right)\nonumber\\
        &+\eta_1\|\gamma_+\nabla n_j^++(\gamma_+)^{-1}\nabla n_j^-\|^2_{L^2}+\eta_1 \tilde{\gamma}\|\nabla n_j^-\|_{L^2}^2 -\frac{\eta_1\gamma_1}{\gamma_2}\|\text{div}u_j^+\|^2_{L^2}  -\frac{\eta_1\gamma_4}{\gamma_3}\|\text{div}u_j^-\|^2_{L^2}\\
        &-\frac{\eta_1\nu_1^+}{\gamma_2}(\Delta u_j^+,\nabla n_j^+)-\frac{\eta_1\nu_2^+}{\gamma_2}(\nabla\text{div} u_j^+,\nabla n_j^+)\nonumber\\
        &-\frac{\eta_1\nu_1^-}{\gamma_3}(\Delta u_j^-,\nabla n_j^-)-\frac{\eta_1\nu_2^-}{\gamma_3}(\nabla\text{div} u_j^-,\nabla n_j^-)\nonumber\\
       \lesssim  &\,\|\dot{\Delta}_j(\mathcal{F}_1^+,\mathcal{F}_2^+,\mathcal{F}_1^-,\mathcal{F}_2^-)\|_{L^2}\|(n_j^+,u_j^+,n_j^-,u_j^-)\|_{L^2},\nonumber
    \end{align}
where we use the fact
\begin{equation}
    \begin{aligned}
        &\frac{ \gamma_1}{\gamma_2}\|\nabla n_j^+\|^2_{L^2}+\frac{ \gamma_4}{\gamma_3}\|\nabla n_j^-\|^2_{L^2}+2 (\nabla n_j^-|\nabla n_j^+) \\
        =&\,\|\gamma^+\nabla n_j^++(\gamma^+)^{-1}\nabla n_j^-\|^2_{L^2}+\widetilde{\gamma}\|\nabla n_j^-\|^2_{L^2}
    \end{aligned}
\end{equation}
with 
\begin{align}\label{condition:gamma}
    \tilde{\gamma}=\gamma_-^2-\frac{1}{\gamma_+^2}=\frac{\gamma_1\gamma_4-\gamma_2\gamma_3}{\gamma_1\gamma_3}=-\frac{\beta_1\beta_4-\beta_2\beta_3}{\gamma_1\gamma_3\sqrt{\beta_1\beta_4}}>0.
\end{align}
Let
\begin{equation}
\begin{aligned}
    \mathcal{L}_j(t)&= \frac{1}{2}\left\{\|\gamma^+ n_j^++(\gamma^+)^{-1}n_j^-\|^2_{L^2}+\widetilde{\gamma}\|n_j^-\|^2_{L^2}  +\frac{\gamma_1}{\gamma_2}\|u_j^+\|^2_{L^2}+\frac{\gamma_4}{\gamma_3}\|u_j^-\|^2_{L^2}\right. \\
       &\quad\quad \left.+ \frac{2\eta_1}{\gamma_2}(u_j^+|\nabla n_j^+) + \frac{2\eta_1}{\gamma_3}(u_j^-|\nabla n_j^-) \right\}
\end{aligned}
\end{equation}
and 
\begin{equation}
    \begin{aligned}
        \mathcal{D}_j(t)= &\frac{\gamma_1}{\gamma_2}\left(\nu_1^+\|\nabla u_j^+\|^2_{L^2}+\nu_2^+\|\text{div} u_j^+\|^2_{L^2}\right)
        +\frac{\gamma_4}{\gamma_3}\left(\nu_1^-\|\nabla u_j^-\|^2_{L^2}+\nu_2^-\|\text{div} u_j^-\|^2_{L^2}\right)\\
        &+\eta_1\|\gamma^+\nabla n_j^++(\gamma^+)^{-1}\nabla n_j^-\|^2_{L^2}+\eta_1\widetilde{\gamma}\|\nabla n_j^-\|^2_{L^2}-\frac{\eta_1\gamma_1}{\gamma_2}\|\text{div}u_j^+\|^2_{L^2} -\frac{\eta_2\gamma_4}{\gamma_3}\|\text{div}u_j^-\|^2_{L^2}\\
        &-\frac{\eta_1\nu_1^+}{\gamma_2}(\Delta u_j^+,\nabla n_j^+)-\frac{\eta_1\nu_2^+}{\gamma_2}(\nabla\text{div} u_j^+,\nabla n_j^+) 
        -\frac{\eta_1\nu_1^-}{\gamma_3}(\Delta u_j^-,\nabla n_j^-)-\frac{\eta_1\nu_2^-}{\gamma_3}(\nabla\text{div} u_j^-,\nabla n_j^-).
    \end{aligned}
\end{equation}

Then, we choose a  $\eta_1$ sufficiently small such that
\begin{align}\label{eq:1}
    \mathcal{L}_j^2(t)\sim \| n_j^+ \|^2_{L^2}+ \|n_j^-\|^2_{L^2} + \| u_j^+ \|^2_{L^2}+ \|u_j^-\|^2_{L^2},
\end{align}
and 
\begin{align}\label{eq:2}
     \mathcal{D}_j^2(t)\gtrsim 2^{2j}(\| n_j^+ \|^2_{L^2}+ \|n_j^-\|^2_{L^2} + \| u_j^+ \|^2_{L^2}+ \|u_j^-\|^2_{L^2} ).
\end{align}

Combining \eqref{eq:1}, \eqref{eq:2} with \eqref{eq:0}, we get the following equality
\begin{equation}\label{eq:3}
    \begin{aligned}
        \frac{d}{dt}\mathcal{L}^2_j(t)+2^{2j}\mathcal{L}^2_j(t)\lesssim \|\dot{\Delta}_j(\mathcal{F}_1^+,\mathcal{F}_2^+,\mathcal{F}_1^-,\mathcal{F}_2^-)\|_{L^2}\mathcal{L} _j(t).
    \end{aligned}
\end{equation}
Define $\mathcal{L}^{\delta}_j(t)=(\mathcal{L} _j(t)+\delta^2)^{\frac12}$, by applying Lemma \ref{Lem1.2}, and taking the limit $\delta\to 0$, we  get
\begin{equation}\label{eq:4}
    \begin{aligned}
         &\| (n_j^+,n_j^-,u_j^+,u_j^-) \|_{L^2}+2^{2j}\|(n_j^+,n_j^-,u_j^+,u_j^-)\|_{L^{1}_t(L^2)}\\
        \lesssim&\,\| \dot{\Delta}_j(n_0^+,n_0^-,u_0^+,u_0^-) \|_{ L^2}+\|\dot{\Delta}_j(\mathcal{F}_1^+,\mathcal{F}_2^+,\mathcal{F}_1^-,\mathcal{F}_2^-)\|_{L^{1}_t(L^2)}.
    \end{aligned}
\end{equation}
Multiplying \eqref{eq:4}  by $2^{j(\frac{d}{2}-1)}$, taking the supremum on $[0,t]$ and then summing over $j\leq j_0$, we have
\begin{equation}\label{eq:5}
    \begin{aligned}
        &\| (n^+,n^-,u^+,u^-) \|^{\ell}_{\widetilde{L}^{\infty}_t(\dot{B}_{2,1}^{\frac{d}{2}-1})}+\|(n^+,n^-,u^+,u^-)\|^{\ell}_{L^{1}_t(\dot{B}_{2,1}^{\frac{d}{2}+1})}\\
        \lesssim&\,\| (n_0^+,n_0^-,u_0^+,u_0^-) \|^{\ell}_{ \dot{B}_{2,1}^{\frac{d}{2}-1}}+\|(\mathcal{F}_1^+,\mathcal{F}_2^+,\mathcal{F}_1^-,\mathcal{F}_2^-)\|^{\ell}_{L^{1}_t(\dot{B}_{2,1}^{\frac{d}{2}-1})}.
    \end{aligned}
\end{equation}
Next we begin to bound those nonlinear terms arising in $\mathcal{F}_i^\pm(i=1,2)$.
With regard to $u^\pm\cdot\nabla n^\pm$,  by taking  $s_1 = s_2 = \frac{n}{2}-1$, $t_1 = t_2 =\frac{n}{2}$, $r_1 = r_4 = \infty$, $r_2 = r_3 =
1$, $\gamma=-1$ in \eqref{A2} of Proposition \ref{propA1},  we arrive at
\begin{equation}\label{eq:6}
    \begin{aligned}
        &\|u^\pm\cdot\nabla n^\pm\|^{\ell}_{L^1_t(\dot{B}_{2,1}^{\frac{d}{2}-1})}\\
        \lesssim \;&  \|u^\pm\|_{\widetilde{L}^{\infty}_t(\dot{\mathcal{B}}_{2,\,p}^{\frac{d}{2}-1,\,\frac{d}{p}-1})} \|\nabla n^\pm \|_{L^1_t(\dot{\mathcal{B}}_{2,\,p}^{\frac{d}{2},\, \frac{d}{p}-1})}+\|\nabla n^\pm\|_{\widetilde{L}^{2}_t(\dot{\mathcal{B}}_{2,\,p}^{\frac{d}{2}-1,\,\frac{d}{p}-1})} \|u^\pm \|_{\widetilde{L}^2_t(\dot{\mathcal{B}}_{2,\,p}^{\frac{d}{2},\,\frac{d}{p}})}.
    \end{aligned}
\end{equation}
With regard to $n^+\text{div}u^+$,  by taking  $s_1 = s_2 = \frac{n}{2}-1$, $t_1 = t_2 =\frac{n}{2}$, $r_1 = r_4 = \infty$, $r_2 = r_3 =
1$, $\gamma=0$ in \eqref{A2} of Proposition \ref{propA1}, we have 
 \begin{equation}\label{eq:7}
    \begin{aligned}
        &\|n^\pm\text{div}u^\pm\|^{\ell}_{L^1_t(\dot{B}_{2,1}^{\frac{d}{2}-1})}\\
        \lesssim\; &  \|n^\pm\|_{\widetilde{L}^{\infty}_t(\dot{\mathcal{B}}_{2,\,p}^{\frac{d}{2}-1,\,\frac{d}{p}-1})} \|\text{div}u^\pm \|_{L^1_t(\dot{\mathcal{B}}_{2,\,p}^{\frac{d}{2},\, \frac{d}{p}})}+\|\text{div}u^\pm\|_{\widetilde{L}^{2}_t(\dot{\mathcal{B}}_{2,\,p}^{\frac{d}{2}-1,\,\frac{d}{p}-1})} \|n^\pm\|_{\widetilde{L}^2_t(\dot{\mathcal{B}}_{2,\,p}^{\frac{d}{2},\,\frac{d}{p}})}.
    \end{aligned}
\end{equation}
Thus we obtain
\begin{equation}
    \begin{aligned}
        \|(\mathcal{F}_1^+,\mathcal{F}_1^-)\|^{\ell}_{L^1_t(\dot{B}_{2,1}^{\frac{d}{2}-1})}\lesssim\,\mathcal{X}_p^2(t).
    \end{aligned}
\end{equation}

With regard to $\mathcal{F}_2^+\sim g^+(n^+,n^-)\,\nabla n^+ + (u^+\cdot\nabla)u^+
      + h^+(n^+,n^-)\,\nabla n^+\nabla u^+
      + l^+(n^+,n^-)\,\Delta u^+ $.  By choosing $s_1 = s_2 = \frac{n}{2}-1$, $t_1 = t_2 =\frac{n}{2}$, $r_1 = r_4 = \infty$, $r_2 = r_3 =
1$, $\gamma=-1$ in \eqref{A2} of Proposition \ref{propA1}, we get
\begin{equation}\label{eq:8}
    \begin{aligned}
        &\|g^+(n^+,n^-)\nabla n^+\|^{\ell}_{L^1_t(\dot{B}_{2,1}^{\frac{d}{2}-1})}\\\leq &\,\|g^+(n^+,n^-)\|_{\widetilde{L}^{\infty}_t(\dot{\mathcal{B}}_{2,\,p}^{\frac{d}{2}-1,\,\frac{d}{p}-1})}\|\nabla n^+\|_{L^1_t(\dot{\mathcal{B}}_{2,\,p}^{\frac{d}{2},\, \frac{d}{p}-1})}+\|\nabla n^+\|_{\widetilde{L}^{2}_t(\dot{\mathcal{B}}_{2,\,p}^{\frac{d}{2}-1,\,\frac{d}{p}-1})}\|g^+(n^+,n^-) \|_{\widetilde{L}^2_t(\dot{\mathcal{B}}_{2,\,p}^{\frac{d}{2},\,\frac{d}{p}})}\\
        \lesssim&\,(1+\|(n^+,n^-)\|_{\widetilde{L}^{\infty}_t(\dot{\mathcal{B}}_{2,\,p}^{\frac{d}{2}-1,\,\frac{d}{p}})})^{\frac{n}{2}+1} \Big(\| (n^+,n^-)\|_{\widetilde{L}^{\infty}_t(\dot{\mathcal{B}}_{2,\,p}^{\frac{d}{2}-1,\,\frac{d}{p}})}\|  n^+ \|_{L^1_t(\dot{\mathcal{B}}_{2,\,p}^{\frac{d}{2}+1,\, \frac{d}{p}})}\\
        &\qquad\qquad\qquad\qquad\qquad\qquad\qquad\qquad+\|  n^+\|_{\widetilde{L}^{2}_t(\dot{\mathcal{B}}_{2,\,p}^{\frac{d}{2},\,\frac{d}{p}})}\| (n^+,n^-) \|_{\widetilde{L}^2_t(\dot{\mathcal{B}}_{2,\,p}^{\frac{d}{2},\,\frac{d}{p}})}\Big)\\
        \lesssim&\,(1+\mathcal{X}_p(t) )^{\frac{n}{2}+1}\mathcal{X}_p^2(t),
    \end{aligned}
\end{equation}

 \begin{equation}\label{eq:9}
    \begin{aligned}
        &\|u^+\cdot\nabla u^+ \|^{\ell}_{L^1_t(\dot{B}_{2,1}^{\frac{d}{2}-1})}\\
        \lesssim &  \|u^+\|_{\widetilde{L}^{\infty}_t(\dot{\mathcal{B}}_{2,\,p}^{\frac{d}{2}-1,\,\frac{d}{p}-1})} \|\nabla u^+ \|_{L^1_t(\dot{\mathcal{B}}_{2,\,p}^{\frac{d}{2},\, \frac{d}{p}})}+\|\nabla u^+\|_{\widetilde{L}^{2}_t(\dot{\mathcal{B}}_{2,\,p}^{\frac{d}{2}-1,\,\frac{d}{p}-1})} \|u^+ \|_{\widetilde{L}^2_t(\dot{\mathcal{B}}_{2,\,p}^{\frac{d}{2},\,\frac{d}{p}})}
    \end{aligned}
\end{equation}

and 
\begin{equation}\label{eq:10}
    \begin{aligned}
       &\|h^+(n^+,n^-)\nabla n^+\nabla u^+\|^{\ell}_{L^1_t(\dot{B}_{2,1}^{\frac{d}{2}-1})}\\
        \lesssim &\,\|h^+(n^+,n^-)\|_{\widetilde{L}^{\infty}_t(\dot{\mathcal{B}}_{2,\,p}^{\frac{d}{2}-1,\,\frac{d}{p}-1})}\|\nabla n^+\nabla u^+\|_{L^{1}_t(\dot{\mathcal{B}}_{2,\,p}^{\frac{d}{2},\,\frac{d}{p}-1})}+\|\nabla n^+\nabla u^+\|_{L^1_t(\dot{\mathcal{B}}_{2,\,p}^{\frac{d}{2}-1,\,\frac{d}{p}-1})}\|h^+(n^+,n^-)\|_{\widetilde{L}^{\infty}_t(\dot{\mathcal{B}}_{2,\,p}^{\frac{d}{2},\,\frac{d}{p}})}\\
        \lesssim&\,\|h^+(n^+,n^-)\|_{\widetilde{L}^{\infty}_t(\dot{\mathcal{B}}_{2,\,p}^{\frac{d}{2}-1,\,\frac{d}{p}})}\|\nabla n^+\nabla u^+\|_{L^1_t(\dot{\mathcal{B}}_{2,\,p}^{\frac{d}{2}-1,\,\frac{d}{p}-1})} .
    \end{aligned}
\end{equation}
Next we need to bound $\|\nabla n^+\nabla u^+\|_{L^1_t(\dot{\mathcal{B}}_{2,\,p}^{\frac{d}{2}-1,\,\frac{d}{p}-1})}$. Due to  \eqref{A1} and $\eqref{A2}$, we have
\begin{equation}\label{eq:11}
    \begin{aligned}
        &\|\nabla n^+\nabla u^+\|^{\ell}_{L^1_t(\dot{\mathcal{B}}_{2,1}^{\frac{d}{2}-1 })}\\
        \lesssim&\, \|\nabla u^+\|_{\widetilde{L}^{2}_t(\dot{\mathcal{B}}_{2,\,p}^{\frac{d}{2}-1,\,\frac{d}{p}-1})}\|\nabla n^+ \|_{\widetilde{L}^2_t(\dot{\mathcal{B}}_{2,\,p}^{\frac{d}{2},\,\frac{d}{p}-1})}+\|\nabla n^+\|_{\widetilde{L}^{\infty}_t(\dot{\mathcal{B}}_{2,\,p}^{\frac{d}{2}-1,\,\frac{d}{p}-1})}\|\nabla u^+\|_{L^{1}_t(\dot{\mathcal{B}}_{2,\,p}^{\frac{d}{2},\,\frac{d}{p}})},
    \end{aligned}
\end{equation}
\begin{equation}\label{eq:12}
    \begin{aligned}
         \|\nabla n^+\nabla u^+\|^{h}_{L^1_t(\dot{\mathcal{B}}_{p,1}^{\frac{d}{p}-1 })}\lesssim\, \|\nabla n^+\|_{\widetilde{L}^{\infty}_t(\dot{\mathcal{B}}_{2,\,p}^{\frac{d}{2}-1,\,\frac{d}{p}-1})}  \|\nabla u^+\|_{L^1_t(\dot{\mathcal{B}}_{2,\,p}^{\frac{d}{2},\,\frac{d}{p}})}.
    \end{aligned}
\end{equation}
Inserting \eqref{eq:11} and \eqref{eq:12} into \eqref{eq:10}, with aid of Proposition \ref{propA2}, we can get
\begin{equation}\label{eq:13}
    \begin{aligned}
        &\|h^+(n^+,n^-)\nabla n^+\nabla u^+\|^{\ell}_{L^1_t(\dot{B}_{2,1}^{\frac{d}{2}-1})}\\
        \lesssim&\;\|h^+(n^+,n^-)\|_{\widetilde{L}^{\infty}_t(\dot{\mathcal{B}}_{2,\,p}^{\frac{d}{2}-1,\,\frac{d}{p}})}\Big(\| u^+\|_{\widetilde{L}^{2}_t(\dot{\mathcal{B}}_{2,\,p}^{\frac{d}{2},\,\frac{d}{p}})}\| n^+ \|_{\widetilde{L}^2_t(\dot{\mathcal{B}}_{2,\,p}^{\frac{d}{2},\,\frac{d}{p}})} +\| n^+\|_{\widetilde{L}^{\infty}_t(\dot{\mathcal{B}}_{2,\,p}^{\frac{d}{2},\,\frac{d}{p}})}\|\nabla u^+\|_{L^{1}_t(\dot{\mathcal{B}}_{2,\,p}^{\frac{d}{2},\,\frac{d}{p}})}\Big)\\
        \lesssim&\; (1+\mathcal{X}_p(t) )^{\frac{n}{2}+1}\mathcal{X}_p^2(t).
    \end{aligned}
\end{equation}
For the last term $l^+(n^+,n^-)\Delta u^+$,  choosing $s_1 = s_2 = \frac{n}{2}-1$, $t_1 = t_2 =\frac{n}{2}$, $r_1 = r_4 = \infty$, $r_2 = r_3 =
1$, $\gamma=-1$ in \eqref{A2} of Proposition \ref{propA1}, we get
\begin{equation}\label{eq:14}
    \begin{aligned}
        &\|l^+(n^+,n^-)\Delta u^+\|^{\ell}_{L^1_t(\dot{B}_{2,1}^{\frac{d}{2}-1})}\\
        \lesssim&\,\|l^+(n^+,n^-)\|_{\widetilde{L}^{\infty}_t(\dot{\mathcal{B}}_{2,\,p}^{\frac{d}{2}-1,\,\frac{d}{p}-1})}   \|\Delta u^+\|_{L^{1}_t(\dot{\mathcal{B}}_{2,\,p}^{\frac{d}{2},\,\frac{d}{p}-1})}+\|\Delta u^+\|_{\widetilde{L}^{1}_t(\dot{\mathcal{B}}_{2,\,p}^{\frac{d}{2}-1,\,\frac{d}{p}-1})}   \|l^+(n^+,n^-)\|_{\widetilde{L}^{\infty}_t(\dot{\mathcal{B}}_{2,\,p}^{\frac{d}{2},\,\frac{d}{p}})}\\
        \lesssim&\, (1+\mathcal{X}_p(t) )\mathcal{X}_p^2(t).
    \end{aligned}
\end{equation}
Collecting all estimates \eqref{eq:8}-\eqref{eq:14} together, we have 
\begin{equation}
    \begin{aligned}
        \|\mathcal{F}_2^+\|^{\ell}_{L^1_t(\dot{B}_{2,1}^{\frac{d}{2}-1})}\lesssim\,(1+\mathcal{X}_p(t) )^{\frac{n}{2}+1}\mathcal{X}_p^2(t).
    \end{aligned}
\end{equation}
The terms $ \mathcal{F}_2^-$ can be treated along the same way as $\mathcal{F}_2^+$ above.
Thus one  combine these estimates above to get
\begin{equation}
    \begin{aligned}
        &\| (n^+,n^-,u^+,u^-) \|^{\ell}_{\widetilde{L}^{\infty}_t(\dot{B}_{2,1}^{\frac{d}{2}-1})\cap L^{1}_t(\dot{B}_{2,1}^{\frac{d}{2}+1}) }\lesssim\| (n_0^+,n_0^-,u_0^+,u_0^-) \|^{\ell}_{ \dot{B}_{2,1}^{\frac{d}{2}-1}}+(1+\mathcal{X}_p(t) )^{\frac{n}{2}+1}\mathcal{X}_p^2(t).
    \end{aligned}
\end{equation}

\end{proof}

\subsection{High-frequency estimate}
The main difficulty in the high-frequency analysis comes from the fact that the density variables $n^\pm$ do not possess any direct parabolic smoothing effect. A direct estimate of the transport equations for $n^\pm$ would lead to a loss of one derivative, due to the convection terms $u^\pm\cdot\nabla n^\pm$ in $\mathcal{F}_1^+,\mathcal{F}_1^-$. To overcome this difficulty, we introduce two effective velocities as follows
\begin{equation}\label{effective-velo-w}
\left\{
    \begin{aligned}
        w^{+}=\nabla (-\Delta)^{-1}\left( \frac{\gamma_1}{\nu_1^++\nu_2^+}n^+ +\frac{\gamma_2}{\nu_1^++\nu_2^+}n^--\text{div}u^+\right),\\
          w^{-}=\nabla (-\Delta)^{-1}\left( \frac{\gamma_3}{\nu_1^-+\nu_2^-}n^+ +\frac{\gamma_4}{\nu_1^-+\nu_2^-}n^--\text{div}u^-\right).
    \end{aligned}
    \right.
\end{equation}
The idea is to combine the pressure-gradient terms with the viscous part of the momentum equations so that the new unknowns satisfy heat-type equations.
Applying $\text{div}$ to \eqref{effective-velo-w}, we obtain that
\begin{equation}
\left\{
    \begin{aligned}
        \text{div}u^+-\text{div}w^+=\frac{\gamma_1}{\nu_1^++\nu_2^+}n^+ +\frac{\gamma_2}{\nu_1^++\nu_2^+}n^-,\\
         \text{div}u^--\text{div}w^-=\frac{\gamma_3}{\nu_1^-+\nu_2^-}n^+ +\frac{\gamma_4}{\nu_1^-+\nu_2^-}n^-.
    \end{aligned}
    \right.
\end{equation}
Due to the structural positivity condition
$\gamma_1\gamma_4-\gamma_2\gamma_3>0,$
which follows from the stability assumption $f'(1)<0$.By solving this system for $(n^+,n^-)$, we find that there exist two linear function $\Phi, \Psi$ such that  
\begin{equation}
\left\{
    \begin{aligned}
        n^+= \Phi(\text{div}w^+,\text{div}w^-,\text{div}u^+,\text{div}u^-),\\
         n^-= \Psi(\text{div}w^+,\text{div}w^-,\text{div}u^+,\text{div}u^-).
    \end{aligned}
    \right.
\end{equation}
This avoids the derivative loss caused by the transport terms in the density equations.

\begin{lemma}\label{Lem2.2}
Let  $T>0$ be any given time, and $(n^+,n^-,u^+,u^-) $ be any strong solution to the Cauchy problem \eqref{equ:fornu-pm} for $t\in(0,T)$. Then it holds that
    \begin{equation}
    \begin{aligned}
          & \|(n^+,n^-)\|^h_{\widetilde{L}^{\infty}_t(\dot{B}_{p,1}^{\frac{d}{p} })\cap  L^{1}_t(\dot{B}_{p,1}^{\frac{d}{p}})} +\|(u^+,u^-)\|^h_{\widetilde{L}^{\infty}_t(\dot{B}_{p,1}^{\frac{d}{p}-1})\cap  L^{1}_t(\dot{B}_{p,1}^{\frac{d}{p}+1})}  \\
         \lesssim&\,\|(n_0^+,n_0^-)\|^h_{\dot{B}_{p,1}^{\frac{d}{p} }}+\|(u^+_0,u^-_0)\|^h_{\dot{B}_{p,1}^{\frac{d}{p}-1}}+  (1+\mathcal{X}_p(t) +\mathcal{X}_p^2(t)) \mathcal{X}_p^2(t).
    \end{aligned}
\end{equation}
\end{lemma}
\begin{proof}
It follows from \eqref{equ:fornu-pm} that
\begin{equation} 
 \left\{
   \begin{array}{ll}
     \partial_t w^+-(\nu_1^++\nu_2^+)\Delta w^+  =g_1,\\
     \partial_t w^- -(\nu_1^-+\nu_2^-)\Delta w^-  =g_2,
   \end{array}
 \right.
\end{equation}
where
\begin{equation}\label{def:g1g2}
\left\{
\begin{aligned}
     g_1=&\nabla (-\Delta)^{-1}\left( \frac{\gamma_1}{\nu_1^++\nu_2^+}\mathcal{F}_1^+ +\frac{\gamma_2}{\nu_1^++\nu_2^+}\mathcal{F}_1^--\operatorname{div}\mathcal{F}_2^+\right)\\
     &+\frac{\gamma_1^2}{\nu_1^++\nu_2^+}w^++\frac{\gamma_2\gamma_4}{\nu_1^++\nu_2^+}w^- \\
      &-\frac{\gamma_1^2}{\nu_1^++\nu_2^+}\nabla (-\Delta)^{-1}\left( \frac{\gamma_1}{\nu_1^++\nu_2^+}n^+ +\frac{\gamma_2}{\nu_1^++\nu_2^+}n^- \right)\\
      &-\frac{\gamma_2\gamma_4}{\nu_1^++\nu_2^+}\nabla (-\Delta)^{-1}\left( \frac{\gamma_3}{\nu_1^-+\nu_2^-}n^+ +\frac{\gamma_4}{\nu_1^-+\nu_2^-}n^- \right),\\
    g_2=& \nabla (-\Delta)^{-1}\left( \frac{\gamma_3}{\nu_1^++\nu_2^+}\mathcal{F}_1^+ +\frac{\gamma_4}{\nu_1^-+\nu_2^-}\mathcal{F}_1^--\operatorname{div}\mathcal{F}_2^-\right)\\
    &+\frac{\gamma_1\gamma_3}{\nu_1^-+\nu_2^-}w^++\frac{\gamma_4^2}{\nu_1^-+\nu_2^-}w^-\\
      &-\frac{\gamma_1\gamma_3}{\nu_1^-+\nu_2^-}\nabla (-\Delta)^{-1}\left( \frac{\gamma_1}{\nu_1^++\nu_2^+}n^+ +\frac{\gamma_2}{\nu_1^++\nu_2^+}n^- \right)\\
      &-\frac{\gamma_4^2}{\nu_1^-+\nu_2^-}\nabla (-\Delta)^{-1}\left( \frac{\gamma_3}{\nu_1^-+\nu_2^-}n^+ +\frac{\gamma_4}{\nu_1^-+\nu_2^-}n^- \right).
\end{aligned}
\right.
\end{equation}
Applying Proposition \ref{lame-system-estimate} to the above equations implies that
\begin{equation}\label{ineq:wpm}
    \begin{aligned}
        &\|(w^+,w^-)\|^h_{\widetilde{L}^{\infty}_t(\dot{B}_{p,1}^{\frac{d}{p}-1})\cap  L^{1}_t(\dot{B}_{p,1}^{\frac{d}{p}+1})}\lesssim \|(w^+_0,w^-_0)\|^h_{\dot{B}_{p,1}^{\frac{d}{p}-1}}+\|(w^+,w^-)\|^h_{L^{1}_t(\dot{B}_{p,1}^{\frac{d}{p}-1})}\\
       +&\|(n^+,n^-)\|^h_{L^{1}_t(\dot{B}_{p,1}^{\frac{d}{p}-2})}+\|(\mathcal{F}_1^+,\mathcal{F}_1^-)\|^h_{L^{1}_t(\dot{B}_{p,1}^{\frac{d}{p}-2})}+\|(\mathcal{F}_2^+,\mathcal{F}_2^-)\|^h_{L^{1}_t(\dot{B}_{p,1}^{\frac{d}{p}-1})}.
    \end{aligned}
\end{equation}
Secondly, we see that $(n^+,n^-)$ satisfies the following damped equations in terms of effective velocities
\begin{equation}
\left\{
    \begin{aligned}
      \partial_t n^++u^+\cdot\nabla n^+ +   \frac{\gamma_1}{\nu_1^++\nu_2^+}n^+ +\frac{\gamma_2}{\nu_1^++\nu_2^+}n^-=-n^+\text{div}u^+-\text{div}w^+,\\
        \partial_t n^-+u^-\cdot\nabla n^-+  \frac{\gamma_3}{\nu_1^-+\nu_2^-}n^+ +\frac{\gamma_4}{\nu_1^-+\nu_2^-}n^-=-n^-\text{div}u^--\text{div}w^-.
    \end{aligned}
    \right.
\end{equation}
Applying $\dot{\Delta}_j$ to the above equations, we obtain
\begin{equation}\label{damping-for-npm-local}
\left\{
    \begin{aligned}
      \partial_t n_j^++u^+\cdot\nabla n_j^+ +   \frac{\gamma_1}{\nu_1^++\nu_2^+}n_j^+ +\frac{\gamma_2}{\nu_1^++\nu_2^+}n_j^-=-\dot{\Delta}_j(n^+\text{div}u^+)-\text{div}w_j^++R_j^+,\\
        \partial_t n_j^-+u^-\cdot\nabla n_j^-+  \frac{\gamma_3}{\nu_1^-+\nu_2^-}n_j^+ +\frac{\gamma_4}{\nu_1^-+\nu_2^-}n_j^-=-\dot{\Delta}_j(n^-\text{div}u^-)-\text{div}w_j^-+R_j^-,
    \end{aligned}
    \right.
\end{equation}
where $R_j^\pm=[u^\pm\cdot\nabla,\dot{\Delta}_j]n^\pm$. Multiplying $(\ref{damping-for-npm-local})_1$  by $\frac{1}{\gamma_2}n_j^+|n_j^+|^{p-2}$ and $(\ref{damping-for-npm-local})_2$  by $\frac{1}{\gamma_3}n_j^-|n_j^-|^{p-2}$, and then integrating over  $\mathbb{R}^d\times[0,t]$, we can obtain
%\begin{equation}
  %  \begin{aligned}
  %      &\frac{1}{p\gamma_2}\frac{d}{dt}\|n_j^+\|^p_{L^p}+\frac{\gamma_+^2}{\nu_1^++\nu_2^+}\|n_j^+\|^p_{L^p}\leq \frac{1}{\gamma_2}\left( \frac{\gamma_2}{\nu_1^++\nu_2^+}\|n_j^-\|_{L^p} \right.\\
 %      +&\left.\frac{1}{p}\|\text{div}u^+\|_{L^{\infty}}\|n_j^+\|_{L^p}+\|\dot{\Delta}_j(n^+\text{div}u^+)\|_{L^p}+\|\text{div}w_j^+\|_{L^p}+\|R_j^{+}\|_{L^p}\right)\|n_j^+\|^{p-1}_{L^p}
  %  \end{aligned}
%\end{equation}
\begin{equation}\label{eq:n+}
    \begin{aligned}
        &\frac{1}{p\gamma_2}\|n_j^+\|_{L^p}+\frac{\gamma_+^2}{\nu_1^++\nu_2^+}\int_0^t\|n_j^+\|_{L^p}ds\leq \frac{1}{p\gamma_2}\|\dot{\Delta}_j n_{0}^+\|_{L^p} +\frac{1}{\nu_1^++\nu_2^+}\int_0^t\|n_j^-\|_{L^p}ds \\
       +&\frac{1}{\gamma_2}\left(\int_0^t \frac{1}{p}\|\nabla u^+\|_{L^{\infty}}\|n_j^+\|_{L^p}ds+\int_0^t \|\dot{\Delta}_j(n^+\text{div}u^+)\|_{L^p}+\|\text{div}w_j^+\|_{L^p}ds+\int_0^t \|R_j^+\|_{L^p}ds\right)
    \end{aligned}
\end{equation}
and
\begin{equation}\label{eq:n-}
    \begin{aligned}
        &\frac{1}{p\gamma_3}\|n_j^-\|_{L^p}+\frac{\gamma_-^2}{\nu_1^-+\nu_2^-}\int_0^t\|n_j^-\|_{L^p}ds\leq \frac{1}{p\gamma_3} \|\dot{\Delta}_j n_{0}^-\|_{L^p} +\frac{1}{\nu_1^-+\nu_2^-}\int_0^t\|n_j^+\|_{L^p}ds \\
       +&\frac{1}{\gamma_3}\left(\int_0^t \frac{1}{p}\|\nabla u^-\|_{L^{\infty}}\|n_j^+\|_{L^p}ds+\int_0^t \|\dot{\Delta}_j(n^-\text{div}u^-)\|_{L^p}+\|\text{div}w_j^+\|_{L^p}ds+\int_0^t \|R_j^-\|_{L^p}ds\right).
    \end{aligned}
\end{equation}
From \eqref{condition:gamma}, we notice that $\gamma_1\gamma_4>\gamma_2\gamma_3$. And this implies that $\gamma_+^2\gamma_-^2>1$, then there exists a constant $\tilde{\eta}\in\left(\frac{\nu_1^-+\nu_2^-}{(\nu_1^++\nu_2^+)\gamma_-^2},\frac{\nu_1^-+\nu_2^-}{\nu_1^++\nu_2^+}\gamma_+^2\right)$, such that 
\begin{equation}\label{eq:npm}
    \begin{aligned}
        &\frac{1}{p\gamma_2}\|n_j^+\|_{L^p}+\frac{1}{p\gamma_3}\|n_j^-\|_{L^p}+\left(\frac{\gamma_+^2}{\nu_1^++\nu_2^+}-\frac{\tilde{\eta}}{\nu_1^-+\nu_2^-}\right)\int_0^t\|n_j^+\|_{L^p}ds\\
        +&\left(\frac{\tilde{\eta}\gamma_-^2}{\nu_1^-+\nu_2^-}-\frac{1}{\nu_1^++\nu_2^+}\right)\int_0^t\|n_j^-\|_{L^p}ds\lesssim \|\dot{\Delta}_j n_{0}^+\|_{L^p} +   \|\dot{\Delta}_j n_{0}^-\|_{L^p} \\
        +&\int_0^t \frac{1}{p}\|\nabla u^+\|_{L^{\infty}}\|n_j^+\|_{L^p}ds+\int_0^t \|\dot{\Delta}_j(n^+\text{div}u^+)\|_{L^p}+\|\text{div}w^+\|_{L^p}ds+\int_0^t \|R_j^+\|_{L^p}ds\\
        +&\int_0^t \frac{1}{p}\|\nabla u^-\|_{L^{\infty}}\|n_j^+\|_{L^p}ds+\int_0^t \|\dot{\Delta}_j(n^-\text{div}u^-)\|_{L^p}+\|\text{div}w^+\|_{L^p}ds+\int_0^t \|R_j^-\|_{L^p}ds.
    \end{aligned}
\end{equation}
Now multiplying \eqref{eq:npm}   by $2^{j\frac{d}{p}}$, and then summing over the index $j\geq j_0$, we have
\begin{equation}\label{ineq:npm}
    \begin{aligned}
        & \|(n^+,n^-)\|^h_{\widetilde{L}^{\infty}_t(\dot{B}_{p,1}^{\frac{d}{p} })\cap  L^{1}_t(\dot{B}_{p,1}^{\frac{d}{p}})}\\
         \lesssim&\, \|(n_0^+,n_0^-)\|^h_{\dot{B}_{p,1}^{\frac{d}{p} }} +\|\nabla u^+\|_{L^{1}_t(\dot{B}_{p,1}^{\frac{d}{p}})}\|n^+\|_{ L^{1}_t(\dot{B}_{p,1}^{\frac{d}{p}})} +\|w^+\|^h_{L^{1}_t(\dot{B}_{p,1}^{\frac{d}{p}+1})}\\
         \lesssim&\, \|(n_0^+,n_0^-)\|^h_{\dot{B}_{p,1}^{\frac{d}{p} }}+\|(w^+,w^-)\|^h_{L^{1}_t(\dot{B}_{p,1}^{\frac{d}{p}+1})}+\mathcal{X}_p^{2}(t).
    \end{aligned}
\end{equation}
Multiplying \eqref{ineq:npm} by $\eta_2>0$, and then add this inequality to \eqref{ineq:wpm} together, we can get
\begin{equation}
    \begin{aligned}
          &\|(n^+,n^-)\|^h_{\widetilde{L}^{\infty}_t(\dot{B}_{p,1}^{\frac{d}{p} })\cap  L^{1}_t(\dot{B}_{p,1}^{\frac{d}{p}})}+\|(w^+,w^-)\|^h_{\widetilde{L}^{\infty}_t(\dot{B}_{p,1}^{\frac{d}{p}-1})\cap  L^{1}_t(\dot{B}_{p,1}^{\frac{d}{p}+1})}\\
          \lesssim&\,\|(n_0^+,n_0^-)\|^h_{\dot{B}_{p,1}^{\frac{d}{p} }}+\|(w^+_0,w_0^-)\|^h_{\dot{B}_{p,1}^{\frac{d}{p}-1}}+\eta_2\|(w^+,w^-)\|^h_{L^{1}_t(\dot{B}_{p,1}^{\frac{d}{p}+1})}\\
         &+\|(\mathcal{F}_1^+,\mathcal{F}_1^-)\|^h_{L^{1}_t(\dot{B}_{p,1}^{\frac{d}{p}-2})}+\|(\mathcal{F}_2^+,\mathcal{F}_2^-)\|^h_{L^{1}_t(\dot{B}_{p,1}^{\frac{d}{p}-1})} +\mathcal{X}_p^{2}(t).
    \end{aligned}
\end{equation}
By choosing $\eta_2$ sufficiently small, we can get 
\begin{equation}
    \begin{aligned}
        & \|(n^+,n^-)\|^h_{\widetilde{L}^{\infty}_t(\dot{B}_{p,1}^{\frac{d}{p} })\cap  L^{1}_t(\dot{B}_{p,1}^{\frac{d}{p}})} +\|(w^+,w^-)\|^h_{\widetilde{L}^{\infty}_t(\dot{B}_{p,1}^{\frac{d}{p}-1})\cap  L^{1}_t(\dot{B}_{p,1}^{\frac{d}{p}+1})}  \\
         \lesssim&\,\|(n_0^+,n_0^-)\|^h_{\dot{B}_{p,1}^{\frac{d}{p} }}+\|(u^+_0,u^-_0)\|^h_{\dot{B}_{p,1}^{\frac{d}{p}-1}}+\|(\mathcal{F}_1^+,\mathcal{F}_1^-)\|^h_{L^{1}_t(\dot{B}_{p,1}^{\frac{d}{p}-2})}+\|(\mathcal{F}_2^+,\mathcal{F}_2^-)\|^h_{L^{1}_t(\dot{B}_{p,1}^{\frac{d}{p}-1})}+\mathcal{X}_p^{2}(t).
    \end{aligned}
\end{equation}

Notice that 
\begin{equation}
\left\{
    \begin{aligned}
    w^{+}=\nabla (-\Delta)^{-1}\left( \frac{\gamma_1}{\nu_1^++\nu_2^+}n^+ +\frac{\gamma_2}{\nu_1^++\nu_2^+}n^--\text{div}u^+\right),\\
    w^{-}=\nabla (-\Delta)^{-1}\left( \frac{\gamma_3}{\nu_1^-+\nu_2^-}n^+ +\frac{\gamma_4}{\nu_1^-+\nu_2^-}n^--\text{div}u^-\right),
\end{aligned}
\right.
\end{equation}
and 
\begin{align}
    \|\nabla (-\Delta)^{-1}n^\pm\|_{\widetilde{L}^{\infty}_t(\dot{B}_{p,1}^{\frac{d}{p}-1})\cap  L^{1}_t(\dot{B}_{p,1}^{\frac{d}{p}+1})}\lesssim 2^{-j}\|n^\pm\|_{\widetilde{L}^{\infty}_t(\dot{B}_{p,1}^{\frac{d}{p}-1})\cap  L^{1}_t(\dot{B}_{p,1}^{\frac{d}{p}+1})},
\end{align}
we arrive at, for $j\geq j_0$,
\begin{equation}
    \begin{aligned}
        & \|(n^+,n^-)\|^h_{\widetilde{L}^{\infty}_t(\dot{B}_{p,1}^{\frac{d}{p} })\cap  L^{1}_t(\dot{B}_{p,1}^{\frac{d}{p}})} +\|(u^+,u^-)\|^h_{\widetilde{L}^{\infty}_t(\dot{B}_{p,1}^{\frac{d}{p}-1})\cap  L^{1}_t(\dot{B}_{p,1}^{\frac{d}{p}+1})}  \\
         \lesssim&\,\|(n_0^+,n_0^-)\|^h_{\dot{B}_{p,1}^{\frac{d}{p} }}+\|(u^+_0,u^-_0)\|^h_{\dot{B}_{p,1}^{\frac{d}{p}-1}}+\|(\mathcal{F}_1^+,\mathcal{F}_1^-)\|^h_{L^{1}_t(\dot{B}_{p,1}^{\frac{d}{p}-2})}+\|(\mathcal{F}_2^+,\mathcal{F}_2^-)\|^h_{L^{1}_t(\dot{B}_{p,1}^{\frac{d}{p}-1})}+\mathcal{X}_p^{2}(t).
    \end{aligned}
\end{equation}
Likely, we need to bound those different terms in $\mathcal{F}_1^+, \mathcal{F}_2^+, \mathcal{F}_1^-, \mathcal{F}_2^-$. With regard to $u^\pm\cdot\nabla n^\pm$ and $n^\pm\operatorname{div}u^\pm$,  by taking  $r_1=r_2=2$, $\sigma=\frac{n}{p}$, $\tau=\frac{n}{p}-1$ in \eqref{A1} of Proposition \ref{propA1}, we arrive at
\begin{equation}\label{eq:15}
    \begin{aligned}
        &\|u^+\cdot\nabla n^+\|^h_{L^{1}_t(\dot{B}_{p,1}^{\frac{d}{p}-2})}
        \lesssim\,\|u^+\|_{\widetilde{L}^{2}_t(\dot{\mathcal{B}}_{2,p}^{\frac{d}{2},\,\frac{d}{p} })}\|\nabla n^+\| _{\widetilde{L}^{2}_t(\dot{\mathcal{B}}_{2,p}^{\frac{d}{2}-1,\,\frac{d}{p}-1})} 
        \lesssim\, \mathcal{X}_p^{2}(t)
    \end{aligned}
\end{equation}
and  
\begin{equation}\label{eq:16}
    \begin{aligned}
        &\|n^+\operatorname{div}u^+\|^h_{L^{1}_t(\dot{B}_{p,1}^{\frac{d}{p}-2})}
        \lesssim\,\|n^+\|_{\widetilde{L}^{2}_t(\dot{\mathcal{B}}_{2,\,p}^{\frac{d}{2},\,\frac{d}{p} })}\|\operatorname{div}u^+\| _{\widetilde{L}^{2}_t(\dot{\mathcal{B}}_{2,\,p}^{\frac{d}{2}-1,\,\frac{d}{p}-1})} 
        \lesssim\, \mathcal{X}_p^{2}(t).
    \end{aligned}
\end{equation}
Then combining \eqref{eq:15} and \eqref{eq:16}, we have 
\begin{equation}
    \begin{aligned}
   \|(\mathcal{F}_1^+,\mathcal{F}_1^-)\|^h_{L^{1}_t(\dot{B}_{p,1}^{\frac{d}{p}-2})}\lesssim\, \mathcal{X}_p^{2}(t).
\end{aligned}
\end{equation}

With regard to $\mathcal{F}_2^+\sim g^+(n^+,n^-)\,\nabla n^+ + (u^+\cdot\nabla)u^+
      + h^+(n^+,n^-)\,\nabla n^+\nabla u^+
      + l^+(n^+,n^-)\,\Delta u^+ $. For $ g^+(n^+,n^-)\,\nabla n^+,$ by taking  $r_1=r_2=2$, $\sigma=\frac{d}{p}$, $\tau=\frac{d}{p}-1$ in \eqref{A1} of Proposition \ref{propA1}, we have
\begin{equation}\label{eq:17}
    \begin{aligned}
        &\|g^+(n^+,n^-)\nabla n^+\|^{h}_{L^1_t(\dot{B}_{p,1}^{\frac{d}{p}-1})}
        \lesssim\,\|g^+(n^+,n^-)\|_{\widetilde{L}^{2}_t(\dot{\mathcal{B}}_{2,\,p}^{\frac{d}{2} ,\,\frac{d}{p}})}\|\nabla n^+\| _{\widetilde{L}^{2}_t(\dot{\mathcal{B}}_{2,\,p}^{\frac{d}{2}-1,\,\frac{d}{p}-1})}
         \\
        \lesssim&\,(1+\mathcal{X}_p(t))^{\frac{n}{2}+1}\mathcal{X}_p^2(t).
    \end{aligned}
\end{equation}
For $u^+\cdot\nabla u^+$, taking  $r_1=\infty, r_2=1$, $\sigma=\frac{d}{p}-1$, $\tau=\frac{d}{p}$ in \eqref{A1} of Proposition \ref{propA1}, we have
\begin{equation}
    \begin{aligned}\label{eq:18}
        \|u^+\cdot\nabla u^+\|^{h}_{L^1_t(\dot{B}_{p,1}^{\frac{d}{p}-1})}
        \lesssim\|u^+\|_{\widetilde{L}^{\infty}_t(\dot{\mathcal{B}}_{2,\,p}^{\frac{d}{2}-1,\,\frac{d}{p}-1})}   \|\nabla u^+\|_{L^1_t(\dot{\mathcal{B}}_{2,\,p}^{\frac{d}{2},\,\frac{d}{p}})}\lesssim \mathcal{X}_p^2(t).
    \end{aligned}
\end{equation}
For $h^+(n^+,n^-)\nabla n^+\nabla u^+$,   taking  $r_1=\infty, r_2=1$, $\sigma=\frac{d}{p}$, $\tau=\frac{d}{p}-1$ in \eqref{A1} of Proposition \ref{propA1},  and  using same process as  in \eqref{eq:11}-\eqref{eq:12}, we have
\begin{equation}\label{eq:19}
    \begin{aligned}
        \|h^+(n^+,n^-)\nabla n^+\nabla u^+\|^{h}_{L^1_t(\dot{B}_{p,1}^{\frac{d}{p}-1})}
        \lesssim&\,\|h^+(n^+,n^-)\|_{\widetilde{L}^{\infty}_t(\dot{\mathcal{B}}_{2,\,p}^{\frac{d}{2},\,\frac{d}{p}})}\|\nabla n^+\nabla u^+\|_{L^1_t(\dot{\mathcal{B}}_{2,\,p}^{\frac{d}{2}-1,\,\frac{d}{p}-1})}\\
        \lesssim&\,(1+\mathcal{X}_p(t))^{\frac{d}{2}+1}\mathcal{X}_p^2(t).
    \end{aligned}
\end{equation}
For the last term $l^+(n^+,n^-)\Delta u^+$, taking  $r_1=\infty, r_2=1$, $\sigma=\frac{d}{p}$, $\tau=\frac{d}{p}-1$ in \eqref{A1} of Proposition \ref{propA1},
\begin{equation}\label{eq:20}
    \begin{aligned}
        \|l^+(n^+,n^-)\Delta u^+\|^{h}_{L^1_t(\dot{B}_{p,1}^{\frac{d}{p}-1})}
        \lesssim&\,\|l^+(n^+,n^-)\|_{\widetilde{L}^{\infty}_t(\dot{\mathcal{B}}_{2,\,p}^{\frac{d}{2},\,\frac{d}{p}})}   \|\Delta u^+\|_{L^1_t(\dot{\mathcal{B}}_{2,\,p}^{\frac{d}{2}-1,\,\frac{d}{p}-1})}\\
        \lesssim &\,(1+\mathcal{X}_p(t))^{\frac{d}{2}+1}\mathcal{X}_p^2(t).
    \end{aligned}
\end{equation}
Thus combining \eqref{eq:17}-\eqref{eq:20}, we get 
\begin{equation}\label{eq:21}
    \begin{aligned}
        \|\mathcal{F}_2^+\|^{h}_{L^1_t(\dot{B}_{p,1}^{\frac{d}{p}-1})}
        \lesssim \,(1+\mathcal{X}_p(t))^{\frac{d}{2}+1}\mathcal{X}_p^2(t).
    \end{aligned}
\end{equation}
The term $ \mathcal{F}_2^-$ can  be treated along the same way as $ \mathcal{F}_2^+$ above.
Thus one  combine these estimates above to get
\begin{equation}\label{eq:22}
    \begin{aligned}
          & \|(n^+,n^-)\|^h_{\widetilde{L}^{\infty}_t(\dot{B}_{p,1}^{\frac{d}{p} })\cap  L^{1}_t(\dot{B}_{p,1}^{\frac{d}{p}})} +\|(u^+,u^-)\|^h_{\widetilde{L}^{\infty}_t(\dot{B}_{p,1}^{\frac{d}{p}-1})\cap  L^{1}_t(\dot{B}_{p,1}^{\frac{d}{p}+1})}  \\
         \lesssim&\,\|(n_0^+,n_0^-)\|^h_{\dot{B}_{p,1}^{\frac{d}{p} }}+\|(u^+_0,u^-_0)\|^h_{\dot{B}_{p,1}^{\frac{d}{p}-1}}+(1+\mathcal{X}_p(t))^{\frac{d}{2}+1}\mathcal{X}_p^2(t).
    \end{aligned}
\end{equation}

\end{proof}

\subsection{Proof of Theorem \ref{Thm1}}
\subsubsection{Existence}

To prove the global existence, one can construct a local Friedrichs approximation (cf. \cite{bahouri1}),  extend the local approximate sequence to a global one by the uniform-in-time a priori estimates  
established in Proposition \ref{Prop2.1}, and then show the convergence of the approximate sequence to the expected global solution to the Cauchy problem for \eqref{equ:fornu-pm-S}. More precisely, consider the approximate problems ($q\geq1$):
\begin{equation} \label{sysm-Friedrich}
\left\{
\begin{aligned}
&\partial_t n^{+}+ \dot{\mathbb{E}}_q(u^+\cdot \nabla n^{+})+\dot{\mathbb{E}}_q ((1+n^{+})\operatorname{div} u^{+})=0, \\
&\partial_t u^+ +  \beta_1 \nabla \dot{\mathbb{E}}_qn^+ +  \beta_2 \dot{\mathbb{E}}_q \nabla n^- - \nu_1^+\Delta \dot{\mathbb{E}}_qu^+ -\nu_2^+\nabla\operatorname{div}\dot{\mathbb{E}}_qu^+  =\dot{\mathbb{E}}_q S_2^+, \\
&\partial_t n^{-}+ \dot{\mathbb{E}}_q(u^-\cdot \nabla n^{-})+\dot{\mathbb{E}}_q ((1+n^{-})\operatorname{div}  u^{-})=0, \\
&\partial_t u^-  + \beta_3  \nabla \dot{\mathbb{E}}_q n^+ +  \beta_4 \nabla\dot{\mathbb{E}}_q n^- -\nu_1^-\Delta \dot{\mathbb{E}}_q u^- -\nu_2^-\nabla\operatorname{div}\dot{\mathbb{E}}_q u^-  = \dot{\mathbb{E}}_qS_2^-,\\
&(\dot{\mathbb{E}}_q n^+,\dot{\mathbb{E}}_qu^+,\dot{\mathbb{E}}_q n^-,\dot{\mathbb{E}}_qu^-)(x,0)=(\dot{\mathbb{E}}_q n_0^+,\dot{\mathbb{E}}_qu_0^+,\dot{\mathbb{E}}_q n_0^-,\dot{\mathbb{E}}_qu^-_0).
\end{aligned}
\right.
\end{equation}
Since \eqref{sysm-Friedrich} is a system of ordinary differential equations  in 
$L^2_q\times L^2_q\times L^2_q\times L^2_q $ and is locally Lipschitz with respect to the variable $( n^+,n^-, u^+,u^-)$ for every $q\geq 1$.  By virtue of the Cauchy-Lipschitz theorem 
(see \cite[Theorem 3.2]{bahouri1}),  there exists a maximal time $T_q>0$ such that the problem \eqref{sysm-Friedrich} admits a unique solution $(n^{+,q},n^{-,q}, u^{+,q},u^{-,q})\in C([0,T_q);L^2_q)$. Then, $(n^{+,q},n^{-,q}, u^{+,q},u^{-,q})$   solves \eqref{sysm-Friedrich} on $[0,T_q)\times \mathbb{R}^d$. 

Define the maximal time
\begin{align}
	T_q:=\sup\{t\geq0~:~& \mathcal{X}_p^q(t)\leq 2C_0\var_0  \},
\end{align}
where $\mathcal{X}_p^q(t)$ is the exact same as \eqref{def-E} but for  $(n^{+,q},n^{-,q}, u^{+,q}, u^{-,q})$, $\mathcal{X}_p(0)$ is defined in \eqref{def-X0} independent of $q$, and $C_0$ is given by Proposition \ref{Prop2.1}. It is clear that $0<T_q\leq T^*_q.$

We claim that $T_q= T^*_q$. Suppose, on the contrary, that $T_q< T^*_q$. Since $(n^{+,q},n^{-,q}, u^{+,q}, u^{-,q})=\dot{\mathbb{E}}_q (n^{+,q},n^{-,q}, u^{+,q}, u^{-,q})$, the orthogonal
projector $\dot{\mathbb{E}}_q$ does not affect the energy estimates for \eqref{sysm-Friedrich} established in Lemmas \ref{Lem2.1}--\ref{Lem2.2}. Moreover, since  $\mathcal{X}_p(0)\leq \var_0$, by choosing $\var_0=\frac{1}{4C_0}$, we obtain $$\|(n^{+,q},n^{-,q})\|_{L^\infty([0,T_q]\times\mathbb{R}^d)}\leq \mathcal{X}^q_p(t)\leq \frac12.$$  Therefore, Proposition \ref{Prop2.1} yields
\begin{align}\label{inequ-for-E1}
	 \mathcal{X}_p^q(t)\lesssim \mathcal{X}_p(0)  +(1+\mathcal{X}_p^q(t))^{\frac{d}{2}+1} \mathcal{X}_p^q(t)^2.
\end{align}
Choosing $\var_0=\min\{\var_0^*,\frac{1}{4C_0}\}$ sufficiently small so that
\begin{align*}
    4(1+2C_0\var_0)^{\frac{d}{2}+1}C_0^2\var_0\leq \frac12,
\end{align*}
we infer that \begin{align}\label{inequ-for-E2}
	\mathcal{X}_p^q(t)\leq C_0\var_0 +4(1+2C_0\var_0)^{\frac{d}{2}+1}C_0^3\var_0^2 \leq \frac{3}{2} C_0\var_0, \quad 0<t<T_q.
\end{align}

By \eqref{inequ-for-E2} and the continuity in time of \(\mathcal{X}_p^q(t)\), it follows that \(T_q\) cannot be the maximal time such that
 $\mathcal{X}(t)\leq 2C_0\var_0$ holds.  This contradicts the definition of $T_q$. Hence $T_q=T_q^*$.

If $T^*_q<\infty$, by \eqref{inequ-for-E2} and $T_q= T^*_q$, we can take $(n^{+,q},n^{-,q}, u^{+,q}, u^{-,q})(t)$ for $t$ sufficiently close to $T^*_q$ as the new initial data and
obtain the existence from $t$ to some $t+\eta^*>T^*_q$ with a sufficiently small constant $\eta^*>0$ by the Cauchy–Lipschitz theorem, which contradicts the definition of $T^*_q$. Therefore, we have $T^*_q=\infty$, and $\{(n^{+,q},n^{-,q}, u^{+,q}, u^{-,q})\}_{q\ge1}$  is a global-in-time approximate sequence.

Then, utilizing the uniform bound $\mathcal{X}^q(t)\leq C_0\mathcal{X}_{0} $ following a standard compactness argument (see \cite[Pages 458-459]{bahouri1}), one can show that as $q\rightarrow \infty$, $(n^{+,q},\theta^q, u^{+,q}, u^{-,q})$ converges strongly to a limit $(n^+,n^-,u^+,u^-)$ which is a global-in-time solution to the Cauchy problem 
\eqref{equ:fornu-pm-S} satisfying $\mathcal{X}(t)\lesssim \mathcal{X}_0$.
 
\subsubsection{Uniqueness} 
For technical reasons, in this paper, we restrict ourselves to the case
$2\leq p\leq d$. 
It is enough to consider the reformulated system \eqref{equ:fornu-pm-S}. Let $T>0$ be given, and let $(n_1^+, n_1^-,u_1^+,u_1^-)$ and  $(n_2^+, n_2^-,u_2^+,u_2^-)$  be two solutions 
 to \eqref{equ:fornu-pm-S} on $[0,T]\times\mathbb{R}^d$ with the same initial data $(n_0^+, n_0^-,u_0^+,u_0^-)$, both satisfying the regularity assumptions stated in Theorem \ref{Thm1}. Then the difference
$$
(\widetilde{n}^+,\widetilde{n}^-,\widetilde{u}^+,\widetilde{u}^-)=(n_1^+-n_2^+, n_1^--n_2^-,u_1^+-u_2^+,u_1^--u_2^-)
$$
solves
\begin{equation}  \label{equ:fornu-pm-error}
\left\{\begin{array}{l}
\partial_{t} \widetilde{n}^{+}+u_2^+\cdot\nabla \widetilde{n}^+ =\widetilde{S}_1^+, \\
\partial_{t} \widetilde{u}^{+} -\nu_{1}^{+} \Delta \widetilde{u}^{+}-
\nu_{2}^{+} \nabla \operatorname{div} \widetilde{u}^{+}=\widetilde{S}_2^+, \\
\partial_{t} \widetilde{n}^{-}+u_2^-\cdot\nabla \widetilde{n}^-=\widetilde{S}_1^-, \\
\partial_{t} \widetilde{u}^{-} -\nu_{1}^{-} \Delta \widetilde{u}^{-}-
\nu_{2}^{-} \nabla \operatorname{div} \widetilde{u}^{-}=\widetilde{S}_2^-,
\end{array}\right.
\end{equation}
where
\begin{equation}
    \left\{
\begin{aligned}
    \widetilde{S}_1^{+}=&- \widetilde{u}^+\cdot\nabla n^+_1-\operatorname{div}\widetilde{u}^+-\widetilde{n}^+\operatorname{div}u_2^+-n_1^+\operatorname{div}\widetilde{u}^+,\\
    \widetilde{S}_2^+=&-\beta_1\nabla \widetilde{n}^+-\beta_2\nabla \widetilde{n}^- +S_2^+(n_2^+,n_2^-,u_2^+,u_2^-)-S_2^+(n_1^+,n_1^-,u_1^+,u_1^-),\\
    \widetilde{S}_1^-=&- \widetilde{u}^-\cdot\nabla n^-_1-\operatorname{div}\widetilde{u}^--\widetilde{n}^-\operatorname{div}u_2^--n_1^-\operatorname{div}\widetilde{u}^-,\\
    \widetilde{S}_2^-=&-\beta_3\nabla \widetilde{n}^+-\beta_4\nabla \widetilde{n}^-+S_2^+(n_2^+,n_2^-,u_2^+,u_2^-)-S_2^+(n_1^+,n_1^-,u_1^+,u_1^-).\\
\end{aligned}
    \right.
\end{equation}
In the following, we denote 
\begin{align}
    V_i^\pm(t)=\int_0^t \|u_i^\pm(\tau)\|_{\dot{B}_{p,1}^{\frac{d}{p}+1}}d\tau\ \mbox{ for } i=1,2
\end{align}
and we denote by $A_t$ a constant depending on $\|n_i^\pm\|_{\widetilde{L}_t^\infty\dot{B}_{p,1}^{\frac{d}{p}} }$ for $i=1,2.$
Due to the embedding $\mathcal{E}_{2,p}^{\frac{d}{p},+}\cup\mathcal{E}_{2,p}^{\frac{d}{p},-}\subset \mathcal{E}_{2,d}^{1,+}\cup \mathcal{E}_{2,d}^{1,-}$, it suffices to prove the uniqueness in $\mathcal{E}_{2,d}^{1,\pm}$. Therefore,
we can take $ p = d$ in the subsequent process.

Applying Proposition 3.4 gives
\begin{equation}\label{U1}
    \begin{aligned}
        \|\widetilde{n}^\pm\|_{\dot{B}_{p,\infty}^0}\lesssim e^{CV_2^\pm(t)}\int_0^t\| \widetilde{S}_1^\pm(\tau)\|_{\dot{B}_{p,\infty}^0}d\tau.
    \end{aligned}
\end{equation}
By Lemma \ref{Lem1.2}, we obtain
\begin{equation}
    \begin{aligned}
      \| \widetilde{S}_1^\pm(\tau)\|_{\dot{B}_{p,\infty}^0}\lesssim  \|u_2^\pm\|_{\dot{B}_{p,1}^2} \|\widetilde{n}^\pm\|_{\dot{B}_{p,\infty}^0}+(1+\|n_1^\pm\|_{\dot{B}_{p,1}^1})\|\widetilde{u}^\pm\|_{\dot{B}_{p,1}^1}.
    \end{aligned}
\end{equation}
Plugging it into \eqref{U1}, then applying Gr\"onwall’s inequality to get
\begin{equation}\label{U4}
    \begin{aligned}
     \|(\widetilde{n}^+,\widetilde{n}^-)\|_{\dot{B}_{p,\infty}^0}\lesssim   e^{C\max\{V_2^+(t),V_2^-(t)\}}\int_0^t(1+\|(n_1^+,n_1^-)\|_{\dot{B}_{p,1}^1})\|(\widetilde{u}^+,\widetilde{u}^-)\|_{\dot{B}_{p,1}^1}d\tau.
    \end{aligned}
\end{equation}
Applying Proposition  \ref{lame-system-estimate}   gives
\begin{equation}\label{U2}
    \begin{aligned}
\|(\widetilde{u}^+,\widetilde{u}^-)\|_{\widetilde{L}_t^1(\dot{B}_{p,\infty}^1)}  +\|(\widetilde{u}^+,\widetilde{u}^-)\|_{\widetilde{L}_t^2(\dot{B}_{p,\infty}^0)}  \lesssim \int_0^t \|( \widetilde{S}_2^+ ,\widetilde{S}_2^-) \|_{ \dot{B}_{p,\infty}^{-1}}d\tau.
    \end{aligned}
\end{equation}
From Lemma \ref{Lem1.2} and Proposition \ref{propA2}, we have 
\begin{equation}\label{U3}
    \begin{aligned}
    \|( \widetilde{S}_2^+ ,\widetilde{S}_2^-) \|_{ \dot{B}_{p,\infty}^{-1}}\lesssim& \|(u_1^+,u_1^-,u_2^+,u_2^-)\|_{ \dot{B}_{p,1}^{1}} \|(\widetilde{u}^+,\widetilde{u}^-)\|_{ \dot{B}_{p,\infty}^{0}} \\
    &+ A_t \|(n_1^+,n_1^-)\|_{ \dot{B}_{p,1}^{1}}\|(\widetilde{u}^+,\widetilde{u}^-)\|_{ \dot{B}_{p,\infty}^{1}}\\
    &+A_t (1+\|(u_2^+,u_2^-)\|_{\dot{B}_{p,1}^{2}})\|(\widetilde{n}^+,\widetilde{n}^-)\|_{\dot{B}_{p,\infty}^{0}}d\tau.
    \end{aligned}
\end{equation}
We take $T$ small enough such that 
\begin{align}
    C\|(u_1^+,u_1^-,u_2^+,u_2^-)\|_{\widetilde{L}_t^2(\dot{B}_{p,1}^{1})}\leq \frac{1}{2}.
\end{align}
Thus, plugging \eqref{U3} into \eqref{U2} and using Gr\"onwall’s inequality again, we deduce that for any  $t\in [0,T],$
\begin{equation}\label{U5}
    \begin{aligned}
      \|(\widetilde{u}^+,\widetilde{u}^-)\|_{\widetilde{L}_t^1(\dot{B}_{p,\infty}^1)}\leq A_t\int_0^t (1+\|(u_2^+,u_2^-)\|_{\dot{B}_{p,1}^{2}})\|(\widetilde{n}^+,\widetilde{n}^-)\|_{\dot{B}_{p,\infty}^{0}}d\tau. \end{aligned}
\end{equation}

\begin{lemma}\emph{\cite{Danchin2003}}\label{lem:log_inequality}
Let $s\in\mathbb{R}$. Then for any $1\le p,r\le +\infty$ and $0<\varepsilon\le 1$, we have
\[
\|f\|_{\widetilde {L}_t^{ r}\dot B^{s}_{p,1}}
\le
C\,\frac{\|f\|_{\widetilde {L}_t^{ r}\dot B^{s}_{p,\infty}}}{\varepsilon}
\log\!\left(
e+\frac{\|f\|_{\widetilde {L}_t^{ r}\dot B^{s-\varepsilon}_{p,\infty}}
+\|f\|_{\widetilde {L}_t^{ r}\dot B^{s+\varepsilon}_{p,\infty}}}
{\|f\|_{\widetilde {L}_t^{ r}\dot B^{s}_{p,\infty}}}
\right).
\]
\end{lemma}
By applying Lemma \ref{lem:log_inequality} with $s = r =\varepsilon= 1 $ and $f = \widetilde{u}^\pm$, we obtain
\begin{equation}
    \begin{aligned}
        \|\widetilde{u}^\pm\|_{\widetilde {L}_t^{1}(\dot B^{1}_{p,1})}
\le
C\, \|\widetilde{u}^\pm\|_{\widetilde {L}_t^{1}(\dot B^{1}_{p,\infty})}
\log\!\left(
e+\frac{\|\widetilde{u}^\pm\|_{\widetilde {L}_t^{ 1}(\dot B^{0}_{p,\infty})}
+\|\widetilde{u}^\pm\|_{\widetilde {L}_t^{1}(\dot B^{2}_{p,\infty})}}
{\|\widetilde{u}^\pm\|_{\widetilde {L}_t^{1}(\dot B^{1}_{p,\infty})}}
\right),
    \end{aligned}
\end{equation}
 which together with \eqref{U4} and \eqref{U5}, we conclude that
 \begin{equation}
     \begin{aligned}
      \|(\widetilde{u}^+,\widetilde{u}^-)\|_{\widetilde{L}_t^1(\dot{B}_{p,\infty}^1)}\lesssim &A_te^{C\max\{V_2^+(t),V_2^-(t)\}}\int_0^t (1+\|(u_2^+,u_2^-)\|_{\dot{B}_{p,1}^{2}})\|(\widetilde{n}^+,\widetilde{n}^-)\|_{\dot{B}_{p,\infty}^{0}} \\
      &\quad\quad\quad\quad\quad\quad\quad\quad\quad\times \log\!\left(
e+\frac{C_\tau^++C_\tau^-}
{\|(\widetilde{u}^+,\widetilde{u}^-)\|_{\widetilde {L}_t^{1}(\dot B^{1}_{p,\infty})}}
\right)d\tau,
     \end{aligned}
 \end{equation}
where $C_\tau^\pm=\|\widetilde{u}^\pm\|_{\widetilde {L}_t^{ 1}(\dot B^{0}_{p,\infty})}
+\|\widetilde{u}^\pm\|_{\widetilde {L}_t^{1}(\dot B^{2}_{p,\infty})} $. Noting  that $\|(u_2^+,u_2^-)\|_{\dot{B}_{p,1}^{2}}$
is integrable on $ [0, \infty]$ and 
\begin{align}
    \int_0^1\frac{dr}{r\log\!(r+(C_\tau^++C_\tau^-)r^{-1})}=+\infty,
\end{align}
applying the Osgood lemma leads to the conclusion that $(\widetilde{n}^+,\widetilde{n}^-,\widetilde{u}^+,\widetilde{u}^-)=0$ on $[0,T]$, and a continuity argument ensures that $(n_1^+,n_1^-,u_1^+,u_1^-)=(n_2^+,n_2^-,u_2^+,u_2^-)$ on $[0,\infty).$
\section{Time decay rate for the  non-conservative compressible two-fluid}
In this section, we aim at proving the
time-weighted energy inequality
\eqref{ineqtime} based on the global existence of solution to  the  non-conservative compressible two-fluid.
Recall that  the time-weighted energy functional is defined as
\begin{equation}\label{def-Y}
    \begin{aligned}
    \mathcal{Y}_p(t)\triangleq&\sup_{\sigma\in[\sigma_1+\varepsilon,\frac{d}{2}+1]}\|\langle \tau \rangle^{\frac{\sigma-\sigma_1}{2}}(n^+,n^-,u^+,u^-)\|_{L_t^{\infty}(\dot{B}_{2,1}^{\sigma})}^{\ell}+ \|\langle \tau \rangle^{\alpha}(n^+,n^-) \|_{\widetilde{L}_t^{\infty}(\dot{B}_{p,1}^{\frac{d}{p}})}^{h} \\
    &+\|\langle \tau\rangle^{\alpha}(u^+,u^-)\|_{\widetilde{L}_t^{\infty}(\dot{B}_{p,1}^{\frac{d}{p}-1})}^{h}+\| \tau^{\alpha}\nabla(u^+,u^-)\|_{\widetilde{L}_t^{\infty}(\dot{B}_{p,1}^{\frac{d}{p}})}^{h}
\end{aligned}
\end{equation}
with $\alpha=\frac{d}{2}+\frac{1}{2}-\sigma_1-\epsilon$ for $\epsilon>0$ sufficiently small. Subsection \ref{subsection3.1} is devoted to the low-frequency estimate. The proof is similar to \cite{CMP2019Xu}. In Subsection \ref{subsection3.2},  we develop “two effective
velocities” and obtain the upper bound for second and third terms in $\mathcal{Y}_p(t)$. Additionally, to close the
high-frequency estimates, a crucial observation enables us to
establish gain of regularity and decay altogether for the velocity, which strongly
depends on Proposition \ref{lame-system-estimate}.
\subsection{Decay estimates for the low frequencies of \texorpdfstring{$(\nabla n^+,\nabla n^-,u^+,u^-)$}{(n±, u±)}}\label{subsection3.1}
\begin{proposition} If $p$ satisfies \eqref{range-p}, then it holds that  for all $t\geq 0,$ 
    \begin{equation}\label{eq:low-frequency}
        \begin{aligned}
            \|(n^+,n^-,u^+,u^-)\|_{L_t^{\infty}(\dot{B}_{2,1}^{\sigma})}^{\ell}\lesssim \langle \tau \rangle^{-\frac{\sigma-\sigma_1}{2}}( \| (n_0^+,n_0^-,u_0^+,u_0^-) \|^{\ell}_{ \dot{B}_{2,\infty}^{\sigma_1}}+ \mathcal{X}_p^2(t)+\mathcal{Y}_p^2(t) ),
        \end{aligned}
    \end{equation}
    provided that $\sigma_1<\sigma\leq \frac{d}{2}+1$, where $\mathcal{X}_p(t)$ and $\mathcal{Y}_p(t)$ have been defined by \eqref{def-E} and \eqref{def-Y}, respectively.
\end{proposition}
\begin{proof}
By \eqref{eq:3}, we infer that
    \begin{equation} \label{eq:dyadic-energy}
    \begin{aligned}
         \frac{d}{dt}\| (n_j^+,n_j^-,u_j^+,u_j^-) \|_{L^2}+c_02^{2j}\|(n_j^+,n_j^-,u_j^+,u_j^-)\|_{L^2}
    \lesssim\|\dot{\Delta}_j(\mathcal{F}_1^+,\mathcal{F}_2^+,\mathcal{F}_1^-,\mathcal{F}_2^-)\|_{L^2},
    \end{aligned}
\end{equation}
for some constant $c_0>0$ independent of $j$ and $t$.
Integrating \eqref{eq:dyadic-energy} in time and applying the Duhamel principle yields
    \begin{equation} \label{eq:duhamel}
    \begin{aligned}
          \| (n_j^+,n_j^-,u_j^+,u_j^-) \|_{L^2} 
        \lesssim& e^{-c_02^{2j}t}\| \dot{\Delta}_j(n_0^+,n_0^-,u_0^+,u_0^-) \|_{ L^2}\\
        &+\int_0^te^{-c_02^{2j}(t-\tau)}\|\dot{\Delta}_j(\mathcal{F}_1^+,\mathcal{F}_2^+,\mathcal{F}_1^-,\mathcal{F}_2^-)\|_{L^2}d\tau.
    \end{aligned}
\end{equation}
We now estimate the two terms on the right-hand side of \eqref{eq:duhamel} at low frequencies, i.e. for $j\le j_0$.
For the contribution of the initial data, we obtain
\begin{equation}
    \begin{aligned}
        &\langle t \rangle^{\frac{\sigma-\sigma_1}{2}}\sum_{j\leq j_0}e^{-c_02^{2j}t}2^{j\sigma}\| \dot{\Delta}_j(n_0^+,n_0^-,u_0^+,u_0^-) \|_{ L^2}\\
        &\lesssim\langle t \rangle^{\frac{\sigma-\sigma_1}{2}}\sum_{j\leq j_0}e^{-c_02^{2j}t}2^{j(\sigma-\sigma_1)}\| (n_0^+,n_0^-,u_0^+,u_0^-) \|_{ \dot{B}_{2,\infty}^{\sigma_1}}\\
        &\lesssim \| (n_0^+,n_0^-,u_0^+,u_0^-) \|^{\ell}_{ \dot{B}_{2,\infty}^{\sigma_1}}.
    \end{aligned}
\end{equation}
Next, for the second term in \eqref{eq:duhamel}, summing over $j\le j_0$ gives
\begin{equation}\label{eq:low-duhamel}
    \begin{aligned}
        &\sum_{j\leq j_0}2^{j\sigma}\int_0^te^{-c_02^{2j}(t-\tau)}\|\dot{\Delta}_j(\mathcal{F}_1^+,\mathcal{F}_2^+,\mathcal{F}_1^-,\mathcal{F}_2^-)\|_{L^2}d\tau\\
        &\lesssim\int_0^t\langle t-\tau \rangle^{-\frac{\sigma-\sigma_1}{2}}\|(\mathcal{F}_1^+,\mathcal{F}_2^+,\mathcal{F}_1^-,\mathcal{F}_2^-)\|^{\ell}_{ \dot{B}_{2,\infty}^{\sigma_1}}d\tau.
    \end{aligned}
\end{equation}
To proceed, we recall the schematic structure of the nonlinear terms. More precisely,
\begin{equation}\label{eq:F1F2-structure}
\left\{
\begin{aligned}
\mathcal{F}_1^+ &\sim n^+\,\mathrm{div}\,u^+ + u^+\cdot\nabla n^+,\\
\mathcal{F}_2^+ &\sim g^+(n^+,n^-)\,\nabla n^+ + (u^+\cdot\nabla)u^+
      + h^+(n^+,n^-)\,\nabla n^+\nabla u^+
      + l^+(n^+,n^-)\,\Delta u^+ .
\end{aligned}
\right.
\end{equation}
Accordingly, it is convenient to decompose $\mathcal{F}_1$ and $\mathcal{F}_2$ into their low- and high-frequency parts:
\begin{align*}
(\mathcal{F}_1^+)^{\ell} &\sim n^+\,\mathrm{div}(u^+)^{\ell} + u^+\cdot\nabla (n^+)^{\ell},
\qquad
(\mathcal{F}_1^+)^{h} \sim n^+\,\mathrm{div}(u^+)^{h} + u^+\cdot\nabla (n^+)^{h},
\end{align*}
and
\begin{align*}
(\mathcal{F}_2^+)^{\ell} &\sim g^+(n^+,n^-)\,\nabla (n^+)^{\ell} + (u^+\cdot\nabla)(u^+)^{\ell}
+ h^+(n^+,n^-)\,\nabla n^+\nabla (u^+)^{\ell} + l^+(n^+,n^-)\,\Delta (u^+)^{\ell},\\
(\mathcal{F}_2^+)^{h} &\sim g^+(n^+,n^-)\,\nabla (n^+)^{h} + (u^+\cdot\nabla)(u^+)^{h}
+ h^+(n^+,n^-)\,\nabla n^+\nabla (u^+)^{h} + l^+(n^+,n^-)\,\Delta (u^+)^{h},
\end{align*}
where, for $z=n^{\pm},u^{\pm}$,
\[
z^{\ell}\triangleq\sum_{j\le j_0}\dot{\Delta}_j z,\qquad
z^{h}\triangleq z-z^{\ell}.
\]
Then \eqref{eq:low-frequency} holds if the following two lemmas are true.
\end{proof}

\begin{lemma} If $p$ satisfy \eqref{range-p}  and $\sigma_0\leq\sigma_1<\frac{d}{2}-1 $, then it holds that for all $t\geq 0,$
    \begin{equation}
        \begin{aligned}
             \int_0^t\langle t-\tau \rangle^{-\frac{\sigma-\sigma_1}{2}}\|(\mathcal{F}_1^+,\mathcal{F}_2^+,\mathcal{F}_1^-,\mathcal{F}_2^-)^{\ell}\|^{\ell}_{ \dot{B}_{2,\infty}^{\sigma_1}}d\tau\lesssim \langle \tau \rangle^{-\frac{\sigma-\sigma_1}{2}}(\mathcal{X}_p^2(t)+\mathcal{Y}_p^2(t)  ),
        \end{aligned}
    \end{equation}
     where $\mathcal{X}_p(t)$ and $\mathcal{Y}_p(t)$ have been defined by \eqref{def-E} and \eqref{def-Y}, respectively.
   \end{lemma}
\begin{proof}
Following the approach of  Lemma 4.2 in \cite{CMP2019Xu}, we can estimate $(\mathcal{F}_1^+)^{\ell},(\mathcal{F}_2^+)^{\ell},(\mathcal{F}_1^-)^{\ell},$ and $(\mathcal{F}_2^-)^{\ell}$. Hence, the reader may apply the same procedure to each nonlinear term to derive the desired bounds.
\end{proof}

\begin{lemma} If $p$ satisfy \eqref{range-p}  and $\sigma_0\leq\sigma_1<\frac{d}{2}-1 $, then it holds that for all $t\geq 0,$
    \begin{equation}
        \begin{aligned}
             \int_0^t\langle t-\tau \rangle^{-\frac{\sigma-\sigma_1}{2}}\|(\mathcal{F}_1^+,\mathcal{F}_2^+,\mathcal{F}_1^-,\mathcal{F}_2^-)^{h}\|^{\ell}_{ \dot{B}_{2,\infty}^{\sigma_1}}d\tau\lesssim \langle \tau \rangle^{-\frac{\sigma-\sigma_1}{2}}(\mathcal{X}_p^2(t)+\mathcal{Y}_p^2(t) ),
        \end{aligned}
    \end{equation}
     where $\mathcal{X}_p(t)$ and $\mathcal{Y}_p(t)$ have been defined by \eqref{def-E} and \eqref{def-Y}, respectively.
   \end{lemma}
\begin{proof}
We first distinguish between the non-oscillation and oscillation cases. In the non-oscillation case, we further split according to the range of $\sigma_1$, namely  $\sigma_0\le \sigma_1<\frac d2-\frac dp$  and $\frac d2-\frac dp\le \sigma_1<\frac d2-1$. Following the strategy of \cite{CMP2019Xu}, we then derive the nonlinear estimates for  $(\mathcal{F}_1^+)^{h},(\mathcal{F}_2^+)^{h},(\mathcal{F}_1^-)^{h}$, and $ (\mathcal{F}_2^-)^{h} $ in these three situations by applying Lemmas 4.3, 4.4, and 4.5 in \cite{CMP2019Xu}, respectively. The same procedure can be applied to each nonlinear terms, yielding the desired bounds.
   
\end{proof}
   
\subsection{Decay Estimates for the high frequencies of \texorpdfstring{$(\nabla n^+,\nabla n^-,u^+,u^-)$}{(n±, u±)}}\label{subsection3.2}

\begin{lemma}If $p$ satisfy \eqref{range-p}  and $\sigma_0\leq\sigma_1<\frac{d}{2}-1 $, then it holds that for all $T\geq 0,$
 \begin{equation}\label{eq:high-frequency}
    \begin{aligned}
        &\|\langle \tau\rangle^{\alpha}(\nabla n^+,\nabla n^-,u^+,u^-)\|^h_{\widetilde{L}^{\infty}_t(\dot{B}_{p,1}^{\frac{d}{p}-1})}+ \| \tau^{\alpha}\nabla(u^+,u^-)\|_{\widetilde{L}_t^{\infty}(\dot{B}_{p,1}^{\frac{d}{p}})}^{h}\\
        \lesssim& \|(\nabla n_0^+,\nabla n_0^-,u_0^+,u_0^-)\|^h_{\dot{B}_{p,1}^{\frac{d}{p}-1}}+\mathcal{X}_p(t)\mathcal{Y}(t)+ \mathcal{Y}_p^2(t).
\end{aligned}
\end{equation}
with $\alpha=\frac{d}{2}+\frac{1}{2}-\sigma_1-\epsilon$ for $\epsilon>0$ sufficiently small,  where $\mathcal{X}_p(T)$ and $\mathcal{Y}_p(T)$ have been defined by \eqref{def-E} and \eqref{def-Y}, respectively.
\end{lemma}

 \begin{proof}
     
Acting $\mathcal{P} $  onto $(\ref{equ:fornu-pm})_2$ and $(\ref{equ:fornu-pm})_4$, and then applying $\dot{\Delta}_j$ to the resulting systems, we obtain 
\begin{equation} 
 \left\{
   \begin{array}{ll}
     \partial_t \mathcal{P}u_j^+-\nu_1^+ \Delta \mathcal{P}u_j^+  =\mathcal{P}\dot{\Delta}_j\mathcal{F}_2^+,\\
     \partial_t\mathcal{P}u_j^- -\nu_1^-\Delta \mathcal{P}u_j^-  =\mathcal{P}\dot{\Delta}_j\mathcal{F}_2^-.
   \end{array}
 \right.
\end{equation}
Then multiplying each component of the above equation by $\mathcal{P}u_j^+|\mathcal{P}u_j^+|^{p-2}$ and $\mathcal{P}u_j^-|\mathcal{P}u_j^-|^{p-2}$ to get 
\begin{equation}\label{eq-energyforPu}
    \begin{aligned}
         &\frac{1}{p}\frac{d}{dt}\|(\mathcal{P}u_j^+,\mathcal{P}u_j^-) \|_{L^p}^p+c_p\nu_12^{2j}\|(\mathcal{P}u_j^+,\mathcal{P}u_j^-) \|_{L^p}^p\leq \|\mathcal{P}\dot{\Delta}_j(\mathcal{F}_2^+,\mathcal{F}_2^-)\|_{L^p}\|(\mathcal{P}u_j^+,\mathcal{P}u_j^-) \|_{L^p}^{p-1},
    \end{aligned}
\end{equation}
where $\nu_1=\min\{\nu_1^+,\nu_1^-\}.$ Let $\|f\|_{\epsilon,L^p}=(\|f\|^p_{L^p}+\epsilon^p)^{1/p}$. Then by Lemma \ref{Lem1.2}, we get
\begin{equation}
    \begin{aligned}
         &\frac{d}{dt}\|(\mathcal{P}u_j^+,\mathcal{P}u_j^-) \|_{\epsilon,L^p}+c_p\nu_12^{2j}\|(\mathcal{P}u_j^+,\mathcal{P}u_j^-) \|_{\epsilon,L^p}\leq \|\mathcal{P}\dot{\Delta}_j(\mathcal{F}_2^+,\mathcal{F}_2^-)\|_{L^p}+c_p\nu_12^{2j}\epsilon.
    \end{aligned}
\end{equation}
At this point, following Haspot’s method in \cite{MR2823668,MR2847531}, we introduce the effective
velocity $w^+$ and $w^-$ defined in \eqref{effective-velo-w}. It is clear that $w^+$ and $w^-$  fulfill
\begin{equation} 
 \left\{
   \begin{array}{ll}
     \partial_t w^+-(\nu_1^++\nu_2^+)\Delta w^+  =g_1,\\
     \partial_t w^- -(\nu_1^-+\nu_2^-)\Delta w^-  =g_2,
   \end{array}
 \right.
\end{equation}
Hence, arguing exactly as for proving \eqref{eq-energyforPu},
\begin{equation}
    \begin{aligned}
        &\frac{1}{p}\frac{d}{dt}\|(w^+_j,w^-_j)\|_{L^p}^p+c_p2^{2j}\|(w^+_j,w^-_j)\|_{L^p}^p\\
        &\leq \left(\|\dot{\Delta}_j(\mathcal{F}_2^+,\mathcal{F}_2^-)\|_{L^p}+2^{-j}\|\dot{\Delta}_j(\mathcal{F}_1^+,\mathcal{F}_1^-)\|_{L^p}+\|(w^+_j,w^-_j)\|_{L^p}+2^{-j}\|(n^+_j,n^-_j)\|_{L^p}\right)\|(w^+_j,w^-_j)\|_{L^p}^{p-1}   .
        \end{aligned}
\end{equation}
By using Lemma \ref{Lem1.2}, we get
\begin{equation}
    \begin{aligned}
        & \frac{d}{dt}\|(w^+_j,w^-_j)\|_{\epsilon,L^p} +c_p2^{2j}\|(w^+_j,w^-_j)\|_{\epsilon,L^p} \\
        &\leq  \|\dot{\Delta}_j(\mathcal{F}_2^+,\mathcal{F}_2^-)\|_{L^p}+2^{-j}\|\dot{\Delta}_j(\mathcal{F}_1^+,\mathcal{F}_1^-)\|_{L^p}+\|(w^+_j,w^-_j)\|_{L^p} +2^{-j}\|(n^+_j,n^-_j)\|_{L^p}+c_p2^{2j}\epsilon.
        \end{aligned}
\end{equation}
In terms of $w^\pm$, the function $n^\pm$ satisfy the following two damped transport equations
\begin{equation}\label{damping-for-npm}
\left\{
    \begin{aligned}
      \partial_t  n^++u^+\cdot\nabla n^+ +   \frac{\gamma_1}{\nu_1^++\nu_2^+} n^+ +\frac{\gamma_2}{\nu_1^++\nu_2^+} n^-=- n^+\text{div}u^+- \text{div}w^+,\\
        \partial_t  n^-+u^-\cdot\nabla n^-+  \frac{\gamma_3}{\nu_1^-+\nu_2^-} n^+ +\frac{\gamma_4}{\nu_1^-+\nu_2^-} n^-=- n^-\text{div}u^--\text{div}w^-.
    \end{aligned}
    \right.
\end{equation}
Applying $\partial_k\dot{\Delta}_j$ to \eqref{damping-for-npm}, we obtain
\begin{equation}\label{damping-for-npm-local-k}
\left\{
    \begin{aligned}
      \partial_t\partial_k n_j^++u^+\cdot\nabla n_j^+ +   \frac{\gamma_1}{\nu_1^++\nu_2^+}\partial_kn_j^+ +\frac{\gamma_2}{\nu_1^++\nu_2^+}\partial_kn_j^-=-\partial_k\dot{\Delta}_j(n^+\text{div}u^+)-\partial_k\text{div}w_j^++R_j^{k,+},\\
        \partial_t \partial_kn_j^-+u^-\cdot\nabla n_j^-+  \frac{\gamma_3}{\nu_1^-+\nu_2^-}\partial_kn_j^+ +\frac{\gamma_4}{\nu_1^-+\nu_2^-}\partial_kn_j^-=-\partial_k\dot{\Delta}_j(n^-\text{div}u^-)-\partial_k\text{div}w_j^-+R_j^{k,-},
    \end{aligned}
    \right.
\end{equation}
where $R_j^{k,\pm}=[u^\pm\cdot\nabla,\partial_k\dot{\Delta}_j]n^\pm$. Multiplying $(\ref{damping-for-npm-local-k})_1$  by $\frac{1}{\gamma_2}\partial_k n_j^+|\partial_kn_j^+|^{p-2}$ and $(\ref{damping-for-npm-local-k})_2$  by $\frac{1}{\gamma_3}\partial_k n_j^-|\partial_kn_j^-|^{p-2}$, and summing up on $k=1,\cdots,d$, and   applying H\"older and Bernstein inequalities
leads to
\begin{equation}\label{eq:n+1}
    \begin{aligned}
        &\frac{1}{p\gamma_2}\frac{d}{dt}\|\nabla n_j^+\|^p_{L^p}+\frac{\gamma_+^2}{\nu_1^++\nu_2^+}\|\nabla 
 n_j^+\|^p_{L^p}\leq \frac{1}{\gamma_2}\left( \frac{\gamma_2}{\nu_1^++\nu_2^+} \|\nabla  n_j^-\|_{L^p} \right.\\
       +&\left.\frac{1}{p}\|\text{div}u^+\|_{L^{\infty}}\|\nabla n_j^+\|_{L^p}+\|\nabla \dot{\Delta}_j(n^+\text{div}u^+)\|_{L^p}+\|\nabla \text{div}w_j^+\|_{L^p}+\|R_j^{k,+}\|_{L^p}\right)\|\nabla n_j^+\|_{L^p}^{p-1}
    \end{aligned}
\end{equation}
and 
\begin{equation}\label{eq:n+2}
    \begin{aligned}
        &\frac{1}{p\gamma_3}\frac{d}{dt}\|\nabla n_j^-\|^p_{L^p}+\frac{\gamma_-^2}{\nu_1^-+\nu_2^-}\|\nabla n_j^-\|^p_{L^p}\leq \frac{1}{\gamma_3}\left( \frac{\gamma_3}{\nu_1^-+\nu_2^-} \|\nabla n_j^+\|_{L^p}  \right.\\
       +&\left.\frac{1}{p}\|\text{div} u^-\|_{L^{\infty}}\|\nabla n_j^-\|_{L^p}+\|\nabla \dot{\Delta}_j(n^-\text{div}u^-)\|_{L^p}+\|\nabla \text{div}w_j^-\|_{L^p}+\|R_j^{k,-}\|_{L^p}\right)\|\nabla n_j^-\|_{L^p}^{p-1}.
    \end{aligned}
\end{equation}
Once again using Lemma \ref{Lem1.2} gives
\begin{equation}\label{eq:n+11}
    \begin{aligned}
        &\frac{1}{\gamma_2}\frac{d}{dt}\|n_j^+\|_{\epsilon,L^p}+\frac{\gamma_+^2}{\nu_1^++\nu_2^+}\|n_j^+\|_{\epsilon,L^p}\leq \frac{1}{\gamma_2}\left( \frac{\gamma_2}{\nu_1^++\nu_2^+} \|n_j^-\|_{ L^p} \right.\\
       +&\left.\frac{1}{p}\|\text{div}u^+\|_{L^{\infty}}\|n_j^+\|_{L^p}+\|\dot{\Delta}_j(n^+\text{div}u^+)\|_{L^p}+\|\text{div}w_j^+\|_{L^p}+\|R_j^{k,+}\|_{L^p}\right)+\frac{\gamma_+^2}{\nu_1^++\nu_2^+}\epsilon
    \end{aligned}
\end{equation}
and 
\begin{equation}\label{eq:n+22}
    \begin{aligned}
        &\frac{1}{\gamma_3}\frac{d}{dt}\|n_j^-\|_{\epsilon,L^p}+\frac{\gamma_-^2}{\nu_1^-+\nu_2^-}\|n_j^-\|_{\epsilon,L^p}\leq \frac{1}{\gamma_3}\left( \frac{\gamma_3}{\nu_1^-+\nu_2^-} \|n_j^+\|_{L^p}  \right.\\
       +&\left.\frac{1}{p}\|\text{div} u^-\|_{L^{\infty}}\|n_j^-\|_{L^p}+\|\dot{\Delta}_j(n^-\text{div}u^-)\|_{L^p}+\|\text{div}w_j^-\|_{L^p}+\|R_j^{k,-}\|_{L^p}\right) +\frac{\gamma_-^2}{\nu_1^-+\nu_2^-}\epsilon.
    \end{aligned}
\end{equation}
From \eqref{condition:gamma}, we notice that $\gamma_1\gamma_4>\gamma_2\gamma_3$. And this implies that $\gamma_+^2\gamma_-^2>1$, then there exists a constant $\tilde{\eta}\in\left(\frac{\nu_1^-+\nu_2^-}{(\nu_1^++\nu_2^+)\gamma_-^2},\frac{\nu_1^-+\nu_2^-}{\nu_1^++\nu_2^+}\gamma_+^2\right)$, such that 
\begin{equation}\label{eq:n+-1}
    \begin{aligned}
        &\frac{d}{dt}\|\nabla(n_j^+,n_j^-)\|_{\epsilon,L^p}+\kappa\|\nabla(n_j^+,n_j^-)\|_{\epsilon,L^p}\\
        \leq&\,   
       \frac{1}{p}\|\text{div} u^-\|_{L^{\infty}}\|\nabla n_j^-\|_{L^p}+\|\nabla\dot{\Delta}_j(n^-\text{div}u^-)\|_{L^p}+\|\nabla\text{div}w_j^-\|_{L^p}+\|R_j^{k,-}\|_{L^p}\\
       &\, +\frac{1}{p}\|\text{div}u^+\|_{L^{\infty}}\|\nabla n_j^+\|_{L^p}+\|\nabla\dot{\Delta}_j(n^+\text{div}u^+)\|_{L^p}+\|\nabla\text{div}w_j^+\|_{L^p}+\|R_j^{k,+}\|_{L^p}\\
       &\, + \frac{\gamma_+^2}{\nu_1^++\nu_2^+}\epsilon+\frac{\gamma_-^2}{\nu_1^-+\nu_2^-}\epsilon.
    \end{aligned}
\end{equation}
Adding up that inequality 
\begin{equation}
    \begin{aligned}
        &\frac{d}{dt}( \|(\mathcal{P}u_j^+,\mathcal{P}u_j^-) \|_{\epsilon,L^p}+ \|(w^+_j,w^-_j)\|_{\epsilon,L^p}+ \gamma c_p\|\nabla(n_j^+,n_j^-)\|_{\epsilon,L^p})\\
        &\,+c_p\nu_12^{2j}\|(\mathcal{P}u_j^+,\mathcal{P}u_j^-) \|_{\epsilon,L^p}+c_p2^{2j}\|(w^+_j,w^-_j)\|_{\epsilon,L^p}+\gamma c_p\kappa\|\nabla(n_j^+,n_j^-)\|_{\epsilon,L^p}\\
        \leq&\,\|(w^+_j,w^-_j)\|_{L^p}+2^{-2j_0}\|\nabla(n^+_j,n^-_j)\|_{L^p}\\
        &\,+\|\mathcal{P}\dot{\Delta}_j(\mathcal{F}_2,F_4)\|_{L^p}+ \|\dot{\Delta}_j(\mathcal{F}_2,F_4)\|_{L^p}+\|\dot{\Delta}_j(n^+u^+,n^-u^-)\|_{L^p}\\
        &\,+\gamma c_p\left(\frac{1}{p}\|\text{div} u^-\|_{L^{\infty}}\|\nabla n_j^-\|_{L^p}+\|\nabla\dot{\Delta}_j(n^-\text{div}u^-)\|_{L^p}+\|R_j^{k,-}\|_{L^p}\right)\\
       &\,+\gamma c_p\left(\frac{1}{p}\|\text{div}u^+\|_{L^{\infty}}\|\nabla n_j^+\|_{L^p}+\|\nabla\dot{\Delta}_j(n^+\text{div}u^+)\|_{L^p}+\|R_j^{k,+}\|_{L^p}\right)\\
       &\,+C\gamma c_p2^{2j}\|(w_j^+,w_j^-)\|_{L^p}+M_{\epsilon},
    \end{aligned}
\end{equation}
where $M_{\epsilon}= \left(c_p\nu_12^{2j}+ c_p 2^{2j}+  \gamma c_p\frac{\gamma_+^2}{\nu_1^++\nu_2^+} + \gamma c_p\frac{\gamma_-^2}{\nu_1^-+\nu_2^-}\right)\epsilon.$ Choosing $j_0$ suitably large and $\gamma$ sufficiently small, we thus discover that
\begin{equation}
    \begin{aligned}
        &\frac{d}{dt}( \|(\mathcal{P}u_j^+,\mathcal{P}u_j^-) \|_{\epsilon,L^p}+ \|(w^+_j,w^-_j)\|_{\epsilon,L^p}+ \gamma c_p\|\nabla(n_j^+,n_j^-)\|_{\epsilon,L^p})\\
        &\,+c_0\left(\|(\mathcal{P}u_j^+,\mathcal{P}u_j^-) \|_{\epsilon,L^p}+\|(w^+_j,w^-_j)\|_{\epsilon,L^p}+\gamma c_p\|\nabla(n_j^+,n_j^-)\|_{\epsilon,L^p}\right)\\
        \leq&\, \|\dot{\Delta}_j(\mathcal{F}_2^+,\mathcal{F}_2^-)\|_{L^p}+\|\dot{\Delta}_j(n^+u^+,n^-u^-)\|_{L^p}\\
        &\,+\gamma c_p\left(\frac{1}{p}\|\text{div} u^-\|_{L^{\infty}}\|\nabla n_j^-\|_{L^p}+\|\nabla\dot{\Delta}_j(n^-\text{div}u^-)\|_{L^p}+\|R_j^{k,-}\|_{L^p}\right)\\
       &\,+\gamma c_p\left(\frac{1}{p}\|\text{div}u^+\|_{L^{\infty}}\|\nabla n_j^+\|_{L^p}+\|\nabla\dot{\Delta}_j(n^+\text{div}u^+)\|_{L^p}+\|R_j^{k,+}\|_{L^p}\right)\\
       &\,+M_{\epsilon}.
    \end{aligned}
\end{equation}
Integrating in time and having $\epsilon$ tend to $0$, we eventually get 
\begin{equation}
    \begin{aligned}
        &e^{c_0t} (\|(\mathcal{P}u_j^+,\mathcal{P}u_j^-,w^+_j,w^-_j,\nabla n_j^+,\nabla n_j^-) \|_{ L^p} \\
        \lesssim&\|(\mathcal{P}u_j^+,\mathcal{P}u_j^-,w^+_j,w^-_j,\nabla n_j^+,\nabla n_j^-(0)\|_{ L^p} +\int_0^te^{c_0\tau}G^{\pm}_j(\tau)d\tau
    \end{aligned}
\end{equation}
with $G_j^{\pm}=G_{1j}^{\pm}+G_{2j}^{\pm}+G_{3j}^{\pm}+G_{4j}^{\pm}+G_{5j}^{\pm
}$ and 
\begin{align}
    &G_{1j}^{\pm}=\|\dot{\Delta}_j(n^\pm u^\pm)\|_{L^p}, \ G_{2j}^{\pm}=\|\nabla\dot{\Delta}_j(n^\pm\text{div}u^\pm)\|_{L^p}\nonumber\\
    &G_{3j}^{\pm}=\|\text{div}u^\pm\|_{L^{\infty}}\|\nabla n_j^\pm\|_{L^p}, \ G_{3j}^{\pm}=\|R_j^{k,\pm}\|_{L^p}, \ G_{5j}^{\pm}=\|\dot{\Delta}_j(\mathcal{F}_2^+,\mathcal{F}_2^-)\|_{L^p}.
\end{align}

It is clear that
\begin{equation}
    \left\{
\begin{aligned}
    u^+=w^+-\nabla (-\Delta)^{-1}\left( \frac{\gamma_1}{\nu_1^++\nu_2^+}n^+ +\frac{\gamma_2}{\nu_1^++\nu_2^+}n^-\right)+\mathcal{P}u^+,\\
    u^-=w^--\nabla (-\Delta)^{-1}\left( \frac{\gamma_3}{\nu_1^-+\nu_2^-}n^+ +\frac{\gamma_4}{\nu_1^-+\nu_2^-}n^-\right)+\mathcal{P}u^-.
\end{aligned}
\right.
\end{equation}
This leads for $j\geq j_0-1$ to 
\begin{align}
    \|u_j^\pm-(w^\pm_j+\mathcal{P}u_j^\pm)\|_{L^p}\lesssim2^{-2j_0}\|\nabla(n_j^+,n_j^-)\|_{L^p}.
\end{align}
Therefore, there exists a constant $c_0>0$ such that for all $j\geq j_0-1$, we have
\begin{equation}
    \begin{aligned}
        \|(\nabla n_j^+,\nabla n_j^-,u_j^+,u_j^-) \|_{L^2}\lesssim e^{-c_0t}\|\dot{\Delta}_j(\nabla n_0^+,\nabla n_0^-,u_0^+,u_0^-) \|_{L^2}+\int_0^t e^{-c_0(t-\tau)}G^{\pm}_j(\tau)d\tau.
    \end{aligned}
\end{equation}
Now, multiplying both sides by $\langle t\rangle^{\alpha}2^{j(\frac{d}{p}-1)}$, taking the supremum on $t\geq0$, and
summing up over $j\geq j_0-1,$ we have
\begin{equation}
    \begin{aligned}
        &\|\langle \tau\rangle^{\alpha}(\nabla n^+,\nabla n^-,u^+,u^-)\|^h_{\widetilde{L}^{\infty}_t(\dot{B}_{p,1}^{\frac{d}{p}-1})}\lesssim \|(\nabla n_0^+,\nabla n_0^-,u_0^+,u_0^-)\|^h_{\dot{B}_{p,1}^{\frac{d}{p}-1}}\\
        &+\sum_{j\geq j_0-1}\sup_{t\geq0}\left(\langle t\rangle^{\alpha}\int_0^te^{-c_0(t-\tau)}2^{j(\frac{d}{p}-1)}G^{\pm}_j(\tau)d\tau\right) .
\end{aligned}
\end{equation}

To control the nonlinear terms of, we first consider the case  $t\leq 2,$ 
\begin{equation}
    \begin{aligned}
        \sum_{j\geq j_0-1}\sup_{0\leq t\leq2}\left(\langle t\rangle^{\alpha}\int_0^te^{-c_0(t-\tau)}2^{j(\frac{d}{p}-1)}G^{\pm}_j(\tau)d\tau\right)\lesssim\int_0^2\sum_{j\geq j_0-1}2^{j(\frac{d}{p}-1)}G^{\pm}_j(\tau)d\tau.
    \end{aligned}
\end{equation}
Furthermore, it follows from that
\begin{equation}
    \begin{aligned}
        \int_0^2\sum_{j\geq j_0-1}2^{j(\frac{d}{p}-1)}G^{\pm}_j(\tau)d\tau\lesssim \int_0^2\|n^\pm u^\pm\|^h_{\dot{B}_{p,1}^{\frac{d}{p}-1}}+\|\nabla u^\pm\|_{\dot{B}_{p,1}^{\frac{d}{p} }}\|n^\pm\|_{\dot{B}_{p,1}^{\frac{d}{p}}}+\|(\mathcal{F}_2^+,\mathcal{F}_2^-)\|^h_{\dot{B}_{p,1}^{\frac{d}{p}-1}}d\tau.
    \end{aligned}
\end{equation}
Thanks to Lemma \ref{Lemma1.2}, we have 
\begin{equation*}
    \begin{aligned}
       &\|n^\pm u^\pm\|^h_{L^1_t(\dot{B}_{p,1}^{\frac{d}{p}-1})}\lesssim\|n^\pm u^\pm\|^h_{L^1_t(\dot{B}_{p,1}^{\frac{d}{p}})}\lesssim \|n^\pm\|_{\widetilde{L}^2_t(\dot{B}_{p,1}^{\frac{d}{p}})}\|u^\pm\|_{\widetilde{L}^2_t(\dot{B}_{p,1}^{\frac{d}{p}})} 
    \end{aligned}
\end{equation*}
and 
\begin{equation}
    \begin{aligned}
         &\|(\mathcal{F}_2^+,\mathcal{F}_2^-)\|^h_{L^1_t(\dot{B}_{p,1}^{\frac{d}{p}-1})}\\
         \lesssim &\|(u^+,u^-)\|_{\widetilde{L}^\infty_t(\dot{B}_{p,1}^{\frac{d}{p}-1})}\|\nabla (u^+,u^-))\|_{L^1_t(\dot{B}_{p,1}^{\frac{d}{p}})}+\|(n^+,n^-)\|_{\widetilde{L}^\infty_t(\dot{B}_{p,1}^{\frac{d}{p}})}\|(u^+,u^-)\|_{L^1_t(\dot{B}_{p,1}^{\frac{d}{p}+1})}\\ +&\|(n^+,n^-)\|_{\widetilde{L}^2_t(\dot{B}_{p,1}^{\frac{d}{p}})}\|\nabla (n^+,n^-)\|_{\widetilde{L}^2_t(\dot{B}_{p,1}^{\frac{d}{p}-1})}+\|(n^+,n^-)\|^2_{\widetilde{L}^\infty_t(\dot{B}_{p,1}^{\frac{d}{p} })}\|\nabla (u^+,u^-)\|_{L^1_t(\dot{B}_{p,1}^{\frac{d}{p}})}.
    \end{aligned}
\end{equation}
Furthermore, we can conclude that
\begin{equation}
    \begin{aligned}
         \sum_{j\geq j_0-1}\sup_{0\leq t\leq 2}\left(\langle t\rangle^{\alpha}\int_0^te^{-c_0(t-\tau)}2^{j(\frac{d}{p}-1)}G^{\pm}_j(\tau)d\tau\right)\lesssim\mathcal{X}_p^{2}(2).
    \end{aligned}
\end{equation}

Secondly, let us handle the supremum for $  t\geq 2$. To do this,
one can split the integral on $[0, t]$ into integrals $[0, 1]$ and $[1, t]$. The $[0, 1]$ part of the integral may be bounded by $\mathcal{X}^{2}(1)$. Indeed, owing to $e^{-c_0(t-\tau)}\leq e^{-\frac{c_0t}{2}}$ for $ t\geq 2$ and $0\leq \tau \leq1$, we obtain
\begin{equation}
    \begin{aligned}
         &\sum_{j\geq j_0-1}\sup_{t\geq 2}\left(\langle t\rangle^{\alpha}\int_0^1e^{-c_0(t-\tau)}2^{j(\frac{d}{p}-1)}G^{\pm}_j(\tau)d\tau\right)\\
        \lesssim &\sum_{j\geq j_0-1}\sup_{t\geq 2}\left(\langle t\rangle^{\alpha} e^{-\frac{c_0t}{2}}\int_0^1 2^{j(\frac{d}{p}-1)}G^{\pm}_j(\tau)d\tau\right)\\
        \lesssim &\sum_{j\geq j_0-1}\int_0^1 2^{j(\frac{d}{p}-1)}G^{\pm}_j(\tau)d\tau
        \lesssim\mathcal{X}_p^{2}(1)+\mathcal{X}_p^{3}(1).
    \end{aligned}
\end{equation}
In order to handle the $[1, t]$ part of the integral for $t\geq 2$, we start from
\begin{equation}
    \begin{aligned}
        \sum_{j\geq j_0-1}\sup_{t\geq 2}\left(\langle t\rangle^{\alpha}\int_1^te^{-c_0(t-\tau)}2^{j(\frac{d}{p}-1)}G^{\pm}_j(\tau)d\tau\right)\lesssim\sum_{j\geq j_0-1}2^{j(\frac{d}{p}-1)}\sup_{1\leq \tau\leq t}\tau^{\alpha} G^{\pm}_j(\tau).
    \end{aligned}
\end{equation}
For the term $n^+u^+$, we decompose it as follows
$$n^+u^+=n^+(u^+)^h+(n^+)^{h}(u^+)^{\ell}+(n^+)^{\ell}(u^+)^{\ell},$$
\begin{equation*}
    \begin{aligned}
       \|\tau^{\alpha}n^+(u^+)^h \|^h_{\widetilde{L}^\infty_t(\dot{B}_{p,1}^{\frac{d}{p}-1})} \lesssim \|n^+\|_{\widetilde{L}^\infty_t(\dot{B}_{p,1}^{\frac{d}{p}})}\|\tau^{\alpha}(u^+)^h\|_{\widetilde{L}^\infty_t(\dot{B}_{p,1}^{\frac{d}{p}-1})}\lesssim\mathcal{X}_p(t)\mathcal{Y}_p(t),
    \end{aligned}
\end{equation*}
\begin{equation*}
    \begin{aligned}
       \|\tau^{\alpha}(n^+)^{h}(u^+)^{\ell} \|^h_{\widetilde{L}^\infty_T(\dot{B}_{p,1}^{\frac{d}{p}-1})} \lesssim \|\tau^{\alpha}(n^+)^{h}\|_{\widetilde{L}^\infty_t(\dot{B}_{p,1}^{\frac{d}{p}})}\|(u^+)^{\ell}\|_{\widetilde{L}^\infty_t(\dot{B}_{p,1}^{\frac{d}{p}-1})}\lesssim\mathcal{X}_p(t)\mathcal{Y}_p(t).
    \end{aligned}
\end{equation*}
Bernstein inequality and embedding imply that
\begin{equation*}
    \begin{aligned}
        \|\tau^{\alpha}(n^+)^{\ell}(u^+)^{\ell}\|^h_{\widetilde{L}^\infty_t(\dot{B}_{p,1}^{\frac{d}{p}-1})} \lesssim \|\tau^{\alpha}(n^+)^{\ell}(u^+)^{\ell}\|_{\widetilde{L}^\infty_t(\dot{B}_{2,1}^{\frac{d}{2}+1})} 
    \end{aligned}
\end{equation*}
and thus it follows from Lemma \ref{Lemma1.2} ( $\dot{B}_{2,1}^{\frac{d}{2} } \hookrightarrow L^{\infty}$) that
\begin{equation*}
    \begin{aligned}
        \|\tau^{\alpha}(n^+)^{\ell}(u^+)^{\ell}\|_{\widetilde{L}^\infty_t(\dot{B}_{2,1}^{\frac{d}{2}+1})} \lesssim &\|\tau^{\frac{d}{4}-\frac{\sigma_1}{2}}(n^+)^{\ell}\|_{L^\infty_t(\dot{B}_{2,1}^{\frac{d}{2} })}\|\tau^{\frac{d}{4}-\frac{\sigma_1}{2}+\frac{1}{2}-\epsilon}(u^+)^{\ell}\|_{\widetilde{L}^\infty_t(\dot{B}_{p,1}^{\frac{d}{2}+1})}\\
        &+\|\tau^{\frac{d}{4}-\frac{\sigma_1}{2}}(u^+)^{\ell}\|_{L^\infty_t(\dot{B}_{2,1}^{\frac{d}{2} })}\|\tau^{\frac{d}{4}-\frac{\sigma_1}{2}+\frac{1}{2}-\epsilon}(n^+)^{\ell}\|_{\widetilde{L}^\infty_t(\dot{B}_{p,1}^{\frac{d}{2}+1})}.
    \end{aligned}
\end{equation*}
Due to the definition of $\mathcal{Y}_p(t)$, we have 
\begin{equation*}
    \begin{aligned}
        \|\tau^{\frac{d}{4}-\frac{\sigma_1}{2}+\frac{1}{2}-\epsilon}((n^+)^{\ell},(u^+)^{\ell})\|_{\widetilde{L}^\infty_t(\dot{B}_{p,1}^{\frac{d}{2}+1})}\lesssim &\|\tau^{\frac{d}{4}-\frac{\sigma_1}{2}+\frac{1}{2}-\epsilon}((n^+)^{\ell},(u^+)^{\ell})\|_{L^\infty_T(\dot{B}_{p,1}^{\frac{d}{2}+1-2\epsilon})}
        \lesssim \mathcal{Y}_p(t).
    \end{aligned}
\end{equation*}
Therefore, we have 
\begin{equation*}
    \begin{aligned}
        \|\tau^{\alpha}(n^+)^{\ell}(u^+)^{\ell}\|^h_{\widetilde{L}^\infty_t(\dot{B}_{p,1}^{\frac{d}{p}-1})}\lesssim \mathcal{Y}^2(t).
    \end{aligned}
\end{equation*}

We claim that 
\begin{equation}
    \begin{aligned}
        \|\tau^{\alpha}(\mathcal{F}_2^+,\mathcal{F}_2^-)\|^h_{\widetilde{L}^\infty_t(\dot{B}_{p,1}^{\frac{d}{p}-1})}\lesssim  \mathcal{X}_p(t)\mathcal{Y}_p(t)+ \mathcal{Y}_p^2(t) .
    \end{aligned}
\end{equation}
By the definition of $\mathcal{F}_2^+$ and  $ \mathcal{F}_2^-$, we just only deal  with $\mathcal{F}_2^+$. Recall that $\mathcal{F}_2^+\sim g^+(n^+,n^-)\nabla n^+ +(u^+\cdot \nabla )u^++ h^+(n^+,n^-)\nabla n^+\nabla u^+ +l^+(n^+,n^-)\Delta u^+$. Regarding the convection term $(u^+\cdot \nabla )u^+$, we write it as 
$$(u^+\cdot \nabla )u^+=u^+\cdot \nabla (u^+)^h+(u^+ )^h\cdot \nabla(u^+)^{\ell}+(u^+ )^\ell\cdot \nabla(u^+)^{\ell}.$$
Then we get
\begin{equation*}
    \begin{aligned}
         \|\tau^{\alpha}u^+\cdot \nabla (u^+)^h\|^h_{\widetilde{L}^\infty_t(\dot{B}_{p,1}^{\frac{d}{p}-1})}\lesssim \| u^+ \|_{\widetilde{L}^\infty_t(\dot{B}_{p,1}^{\frac{d}{p}-1})}\| \tau^{\alpha}\nabla(u^+)^h\|_{\widetilde{L}^\infty_t(\dot{B}_{p,1}^{\frac{d}{p} })}\lesssim\mathcal{X}_p(t)\mathcal{Y}_p(t),
    \end{aligned}
\end{equation*}
\begin{equation*}
    \begin{aligned}
         \|\tau^{\alpha}(u^+)^h\cdot \nabla (u^+)^\ell\|^h_{\widetilde{L}^\infty_t(\dot{B}_{p,1}^{\frac{d}{p}-1})}&\lesssim \| \tau^{\frac{d}{4}-\frac{\sigma_1}{2}-\frac{\epsilon}{2}}(u^+)^h\|_{\widetilde{L}^\infty_t(\dot{B}_{p,1}^{\frac{d}{p}-1})}\| \tau^{\frac{d}{4}-\frac{\sigma_1}{2}+\frac{1}{2}-\frac{\epsilon}{2}}\nabla(u^+)^\ell\|_{\widetilde{L}^\infty_t(\dot{B}_{p,1}^{\frac{d}{p} })}\\
         &\lesssim \mathcal{Y}_p^2(t)
    \end{aligned}
\end{equation*}
and
\begin{equation*}
    \begin{aligned}
        \|\tau^{\alpha}(u^+ )^\ell\cdot \nabla(u^+)^{\ell}\|^h_{\widetilde{L}^\infty_t(\dot{B}_{p,1}^{\frac{d}{p}-1})}&\lesssim \| \tau^{\frac{d}{4}-\frac{\sigma_1}{2}-\frac{\epsilon}{2}}(u^+)^\ell\|_{\widetilde{L}^\infty_t(\dot{B}_{p,1}^{\frac{d}{p} })}\| \tau^{\frac{d}{4}-\frac{\sigma_1}{2}+\frac{1}{2}-\frac{\epsilon}{2}}\nabla(u^+)^\ell\|_{\widetilde{L}^\infty_t(\dot{B}_{p,1}^{\frac{d}{p} })}\\
        &\lesssim \mathcal{Y}_p^2(t).
    \end{aligned}
\end{equation*}
To bound the term containing $g^+(n^+,n^-)\nabla n^+$, we have the following 
\begin{equation*}
    \begin{aligned}
        \|\tau^{\alpha}g^+(n^+,n^-)\nabla (n^+)^{\ell}\|^h_{\widetilde{L}^\infty_t(\dot{B}_{p,1}^{\frac{d}{p}-1})}&\lesssim\| \tau^{\frac{d}{4}-\frac{\sigma_1}{2}-\frac{\epsilon}{2}}(n^+,n^-)\|_{\widetilde{L}^\infty_t(\dot{B}_{p,1}^{\frac{d}{p} })}\| \tau^{\frac{d}{4}-\frac{\sigma_1}{2}+\frac{1}{2}-\frac{\epsilon}{2}}\nabla(n^+)^\ell\|_{\widetilde{L}^\infty_t(\dot{B}_{p,1}^{\frac{d}{p} })}\\
        &\lesssim\mathcal{Y}_p^2(t),
    \end{aligned}
\end{equation*}
\begin{equation*}
    \begin{aligned}
        \|\tau^{\alpha}g^+(n^+,n^-)\nabla (n^+)^{h}\|^h_{\widetilde{L}^\infty_t(\dot{B}_{p,1}^{\frac{d}{p}-1})}&\lesssim\|  (n^+,n^-)\|_{\widetilde{L}^\infty_t(\dot{B}_{p,1}^{\frac{d}{p} })}\| \tau^{\alpha}\nabla(n^+)^h\|_{\widetilde{L}^\infty_t(\dot{B}_{p,1}^{\frac{d}{p}-1 })}\\
        &\lesssim\mathcal{Y}_p^2(t).
    \end{aligned}
\end{equation*}
The term $h^+(n^+,n^-)\nabla n^+\nabla u^+$ and $l^+(n^+,n^-)\Delta u^+$ can be treated similarily. So we have
\begin{equation*}
    \begin{aligned}
        \|\tau^{\alpha}h^+(n^+,n^-)\nabla n^+\nabla u^+\|^h_{\widetilde{L}^\infty_t(\dot{B}_{p,1}^{\frac{d}{p}-1})}+& \|\tau^{\alpha}l^+(n^+,n^-)\Delta u^+\|^h_{\widetilde{L}^\infty_t(\dot{B}_{p,1}^{\frac{d}{p}-1})}\\
        \lesssim& \mathcal{X}_p (t)\mathcal{Y}_p^2(t)+\mathcal{X}_p^2 (t)\mathcal{Y}_p (t).
    \end{aligned}
\end{equation*}
With the same argument for $ \mathcal{F}_2^-$.  Finally, we obtain
\begin{equation}
    \begin{aligned}
        \|\tau^{\alpha}(\mathcal{F}_2^+,\mathcal{F}_2^-)\|^h_{\widetilde{L}^\infty_t(\dot{B}_{p,1}^{\frac{d}{p}-1})}\lesssim\mathcal{X}_p(t)\mathcal{Y}_p(t) +\mathcal{Y}_p^2(t).
    \end{aligned}
\end{equation}
Similarly, we have 
\begin{equation}
    \begin{aligned}
        \sum_{j\geq j_0-1}2^{j(\frac{d}{p}-1)}\sup_{1\leq t\leq T}\tau^{\alpha} G^{\pm}_3(t)\lesssim \mathcal{X}_p(t)\mathcal{Y}_p(t)+\mathcal{Y}_p^2(t).
    \end{aligned}
\end{equation}
For commutator terms, by Lemma \ref{lemcommutator}, we have 
\begin{equation}
    \begin{aligned}
        \sum_{j\geq j_0-1}2^{j(\frac{d}{p}-1)}\sup_{1\leq t\leq T}\tau^{\alpha} \|R^{k,\pm}_j\|_{L^p}\lesssim\mathcal{X}_p(t)\mathcal{Y}_p(t)+\mathcal{Y}_p^2(t).
    \end{aligned}
\end{equation}

Putting all the above estimates together, we thus obtain
\begin{equation}
    \begin{aligned}
        \|\langle \tau\rangle^{\alpha}(\nabla n^+,\nabla n^-,u^+,u^-)\|^h_{\widetilde{L}^{\infty}_t(\dot{B}_{p,1}^{\frac{d}{p}-1})}&\lesssim \|(\nabla n_0^+,\nabla n_0^-,u_0^+,u_0^-)\|^h_{\dot{B}_{p,1}^{\frac{d}{p}-1}}+\mathcal{X}_p(t)\mathcal{Y}_p(t) +  \mathcal{Y}_p^2(t) .
\end{aligned}
\end{equation}
A new observation is involved in the proof, which enables us to obtain more decay for the improved regularity of velocity.
Precisely, it follows from the velocity
equations in \eqref{equ:fornu-pm} that
\begin{equation} \label{sys:u+u-}
\left\{\begin{array}{l}
\partial_{t} u^{+}-\nu_{1}^{+} \Delta u^{+}-
\nu_{2}^{+} \nabla \operatorname{div} u^{+}=\mathcal{F}_2^+-\gamma_1 \nabla n^{+}-\gamma_2 \nabla n^{-}, \\
\partial_{t} u^{-}-\nu_{1}^{-} \Delta u^{-}-
\nu_{2}^{-} \nabla \operatorname{div} u^{-}=\mathcal{F}_2^--\gamma_3 \nabla n^{+}-\gamma_4 \nabla n^{-}.
\end{array}\right.
\end{equation}
To obtain desired estimates, we multiply \eqref{sys:u+u-}  by $t^{\alpha}$ to  get
\begin{equation} \label{sys:u+u-t}
\left\{\begin{array}{l}
\partial_{t} (t^\alpha u^{+})-\nu_{1}^{+} \Delta (t^\alpha u^{+})-
\nu_{2}^{+} \nabla \operatorname{div} (t^\alpha u^{+})=t^\alpha \mathcal{F}_2^+-t^\alpha(\gamma_1 \nabla n^{+}+\gamma_2 \nabla n^{-})+\alpha t^{\alpha-1}u^+ , \\
\partial_{t}(t^\alpha u^{-})-\nu_{1}^{-} \Delta (t^\alpha u^{-})-
\nu_{2}^{-} \nabla \operatorname{div} (t^\alpha u^{-})=t^\alpha \mathcal{F}_2^--t^\alpha(\gamma_3 \nabla n^{+}+\gamma_4 \nabla n^{-})+\alpha t^{\alpha-1}u^- ,\\
t^\alpha(u^+,u^-)|_{t=0}=0.
\end{array}\right.
\end{equation}
Employing Proposition \ref{lame-system-estimate} for \eqref{sys:u+u-t}, we obtain
\begin{equation}
    \begin{aligned}
        \| \tau^{\alpha}\nabla(u^+,u^-)\|_{\widetilde{L}_t^{\infty}(\dot{B}_{p,1}^{\frac{d}{p}})}^{h}&\lesssim \|\tau^\alpha (\mathcal{F}_2^+,\mathcal{F}_2^-)\|^h_{\widetilde{L}_t^{\infty}(\dot{B}_{p,1}^{\frac{d}{p}-1})}+\|\tau^{\alpha-1}(u^+,u^-)\|^h_{\widetilde{L}_t^{\infty}(\dot{B}_{p,1}^{\frac{d}{p}-1})}\\
        &\quad+\|\tau^\alpha\nabla(n^+,n^-)\|^h_{\widetilde{L}_t^{\infty}(\dot{B}_{p,1}^{\frac{d}{p}-1})}.
    \end{aligned}
\end{equation}
We notice that 

\begin{equation*}
\left\{
    \begin{aligned}
       &\|\tau^\alpha\nabla(n^+,n^-)\|^h_{\widetilde{L}_t^{\infty}(\dot{B}_{p,1}^{\frac{d}{p}-1})}\lesssim\|\langle\tau\rangle^\alpha\nabla(n^+,n^-)\|^h_{\widetilde{L}_t^{\infty}(\dot{B}_{p,1}^{\frac{d}{p}-1})},\\
       &\|\tau^{\alpha-1}(u^+,u^-)\|^h_{\widetilde{L}_t^{\infty}(\dot{B}_{p,1}^{\frac{d}{p}-1})}\lesssim\|\langle\tau\rangle^{\alpha}(u^+,u^-)\|^h_{\widetilde{L}_t^{\infty}(\dot{B}_{p,1}^{\frac{d}{p}-1})},
    \end{aligned}
    \right.
\end{equation*}
and 
\begin{equation}
    \begin{aligned}
        \|\tau^\alpha (\mathcal{F}_2^+,\mathcal{F}_2^-)\|^h_{\widetilde{L}^{\infty}(\dot{B}_{p,1}^{\frac{d}{p}-1})}\lesssim \mathcal{X}_p(t)\mathcal{Y}_p(t)+\mathcal{Y}_p^2(t).
    \end{aligned}
\end{equation}
Thus, we have 
\begin{equation}
    \begin{aligned}
         \| \tau^{\alpha}\nabla(u^+,u^-)\|_{\widetilde{L}^{\infty}(\dot{B}_{p,1}^{\frac{d}{p}})}^{h}&\lesssim\|(\nabla n_0^+,\nabla n_0^-,u_0^+,u_0^-)\|^h_{\dot{B}_{p,1}^{\frac{d}{p}-1}} +\mathcal{X}_p(t)\mathcal{Y}_p(t)+ \mathcal{Y}_p^2(t).
    \end{aligned}
\end{equation}\end{proof}

Finally, adding up \eqref{eq:low-frequency} and \eqref{eq:high-frequency} yields for all $t\geq 0,$
\begin{equation}
    \begin{aligned}
    \mathcal{Y}_p(t)   \lesssim\, \mathcal{Y}_p(0)  +\|(\nabla n_0^+,\nabla n_0^-,u_0^+,u_0^-)\|^h_{\dot{B}_{p,1}^{\frac{d}{p}-1}}+\mathcal{X}_p^2(t)+\mathcal{Y}_p^2(t). 
    \end{aligned}
\end{equation}
The global existence result in Theorem \ref{Thm1} ensures that $\mathcal{X}_p(t)\lesssim\, \mathcal{X}_p(0)\ll 1$  and   as 
\begin{align*}
    \|(n_0^+,n_0^-,u_0^+,u_0^-)\|_{\dot{B}_{2,1}^{\frac{d}{2}-1}}\lesssim\, \mathcal{Y}_p(0)\ll1,
\end{align*}
one can conclude that \eqref{ineqtime} is satisfied for all time if $\mathcal{Y}_p(0)$ and $\mathcal{X}_p(0)$ are small enough. This completes the proof of Theorem \ref{Thm2}.

\appendix
\titleformat{\section}
  {\normalfont\Large\bfseries}
  {Appendix~\thesection} % 这里把"Appendix A"作为编号前缀输出
  {1em}{}

 \section[Appendix~\thesection: Some estimates in the hybrid Besov space]
        {Some estimates in the hybrid Besov space}
 
We give some useful estimates, which can be found in \cite{CPAM2010}.
\begin{proposition}\label{propA1}
Let $s_1,s_2,t_1,t_2,\sigma,\tau\in\mathbb{R}$, $2\le p\le 4$ and $1\le r,r_1,r_2,r_3,r_4\le\infty$ with
\[
\frac{1}{r}=\frac{1}{r_1}+\frac{1}{r_2}=\frac{1}{r_3}+\frac{1}{r_4}.
\]
Then we have the following:

\begin{itemize}
  \item If $\sigma,\tau\le d/p$ and $\sigma+\tau>0$, then
  \begin{equation}\label{A1}
  \sum_{2^k>R_0} 2^{k(\sigma+\tau-d/p)} \|\dot\Delta_k(fg)\|_{L_T^r(L^p)}
  \lesssim\|f\|_{\tilde L_T^{r_1}(\dot{\mathcal{B}}_{2,p}^{d/2-d/p+\sigma,\,\sigma})} 
       \|g\|_{\tilde L_T^{r_2}(\dot{\mathcal{B}}_{2,p}^{d/2-d/p+\tau,\, \tau})}.
  \end{equation}

  \item If $s_1,s_2\leq d/p$ and $s_1+t_1>d-2d/p$ with $s_1+t_1=s_2+t_2$ and $\gamma\in\mathbb{R}$, then
  \begin{equation}
  \begin{aligned}\label{A2}
  &\sum_{2^k\le R_0} 2^{k(s_1+t_1-d/2)}\|\dot\Delta_k(fg)\|_{L_T^r( L^2)}\\
  \lesssim&\,   
     \|f\|_{\tilde L_T^{r_1}(\dot{\mathcal{B}}_{2,p}^{s_1,\,s_1-d/2+d/p})}\|g\|_{\tilde L_T^{r_2}(\dot{\mathcal{B}}_{2,p}^{t_1,\,t_1-d/2+d/p+\gamma})}
     +\|g\|_{\tilde L_T^{r_3}(\dot{\mathcal{B}}_{2,p}^{s_2,\,s_2-d/2+d/p})}\|f\|_{\tilde L_T^{r_4}(\dot{\mathcal{B}}_{2,p}^{t_2,\,t_2-d/2+d/p+\gamma})}
     .
     \end{aligned}
  \end{equation}

  \item If $s_1,s_2\leq d/2$ and $s_1+t_1>d/2-d/p$ with $s_1+t_1=s_2+t_2$, then
  \begin{equation}
  \begin{aligned}\label{A3}
  &\sum_{k\in\mathbb{Z}} 2^{k(s_1+t_1-d/2)}\|\dot\Delta_k(fg)\|_{L_T^r( L^2)}\\
  \lesssim&\,  
     \|f\|_{\tilde L_T^{r_1}(\dot{\mathcal{B}}_{2,p}^{s_1,\,s_1-d/2+d/p})}\|g\|_{\tilde L_T^{r_2}(\dot B_{2,1}^{t_1 })}
     +\|g\|_{\tilde L_T^{r_3}(\dot{\mathcal{B}}_{2,p}^{s_2,\,s_2-d/2+d/p})}\|f\|_{\tilde L_T^{r_4}(\dot B_{2,1}^{t_2 })}.
\end{aligned}
  \end{equation}
\end{itemize}
\end{proposition}

\begin{proposition}\label{propA2}
Let $2\le p\le 4$, $s,\sigma>0$, and $s\ge\sigma-d/2+d/p$, $r\ge1$.
Assume that $F\in W_{\mathrm{loc}}^{[s]+2,\infty}\cap W_{\mathrm{loc}}^{[\sigma]+2,\infty}$ with $F(0)=0$.
Then there holds
\begin{equation}\label{A4}
\|F(f)\|_{\tilde L_T^r(\dot{\mathcal{B}}_{2,p}^{s,\,\sigma})}
\lesssim  (1+\|f\|_{L_T^\infty(\dot{\mathcal{B}}_{2,p}^{d/p,\,d/p})})^{\max([s],[\sigma])+1}
   \|f\|_{\tilde L_T^r(\dot B_{2,p}^{s,\,\sigma})}.
\tag{A.4}
\end{equation}
For any $s>0$ and $p\ge1$, there holds
\begin{equation}\label{A5}
\|F(f)\|_{\tilde L_T^r(\dot{\mathcal{B}}_{p,1}^{s})}
\lesssim (1+\|f\|_{L_T^\infty (L^\infty)})^{[s]+1}\|f\|_{\tilde L_T^r(\dot{\mathcal{B}}_{p,1}^{s})}.
\tag{A.5}
\end{equation}
\end{proposition}
\section[Appendix~\thesection: Technical lemmas in Besov spaces]
        {Technical lemmas in Besov spaces}

The following commutator estimates will used to control some nonlinearities in high frequencies.
\begin{lemma}\label{lemcommutator}
    Let $p\in [1,\infty )$ and $-\frac{d}{2}-1\leq s\leq \frac{d}{2}+1$. Then it holds that 
    \begin{align}
        \sum_{j\in \mathbb{Z}}2^{js}\|[v,\dot{\Delta}_j]\partial_i u\|_{L^2}\lesssim\|\nabla v\|_{\dot{B}^{\frac{d}{2}}_{2,1}}\|u\|_{\dot{B}^{s}_{2,1}}, \ i=1,2,\cdots,d,
    \end{align}
    for the commutator $[A,B]=AB-BA.$
\end{lemma}

We state the following result about the continuity for composition functions.
\begin{lemma}\emph{\cite{SWYZ}}\label{Lemma5-7}
    Let $G:I\to \mathbb{R}$ be a  smooth function satisfying $G(0)=0.$ For any $1\leq p\leq \infty,$ $s>0$, and $g\in \dot{B}^{s}_{2,1}\cap L^{\infty}$, there holds that $G(g)\in \dot{B}^{s}_{p,r}\cap L^{\infty}$ and 
    \begin{align}
        \|G(g)\|_{\dot{B}^{s}_{p,r}}\leq C_g\|g\|_{\dot{B}^{s}_{p,r}},\label{F0}
    \end{align}
where the constant $C_g>0$ depends only on $\|g\|_{L^{\infty}}$, $G'$, $s, $and $d$.

In addition, if $g_{1}, g_{2}\in \dot{B}^{s}_{p,1}\cap L^{\infty}$, then it holds
\begin{align}
&\|G(g_{1})-G(g_{2})\|_{\dot{B}^{s}_{p,1}}\leq C_{g_{1},g_{2}}(1+\|(g_{1},g_{2})\|_{\dot{B}^{\frac{d}{2}}_{2,1}})\|g_{1}-g_{2}\|_{\dot{B}^{s}_{p,1}},\quad ~s\in (-\frac{d}{2},\frac{d}{2}],\label{F2}
\end{align}
where the constant $C_{g_{1},g_{2}}>0$ depends only on $\|(g_{1},g_{2})\|_{L^{\infty}}$, $G'$, $s$, $p$ and $d$.
\end{lemma}

We prove the following lemma about the continuity for composition of multi-component functions. It should be noted that \eqref{F11} will be used to deal with the two-dimensional case in $\dot{B}^{0}$.
\begin{lemma}\emph{\cite{SWYZ}}\label{lemma64}
Let $m\in \mathbb{N}$, $s>0$, and $G\in C^{\infty}(\mathbb{R}^{m})$ satisfy $G(0,...,0)=0$. Then for any $f_{i}\in\dot{B}^{s}\cap L^{\infty}$ $(i=1,...,m)$, there exists a constant $C_{f}>0$ depending on $\sum_{i=1}^{m}\|f_{i}\|_{L^{\infty}}$, $F$, $s$, $m$ and $d$ such that
\begin{equation}
\begin{aligned}
\|G(f_{1},...,f_{m})\|_{\dot{B}^{s} }\leq C_{f}\sum_{i=1}^{m}\|f_{i}\|_{\dot{B}^{s} }.\label{F1}
\end{aligned}
\end{equation}
In the case $s>-\frac{d}{2}$ and $f_{i}\in\dot{B}^{s}\cap \dot{B}^{\frac{d}{2}}$, it holds that
\begin{equation}
\begin{aligned}
\|G(f_{1},...,f_{m})\|_{\dot{B}^{s} }\leq C_{f}\big(1+\sum_{i=1}^{m}\|f_{i}\|_{\dot{B}^{\frac{d}{2}}_{2,1}} \big)\sum_{i=1}^{m}\|f_{i}\|_{\dot{B}^{s} }.\label{F11}
\end{aligned}
\end{equation}
Furthermore, for any $f^{1}_{i}, f^{2}_{i}\in\dot{B}^{s}\cap \dot{B}^{\frac{d}{2}}$, we have
\begin{equation}\label{F3}
\begin{aligned}
&\|G(f^{1}_{1},...,f^{1}_{m})-G(f^{1}_{1},...,f^{1}_{m})\|_{\dot{B}^{s} }\leq C^{*}_{f}\big(1+\sum_{i=1}^{m}\|f_{i}\|_{\dot{B}^{s}\cap\dot{B}^{\frac{d}{2}}_{2,1}} \big)\sum_{i=1}^{m}\|f_{i}^{1}-f_{i}^{2}\|_{\dot{B}^{s}\cap\dot{B}^{\frac{d}{2}}_{2,1}}.
\end{aligned}
\end{equation}
Here $C^{*}_{f}>0$ depends on $\sum_{i=1}^{m}\|(f^{1}_{i},f^{2}_{i})\|_{L^{\infty}}$, $F$, $s$, $m$ and $d$.

\end{lemma}

\begin{proposition}\emph{\cite{Danchin2005}}\label{prop:transport_est}
Let 
$s\in\Bigl(-d\min\Bigl(\frac1p,\frac1{p'}\Bigr),\,1+\frac dp\Bigr)$ and 
$1\le p,q\le +\infty,$ 
where $p'$ is the conjugate exponent of  $p$.
Let $v$ be a vector field such that
$\nabla v\in L_T^1\bigl(\dot B^{d/p}_{p,1}\bigr). $
Assume that \(f_0\in \dot B^s_{p,q}\), \(g\in L_T^1\bigl(\dot B^s_{p,q}\bigr)\), and that \(f\) is a solution of the transport equation
\[
\partial_t f+ v\cdot\nabla f = g,
\qquad
f(0,x)=f_0.
\]
Then for \(t\in[0,T]\), there holds
\[
\|f\|_{\widetilde L_t^\infty(\dot B^s_{p,q})}
\le
\exp\!\left(
C\int_0^t \|\nabla v(\tau)\|_{\dot B^{d/p}_{p,1}}\,d\tau
\right)
\left(
\|f_0\|_{\dot B^s_{p,q}}
+\int_0^t \|g(\tau)\|_{\dot B^s_{p,q}}\,d\tau
\right).
\]
\end{proposition}

Finally, we need the optimal regularity estimates for the Lam\'e system below.
\begin{proposition}\emph{\cite{SWYZ}}\label{lame-system-estimate}
Let $T>0,$ $\mu>0$, $2\mu+\lambda>0$, $s\in \mathbb{R},$ $1\leq p,r\leq \infty,$ and $1\leq \varrho_2\leq\varrho_1\leq\infty.$ Assume that $u_0\in \dot{B}^s_{p,r}$ and $f\in \widetilde{L}^{\varrho_2}(0,T;\dot{B}^{s-2+\frac{2}{\varrho_2}}_{p,r})$ hold. If $u$ is a solution of 
$$\begin{cases}
    \partial_t u-\mu\Delta u-(\mu+\lambda)\nabla\operatorname{div}u=f, \quad &x\in\mathbb{R}^d,\quad t\in(0,T), \\
   u(x,0)=u_0(x),\quad &x\in\mathbb{R}^d,
\end{cases}
$$
then $u$ satisfies
\begin{align}
    \min\{\mu,2\mu+\lambda\}^{\frac{1}{\varrho_1}}\|u\|_{\widetilde{L}_T^{\varrho_1}(\dot{B}^{s +\frac{2}{\varrho_1}}_{p,r})}\lesssim\|u_0\|_{\dot{B}^{s}_{p,r}}+\|f\|_{\widetilde{L}_T^{\varrho_2}(\dot{B}^{s-2 +\frac{2}{\varrho_2}}_{p,r})}.
\end{align}
\end{proposition}
%%%%%%%%%%%%%%%%%%%%%%%%%%%%%%%%%
\noindent \textbf{Acknowledgments} 
The authors would like to thank Professor  Ling-Yun Shou for his helpful discussions and suggestions. This work was partially supported by the Zhejiang Provincial Natural Science Foundation of China (Grant No. LZ25A010003) and the National Natural Science Foundation of  China (Grant No. 12301429).  J. Wu is supported by   the National Natural Science Foundation of China  $\#$12671258.

\noindent \textbf{Conflict of interest.} The authors do not have any possible conflicts of interest.
	
\vspace{2mm}
	
\noindent \textbf{Data availability statement.} Data sharing is not applicable to this article as no data sets were generated or analyzed during the current study.

\vspace{0.5cm}
    
%\bigbreak\bigbreak

%\newpage 

(Zhiwen He)\par\nopagebreak
\noindent\textsc{School of Mathematics,  South China Normal University, Guangzhou 510000, P.R. China}

Email address: {\texttt{zhiwen$\_$he@zju.edu.cn}}

\vspace{3ex}

(Shiyu Li)\par\nopagebreak
\noindent\textsc{Department of Mathematical Sciences, School of Science, Zhejiang Sci-Tech University, Hangzhou 310018, P. R. China}

Email address: {\texttt{lishiyu0110@zstu.edu.cn}}

\vspace{3ex}

(Jiayan Wu)\par\nopagebreak
\noindent\textsc{School of Mathematics, South China University of Technology, Guangzhou 510641, P.R. China}

Email address: {\texttt{wujiayan@scut.edu.cn}}

\vspace{3ex}

(Ting Zhang)\par\nopagebreak
\noindent\textsc{School of Mathematics, Zhejiang University, Hangzhou 310058, P.R. China}

Email address: {\texttt{zhangting79$@$zju.edu.cn}}

\vspace{3ex}

\end{document}